\documentclass[12pt]{article}
\usepackage[psamsfonts]{amssymb}
\usepackage[dvips]{graphicx}

\newtheorem{theorem}{Theorem}[section]
\newtheorem{lemma}[theorem]{Lemma}
\newtheorem{corollary}[theorem]{Corollary}
\newtheorem{conjecture}[theorem]{Conjecture}
\newtheorem{proposition}[theorem]{Proposition}
\newtheorem{definition}[theorem]{Definition}

\newcommand{\pl}{{+}}
\newcommand{\oUJ}{{ U_J}}
\newcommand{\oWJ}{{ W_J}}
\newcommand{\oU}{{ U}}
\newcommand{\oW}{{ W}}
\newcommand{\oS}{{\overline S}}
\newcommand{\nill}{{\emptyset}}

\newcommand{\One}{{\mathbf 1}}

\newcommand{\cLa} {{\Lambda}}

\newcommand{\om}{\omega}
\newcommand{\ga}{r}
\newcommand{\la}{\lambda}
\newcommand{\LD}{{\bf L}^2({{\mathbb{R}}}^d)}
\newcommand{\LU}{{\bf L}^1({{\mathbb{R}}}^d)}
\newcommand{\Ld}{{\bf L}^2}

\newcommand{\ld}{{\bf l}^2 }
\newcommand{\lu}{{\bf l}^1 }

\newcommand{\R}{{\mathbb{R}}}
\newcommand{\C}{{\mathbb{C}}}

\newcommand{\cP}{{\mathcal P}}
\newcommand{\cPinf}{{\overline{\mathcal P}_\infty}}
\newcommand{\cLainf}{{{\Lambda}_\infty^\infty}}
\newcommand{\oP}{{\overline{\mathcal P}}}
\newcommand{\tP}{{\widetilde{\mathcal P}}}

\newcommand{\Z}{{\mathbb{Z}}}

\newcommand{\N}{{\mathbb{N}}}

\newcommand{\bCd}{{\mathbf C^2}({{\mathbb{R}}}^d)}

\newcommand{\La}{{\Lambda}}

\newcommand{\cH}{{\mathcal H}}
\newcommand{\rH}{{H}}

\newcommand{\G}{{G}}

\newcommand{\Exp}{{\rm e}}
\newcommand{\Id}{{\bf 1}}
\newcommand{\iint}{{\int\!\!\int}}

\date{\ }
\title{{Group Invariant Scattering}}
\author {St\'ephane Mallat\\
{\it CMAP, Ecole Polytechnique, Palaiseau, France}\\
{\it IHES, Bures-sur-Yvette, France}
\footnote{This work was supported by the ANR-10-BLAN-0126 grant.}}

\begin{document}

\maketitle

\begin{abstract}
This paper constructs translation invariant 
operators on $\LD$, which
are Lipschitz continuous to the action of diffeomorphisms.
A scattering propagator is a path ordered product of
non-linear and non-commuting operators, each of which computes the modulus of
a wavelet transform. A local integration defines a
windowed scattering transform, which is proved to be
Lipschitz continuous to the action of $\bf C^2$ diffeomorphisms. As the window 
size increases, it converges to a wavelet
scattering transform which is translation invariant.
Scattering coefficients also provide
representations of stationary processes. Expected values
depend upon high order 
moments and can discriminate processes having the same
power spectrum. Scattering operators are extended on $\Ld (G)$, where
$G$ is a compact Lie group, and are invariant under the action of $G$.
Combining a scattering on $\LD$ and on $\Ld (SO(d))$ defines a
translation and rotation invariant scattering on $\LD$.
\end{abstract}

\section{Introduction}

Symmetry and invariants, which play a major role in physics
\cite{frankel}, are 
making their way into signal information processing.
The information content of sounds or images is typically not affected under
the action of finite groups such as translations or rotations,
and it is stable to the action of small 
diffeomorphisms that deform signals \cite{Trouve}. 
This motivates the study of translation-invariant representations
of $\LD$ functions,
which are Lipschitz continuous to the action of diffeomorphisms, 
and which keep high-frequency
information to discriminate different types signals.
Invariance to the action of compact Lie groups and
rotations are then studied.

We first concentrate on translation invariance. 
Let $L_c f(x) = f(x-c)$ denote the translation of 
$f \in \LD$ by $c \in \R^d$.
An operator $\Phi$ from 
$\LD$ to a Hilbert space $\cH$ is 
{\it translation-invariant} if
$\Phi(L_c f) = \Phi(f)$ for all $f \in \LD$ and $c \in \R^d$.
Canonical translation invariant
operators satisfy $\Phi(f) = L_a f$ for some $a \in \R^d$ which 
depends upon $f$ \cite{Olver}. 
The modulus of the Fourier transform of $f$ is an example of non-canonical
translation invariant operator. 
However, these translation invariant operators are not Lipschitz
continuous to the action of diffeomorphisms.
Instabilities to deformations are well-known to appear at high frequencies
\cite{Hormander}. The major difficulty is to maintain the Lipschitz
continuity over high frequencies.

To preserve stability in $\LD$ we want 
$\Phi$ to be nonexpansive:
\[
\forall (f,h) \in \LD^2\,\,,\,\,\|\Phi(f) - \Phi(h)\|_\cH \leq \|f-h\|. 
\]
It is then sufficient to verify its Lipschitz
continuity relatively to the action of 
small diffeomorphisms close to translations.
Such a diffeomorphism
transforms $x \in \R^d$ into $x - \tau (x)$, where
$\tau (x) \in \R^d$ is the displacement field.
Let $L_\tau f(x) = f(x-\tau(x))$ denote the action of the diffeomorphism
$\One - \tau$ on $f$. 
Lipschitz stability means that
$\|\Phi(f) - \Phi (L_\tau f) \|$ is bounded by the ``size'' of the 
diffeomorphism and hence by the distance
between the $\One - \tau$ and $\One$, up to a multiplicative constant
multiplied by $\|f\|$. 
Let $|\tau(x)|$ denote the Euclidean norm in $\R^d$,
$|\nabla \tau(x)|$ the sup norm of the matrix $\nabla \tau(x)$, and
$|\rH \tau (x)|$ the sup norm of the Hessian tensor. 
The weak topology on $\bf C^2$ diffeomorphisms defines a distance
between $\One - \tau$ and $\One$, over any compact subset $\Omega$ of $\R^d$,
by:
\begin{equation}
\label{weaksons}
d_\Omega (\One,\One-\tau) = \sup_{x \in \Omega} | \tau (x)| + 
\sup_{x \in \Omega} |\nabla \tau (x)| +
\sup_{x \in \Omega} |\rH \tau (x)| .
\end{equation}
A translation invariant operator $\Phi$ is said to be 
{\it Lipschitz continuous} to the action of ${\bf C}^2$ diffeomorphisms
if for any compact $\Omega \subset \R^d$
there exists $C$ such that for all $f\in \LD$ supported in $\Omega$
and all $\tau \in \bCd$
\begin{equation}
\label{invafns2}
\|\Phi(f) - \Phi(L_{\tau} f) \|_\cH \leq C\,\|f\|\, 
\Bigl(\sup_{x \in \R^d} |\nabla \tau(x)| +
\sup_{x \in \R^d} |\rH \tau(x)| \Bigr)~.
\end{equation}
Since $\Phi$ is translation invariant,
the Lipschitz upper bound does not depend upon the maximum translation 
amplitude  $\sup_x |\tau(x)|$ of the diffeomorphism metric (\ref{weaksons}).
The Lipschitz continuity (\ref{invafns2})
implies that $\Phi$ is invariant to global translations, 
but it is much stronger.
$\Phi$ is almost invariant to ``local translations'' by $\tau (x)$, 
up to the first and second order deformation terms.

High frequency instabilities to deformations can be avoided by
grouping frequencies into 
dyadic packets in $\R^d$, with a wavelet transform. 
However, a wavelet transform is not translation invariant. A translation
invariant operator is constructed with a scattering procedure along multiple
paths, which preserves the Lipschitz stability of wavelets to the action
of diffeomorphisms.
A scattering propagator is first defined as a path ordered product of
non-linear and non-commuting operators, each of which computes the modulus of
a wavelet transform \cite{MallatEUSIPCO}. 
This cascade of convolutions and modulus can also be interpreted as a 
convolutional neural-network \cite{LeCun1}. 
A windowed scattering transform is a nonexpansive operator which
locally integrates the scattering propagator output. For appropriate wavelets, 
the main theorem in Section \ref{recindsnsec}
prove that a windowed scattering 
preserves the norm: $\|\Phi(f) \|_\cH = \|f\|$ for all $f \in \LD$,
and it is Lipschitz continuous to $\bf C^2$ diffeomorphisms.

When the window size increases, 
windowed scattering transforms converge to a translation
invariant scattering transform, defined on 
a path set $\cPinf$ which is not countable.
Section \ref{infdnsdfonsdsec} introduces 
a measure $\mu$ and a metric on 
$\cPinf$, and proves that scattering transforms 
of functions in $\LD$ belong to $\Ld(\cPinf,d\mu)$.
A scattering transform has striking similarities with a Fourier
transform modulus, but a different behavior relatively to the action
of diffeomorphisms. Numerical examples are shown.
An open conjecture remains on conditions for strong convergence in $\LD$.

The representation of stationary processes with the Fourier power spectrum results from the translation invariance of the Fourier modulus. Similarly, 
Section \ref{Stationary} defines an expected scattering transform
which maps stationary processes to an $\ld$ space.
Scattering coefficients depend upon high order moments of stationary
processes, and can thus discriminate processes having same
second-order moments. As opposed to the Fourier spectrum,
a scattering representation is
Lipschitz continuous to random deformations up to a log term. 
For large classes of ergodic processes, it is numerically observed that
the scattering transform of a single realization provides a
mean-square consistent estimator of the expected 
scattering transform.

Section \ref{MultiGroup} extends
scattering operators to build invariants to actions of compact Lie groups $G$.
The left action of $g \in G$ on $f \in \Ld(\G)$ is denoted
$L_g f(r) = f(g^{-1} r)$. An operator $\Phi$ on $\Ld (G)$ is 
invariant to the action of $G$ if $\Phi(L_g f) = \Phi(f)$
for all $f \in \Ld(G)$ and all $g \in G$. 
Invariant scattering operators are constructed on $\Ld(G)$
with a scattering propagator which iterates on a wavelet transform
defined on $\Ld(G)$, and a modulus operator which removes complex phases. 
A translation and rotation
invariant scattering on $\LD$ is obtained by combining a scattering
on $\LD$ and a scattering on $\Ld(SO(d))$.

A software package is available
at {\it www.cmap.polytechnique.fr/scattering},
to reproduce numerical experiments.
Applications to audio and image classification
can be found in \cite{Joakim,Bruna,Bruna2,Sifre}.

{\bf Notations:} $\|\tau\|_\infty  = \sup_{x \in \R^d} |\tau(x)|$,
$\|\Delta \tau \|_\infty =  \sup_{(x,u) \in \R^{2d}} |\tau(x) - \tau(u)|$,
$\|\nabla \tau \|_\infty =  \sup_{x \in \R^d} |\nabla \tau(x)|$ and
$\|H \tau \|_\infty =  \sup_{x \in \R^d} |H \tau(x)|$ where
$|H \tau(x)|$ is the norm of the Hessian tensor.
The inner product of $(x,y) \in \R^{2d}$ is $x \cdot y$.
The norm of $f$ in a Hilbert space is
$\|f\|$ and in $\LD$: $\|f\|^2 = \int |f(x)|^2 \,dx$.
The norm in 
${\bf L^1} (\R^d)$ is
$\|f\|_1 = \int |f(x)|\, dx$. 
We denote $\hat f(\om) = \int f(x)\, e^{-i x \cdot \om}\,d\om$ 
the Fourier transform of $f$.
We denote
$g \circ f(x) = f(g x)$ the action of a group element $g \in G$.
An operator $R$ parametrized
by $p$ is denoted $R[p]$ and 
$R[\Omega] = \{R[p]\}_{p \in \Omega}$.
The sup norm of a linear operator $A$ in $\LD$
is denoted $\|A \|$ and 
the commutator of two operators is $[A\,,\,B] = AB - BA$.

\section{Finite Path Scattering}
\label{recindsnsec}

To avoid high frequency instabilities under the action of diffeomorphisms,
Section \ref{fastalgo} introduces scattering operators, which
iteratively apply wavelet transforms and remove complex phases
with a modulus. Section \ref{upsdfnwd} proves that a scattering is
nonexpansive and preserves $\LD$ norms.
Translation invariance and Lipschitz continuity to deformations
are proved in Section \ref{invasnsf1} and \ref{invasnsf2}.
 
\subsection{From Fourier to Littlewood-Paley Wavelets}
\label{scalittwnwavel}

The Fourier transform modulus $\Phi(f) = |\hat f|$
is translation-invariant. Indeed for
$c \in \R^d$, the translation $L_c f(x) = f(x-c)$
satisfies 
$\widehat {L_c f} (\om) = e^{-i c. \om} \hat f (\om)$ and hence
$|\widehat {L_c f} | = |\hat f|$.
However, deformations lead to well-known instabilities at high frequencies
\cite{Hormander}. This is illustrated with 
a small scaling operator,
$L_\tau f (x) = f(x-\tau(x)) = f((1-s) x)$, for $\tau(x) = s x$ and 
$\|\nabla \tau\|_\infty = |s| < 1$.
If $f(x) = e^{i \xi\cdot x} \,\theta(x)$ then scaling by $1-s$
translates the central frequency $\xi$ to $(1-s) \xi$. If
$\theta$ is regular with a fast decay then
\begin{equation}
\label{instablsndf}
\||\widehat {L_\tau f} | -|\widehat {f} |\| \sim |s|\, |\xi|\,\|\theta\| =
\|\nabla \tau\|_\infty\, |\xi|\,\|f\|\,.
\end{equation}
Since $|\xi|$ can be arbitrarily large, $\Phi(f) = |\hat f|$
does not satisfy the Lipschitz continuity
condition (\ref{invafns2}) when scaling high frequencies.
The frequency displacement from
$\xi$ to $(1-s) \,\xi$ has a small impact 
if sinusoidal waves are replaced by localized functions having a Fourier support
which is wider at high frequencies. This is achieved
by a wavelet transform \cite{Frazier,Meyerbook}, whose properties
are briefly reviewed in this section.

A wavelet transform is constructed by dilating 
a wavelet $\psi \in \LD$ with a scale sequence $\{a^j \}_{j \in \Z}$ for
$a > 1$.
For image processing, usually $a = 2$ \cite{Bruna,Bruna2}.
Audio processing requires a better frequency resolution  
with typically $a \leq 2^{1/8}$ \cite{Joakim}.
To simplify notations, we normalize $a = 2$, with no loss of generality.
Dilated wavelets are also rotated with elements of a finite rotation
group $\G$, which also includes the reflection $-\One$
with respect to $0$: $-\One x = - x$.
If $d$ is even then $G$ is a subgroup of $SO(d)$, and 
if $d$ is odd then $G$ is a finite subgroup of $O(d)$.
A mother wavelet $\psi$ is dilated by $2^{-j}$ and rotated by $r \in \G$
\begin{equation}
\label{motherandfosdf}
\psi_{2^j r }(x) = 2^{dj}\,\psi(2^{j} \,\ga^{-1} x).
\end{equation}
Its Fourier transform is
$\hat \psi_{2^j r} (\om) = \hat \psi(2^{-j} r^{-1} \om)$.
A scattering transform is computed with wavelets that
can be written
\begin{equation}
\label{wavedesin}
\psi(x) = e ^{i \eta \cdot x}\,\theta(x)~~\mbox{and hence}~~
\hat \psi(\om) = \hat \theta (\om - \eta)~,
\end{equation}
where $\hat \theta(\omega)$ is a real function
concentrated in a low frequency ball
centered at $\omega = 0$, whose radius is of the order of $\pi$. 
It results that $\hat \psi(\om)$ is real and concentrated in a frequency ball
of same radius, but centered at $\om = \eta$.
To simplify notations we denote $\la = 2^j r \in 2^\Z \times G$,
with $|\la| = 2^j$. After dilation and rotation, 
$\hat \psi_{\la} (\om) = \hat \theta (\la^{-1} \om - \eta)$ covers
a ball centered at $\la \eta$ with a radius 
proportional to $|\la| = 2^j$.
The index $\la$ thus specifies the frequency localization and spread of
$\hat \psi_\la$.

As opposed to wavelet bases, a Littlewood-Paley 
wavelet transform \cite{Frazier,Meyerbook} is a redundant
representation which computes convolution values
at all $x \in \R^d$, without subsampling:
\begin{equation}
\label{wavedfn0}
\forall x \in \R^d~~,~~W[\la] f(x) = f \star \psi_{\la} (x) =
\int f(u)\ \psi_{\la} (x-u)\,du ~.
\end{equation}
Its Fourier transform is
\[
\widehat {W[\la] f} (\om) = \hat f(\om) \,\hat \psi_{\la} (\om) = 
\hat f(\om) \,\hat \psi (\la^{-1} \om) ~.
\]
If $f$ is real then $\hat f(-\om) = \hat f^*(\om)$ and if $\hat \psi(\om)$
is real then $W[-\la] f = W[\la] f^*$. 
Let $\G^+$ denote the quotient of $\G$ with $\{-\One\,,\,\One \}$,
where two rotations $r$ and $-r$ are equivalent. 
It is sufficient to compute
$W[2^j r] f$ for ``positive'' rotations $r \in G^+$. If $f$ is complex then 
$W[2^j r] f$ must be computed for all 
$r \in G = G^+ \times \{ -\One,\One\}$.

A wavelet transform at a scale $2^J$ only keeps wavelets of frequencies
$2^{j} > 2^{-J}$.
The low frequencies which are not covered by these wavelets are
provided by an averaging over a spatial domain proportional to $2^J$:
\begin{equation}
\label{averaginop}
A_J f = f \star \phi_{2^J}~~\mbox{with}~~\phi_{2^J} (x) = 2^{-dJ} \phi(2^{-J} x)~.
\end{equation}
If $f$ is real then 
the wavelet transform 
$ W_J f = \Big\{A_J f  \,,\,\Bigl(W[\la] f \Bigr)_{\la \in \cLa_J} \Bigr\}$
is indexed by $\cLa_J = \{ \la = 2^jr \,:\,r \in \G^+\,,\,2^j>2^{-J}\}$.
Its norm is 
\begin{equation}
\label{Normwinds}
\| \oWJ f \|^2 = \|A_J f \|^2 + \sum_{\la \in \cLa_J} \|W[\la] f \|^2~.
\end{equation}
If $J = \infty$ then 
$W_\infty f = \Big\{W[\la] f \}_{\la \in \cLa_\infty}$
with $\La_\infty = 2^\Z \times G^+$. Its norm is
$\| \oW_\infty f \|^2 = \sum_{\la \in \cLa_\infty} \|W[\la] f \|^2$.
For complex-valued functions $f$, all rotations in $G$ are included by
defining
$ W_J f = \Big\{A_J f  \,,\,\Bigl(W[\la] f \Bigr)_{\la \in \cLa_J
\atop -\la \in \cLa_J} \Bigr\}$ and
$W_\infty f = \Big\{W[\la] f \}_{\la \in \cLa_\infty \atop - \la \in \cLa_\infty}$.
The following proposition gives a standard Littlewood-Paley 
condition \cite{Frazier} so that $W_J$ is unitary.

\begin{proposition}
\label{Littlewood}
For any $J \in \Z$ or $J = \infty$, 
$W_J$ is unitary in the spaces of 
real-valued or complex-valued functions in $\LD$ if and only if
for almost all $\omega \in \R^d$
\begin{equation}
\label{consenas}
\beta \sum_{j=-\infty}^{\infty}  \sum_{r \in \G}
|\hat \psi (2^{-j} \ga^{-1} \om)|^2  = 1
~~\mbox{and}~~
|\hat \phi(\om)|^2 =  \beta \sum_{j=-\infty}^{0}  \sum_{r \in \G}
|\hat \psi (2^{-j} \ga^{-1} \om)|^2 ~,
\end{equation}
where $\beta = 1$ for complex functions and $\beta=1/2$ for real
functions.
\end{proposition}

{\it Proof:}
If $f$ is complex, $\beta = 1$ and one can verify that
(\ref{consenas}) is equivalent to 
\begin{equation}
\label{consenas0}
\forall J \in \Z~,~|\hat \phi(2^{J} \om)|^2 + 
\sum_{j > -J,r \in G} |\hat \psi (2^{-j} r^{-1} \om)|^2  = 1~.
\end{equation}
Since $\widehat {W[2^j\ga] f}(\om)  = \hat f(\om)\, \hat \psi_{2^j r}(\om)$,
multiplying (\ref{consenas0}) by $|\hat f (\om)|^2$, 
and applying the Plancherel formula proves that $\|W_J f \|^2 = \|f\|^2$.
For $J = \infty$ the same result is obtained by letting $J$ go to $\infty$.

Conversely, if $\|W_J f \|^2 = \|f\|^2$ then (\ref{consenas0}) is
satisfied for almost all $\omega$.
Otherwise, one can construct a function $f \neq 0$ 
where $\hat f$ has a support
in the domain of $\omega$ where (\ref{consenas0}) is not valid.
With the Plancherel formula we verify that 
$\|W_J f \|^2 \neq \|f\|^2$, which contradicts the hypothesis.

If $f$ is real then 
$|\hat f(\om)| = |\hat f(-\om)|$ so $\|W[2^j r] f \| = \|W[-2^j r] f\|$.
Hence $\|W_J f\|$ remains the same if $r$ is restricted to $G^+$ 
and $\psi$ is multiplied by $\sqrt 2$, 
which yields condition (\ref{consenas}) with $\beta = 1/2$.
$\Box$

In all the following, $\hat \psi$ is a real function which satisfies
the unitary condition (\ref{consenas}). It 
implies that $\hat \psi(0) = \int \psi (x)\,dx = 0$ and 
$|\hat \phi(r \om)| = |\hat \phi(\om)|$
for all $r \in \G$. We choose $\hat \phi(\om)$ to be real and
symmetric so that $\phi$ is also real and symmetric and 
$\phi (r x) = \phi(x)$ for all $r \in \G$. 
We also suppose that $\phi$ and 
$\psi$ are twice differentiable and that their decay as 
well as the decay of their partial derivatives of order $1$ and $2$ 
is $O((1 + |x|)^{-d-2})$.

A change of variable in the wavelet transform integral shows
that if $f$ is scaled and rotated,
$2^l g \circ f = f(2^l g x)$ with $2^l g \in 2^\Z \times \G$, then
the wavelet transform is scaled and rotated according to:
\begin{equation}
\label{wavedfn98}
W[\la] ({2^l g} \circ f) = 2^l g \circ W[{2^{-l} g \la}] f~.
\end{equation}
Since $\phi$ is invariant to rotations in $G$ we verify that
$A_J$ commutes with rotations in $G$:
$A_J (g \circ f) = g \circ A_J f $ for all $g \in G$. 

In dimension $d=1$, $\G = \{-\One,\One\}$.
According to (\ref{wavedesin}), 
to build a complex wavelet $\psi$ concentrated on a single
frequency band, we set $\hat \psi(\om) = 0$ for $\om < 0$.
Following (\ref{consenas}), $W_J$ is unitary if and only if
\begin{equation}
\label{1Dsdns}
\beta\,\sum_{j \in \Z} |\hat \psi (2^{-j} |\om|)|^2 = 1~\mbox{and}~
|\hat \phi(\om)|^2 =  \beta\, \sum_{j=-\infty}^{0} |\hat \psi (2^{-j} |\om|)|^2~.
\end{equation}
If $\tilde \psi$ is a real wavelet which generates a 
dyadic orthonormal basis of
$\Ld(\R)$ \cite{Meyerbook} then 
$\hat \psi  = 2 \,\widehat {\tilde \psi}\,1_{\omega > 0}$
satisfies (\ref{consenas}).
Numerical examples in the paper are computed with a complex
wavelet $\psi$ calculated from a
cubic-spline 
orthogonal Battle-Lemari\'e wavelet $\tilde \psi$ \cite{Meyerbook}.

In any dimension $d \geq 2$,  $\hat \psi \in \LD$ can be defined
as a separable product in frequency polar 
coordinates $\om = |\om|\,\eta$, with $\eta$ in the unit sphere
${\bf S}^d$ of $\R^d$: 
\[
\forall (|\om|,\eta) \in \R^+ \times {\bf S}^d~~,~~
\hat \psi(|\om|\, \eta) = \hat \psi(|\om|)\, \gamma(\eta)~.
\]
The one-dimensional function $\hat \psi(|\om|)$ is chosen to
satisfy (\ref{1Dsdns}).
The Littlewood-Paley condition (\ref{consenas}) is then equivalent to
\[
\forall \eta \in {\bf S^d}~~,~~
\sum_{\ga \in G} |\gamma (\ga^{-1} \eta)|^2 = 1~.
\]

\subsection{Path Ordered Scattering}
\label{fastalgo}

Convolutions with wavelets
defines operators which are Lipschitz continuous to the action of 
diffeomorphisms, because 
wavelets are regular and localized functions.
However, a wavelet transform is not invariant to translations, and
$W[\la] f  = f \star \psi_\la$ translates when $f$ is translated.
The main difficulty is to compute translation invariant coefficients,
which remain stable to the action of diffeomorphisms,
and retain high frequency information provided by wavelets.
A scattering operator computes such a translation invariant representation.
We first explain how to build translation invariant coefficients from 
a wavelet transform, while maintaining
stability to the action of diffeomorphisms. Scattering operators
are then defined, and their main properties are summarized.
 
If $U[\la]$ is an operator defined on $\LD$,
not necessarily linear but which commutes
with translations, then
$\int U[\la] f(x)\, dx$ is translation invariant, if finite.
$W[\la] f = f \star \psi_\la$ commutes with translations but
$\int W[\la] f(x)\, dx = 0$ 
because $\int \psi(x)\, dx = 0$. More generally,
one can verify that any linear transformation of $W[\la] f$, which is
translation invariant, is necessarily zero.
To get a non-zero invariant, we 
set $U[\la]f = M[\la] W[\la]f$ where $M[\la]$  
is a non-linear ``demodulation'' which maps $W[\la] f$
to a lower frequency function having a non-zero integral.
The choice of $M[\la]$ must also preserve the Lipschitz continuity
to diffeomorphism actions.

If $\psi(x) = e^{i \eta \cdot x} \theta(x)$ then
$\psi_\la (x) = e^{i \la \eta \cdot x}\, \theta_\la (x)$, and hence
\begin{equation}
\label{wavemodfnsodf}
W[\la] f(x) = \Exp^{i  \la \eta \cdot  x}\, \, \Bigl( f^{\la}  \star \theta_{\la} (x) \Bigr)~
~\mbox{with}~~f^{\la} (x) = \Exp^{-i \la \eta \cdot x} \,f(x)~.
\end{equation}
The convolution $f^{\la}  \star \theta_{\la}$
is a low-frequency filtering because
$\hat \theta_\la (\om) = \hat \theta(\la^{-1} \om)$
covers a frequency
ball centered at $\om=0$, of radius proportional to $|\la|$.
A non-zero invariant can thus be obtained by
canceling the modulation term $\Exp^{i  \la \eta \cdot  x}$ with $M[\la]$. 
A simple example is:
\begin{equation}
\label{isdofns9d0nsdf}
M[\la] h(x) = e^{-i \la \eta \cdot x} e^{-i \Phi( \hat h(\la \eta))} h(x)~
\end{equation}
where $\Phi(\hat h(\la \eta))$ is the complex phase of $\hat h(\la \eta)$.
This non-linear phase registration guarantees that $M[\la]$ 
commutes with translations. It results from (\ref{wavemodfnsodf})
that $\int  {M[\la] W[\la] f(x)}\,dx = |\hat f(\la \eta)|\,|\hat \theta(0)|$.
It recovers the Fourier modulus representation, which is translation 
invariant but not Lipschitz continuous to diffeomorphisms
as shown in (\ref{instablsndf}).
Indeed, the demodulation operator $M[\la]$ in (\ref{isdofns9d0nsdf})
commutes with translations but
does not commute with the action of diffeomorphisms,
and in particular with dilations. 
The commutator norm of $M[\la]$ with a dilation 
is equal to $2$, even for arbitrarily small dilations, which explains
the resulting instabilities.

Lipschitz continuity under the action of diffeomorphisms is preserved if
$M[\la]$ commutes with the action of diffeomorphisms. For $\LD$ stability,
we also impose that $M[\la]$ is nonexpansive.
One can prove \cite{Bruna2}
that $M[\la]$ is then necessarily a pointwise operator, which means that 
$M[\la] h(x)$ only depends on the value of $h$ at $x$.
We further impose that $\| M[\la] h\| = \|h\|$ for all  $h \in \LD$, which then
implies that $|M[\la] h| = |h|$. The most regular functions are
obtained with $M[\la] h = |h|$, which eliminates all phase variations.
We derive from (\ref{wavemodfnsodf}) that
this modulus maps $W[\la] f$ into a lower frequency
envelop:
\[
M[\la] W[\la]  f = |W [\la] f| = |f^{\la}  \star \theta_{\la} |~. 
\]
Lower frequencies created by a modulus result from
interferences.
For example, if $f(x) = \cos(\xi_1 \cdot x) + a \,\cos(\xi_2 \cdot x)$ 
where
$\xi_1$ and $\xi_2$ are in the frequency band covered by $\hat \psi_\la$ then
$|f \star \psi_{\la}(x)| =  2^{-1} |\hat \psi_{\la} (\xi_1) \,  +
a\, \hat \psi_{\la} (\xi_2)\, e^{i (\xi_2-\xi_1) \cdot x}|$
oscillates at the interference frequency $|\xi_2 - \xi_1|$, which is
smaller than $|\xi_1|$ and $|\xi_2|$.

The integration
$\int U[\la]f(x) \, dx = \int |f \star \psi_\la(x)| \, dx$
is translation invariant but it removes all the
high frequencies of $|f \star \psi_\la(x)|$.
To recover these high frequencies, a scattering also computes the wavelet
coefficients of each $ U[\la] f$: 
$\{U[{\la}] f \star \psi_{\la'}\}_{\la'}$.
Translation invariant coefficients are again
obtained with a modulus 
$U[\la'] U[\la] f = |U[\la]f \star \psi_{\la'}|$
and an integration $\int U[\la'] U[\la] f(x) \, dx$.
If $f(x) = \cos(\xi_1 \cdot x) + a \,\cos(\xi_2 \cdot x)$ with $a < 1$, 
$|\xi_2 - \xi_1| \ll |\la|$ and $|\xi_2 - \xi_1|$ in the support of
$\hat \psi_{\la'}$ 
then $U[\la'] U[\la] f$ 
is proportional to $a\, |\psi_\la (\xi_1)|\, |\psi_{\la'} (|\xi_2 - \xi_1|)|$.
The second wavelet $\hat \psi_{\la'}$ captures 
the interferences created by the modulus, 
between the frequency components of $f$ 
in the support of $\hat \psi_{\la}$.
We now introduce the 
scattering propagator, which extends these decompositions.

\begin{definition}
\label{propdefin}
An ordered sequence $p = (\la_1,...,\la_m)$ with
$\la_k \in \cLa_\infty = 2^{\Z} \times G^+$ is called a path. 
The empty path is denoted $p = \nill$. Let $U[\la] f = |f \star \psi_\la|$
for $f \in \LD$.
A scattering propagator is a path ordered product of
non-commutative operators defined by
\begin{equation}
\label{relastn6}
U[p] = U[{\la_m}]\,...\,U[{\la_2}]\,U[{\la_1}]~~,
\end{equation}
with $U[\nill] = Id$. 
\end{definition}

The operator $U[p]$ is well defined on $\LD$ because
$\|U [\la] f\| \leq \|\psi_\la\|_1 \|f\|$ for all $\la \in \cLa_\infty$. 
The scattering propagator is a cascade of convolutions
and modulus:
\begin{equation}
\label{relastn62}
U[p] f = |~ |f \star \psi_{\la_1}| 
\star \psi_{\la_2}| \cdots |\star \psi_{\la_{m}}|~.
\end{equation}
Each $U[\la]$ filters the frequency component in the band
covered by $\hat \psi_\la$, and maps it to lower frequencies
with the modulus. The index sequence $p = (\la_1,...,\la_m)$ is
thus a frequency path variable. 
The scaling and rotation by $2^l g \in 2^\Z \times G$ 
of a path $p$ is written
$2^l g\, p = (2^l g \la_1, ...,2^l g \la_m)$.
The concatenation of two paths is denoted
$p+p' = (\la_1,...,\la_m,\la'_1,...,\la'_{m'})$, in particular
$p+\la= (\la_1,...,\la_m,\la)$.
It results from (\ref{relastn6}) that
\begin{equation}
\label{relastn08f}
U[p+p']  = U[p'] \,U[p]~ .
\end{equation}

Section \ref{scalittwnwavel} explains that 
if $f$ is complex valued then its wavelet transform is
$W_\infty f = \{W[\lambda] f \}_{\lambda \in \Lambda_\infty \atop 
-\lambda \in \Lambda_\infty}$ 
whereas if $f$ is real then
$W_\infty f = \{W[\lambda] f \}_{\lambda \in \Lambda_\infty}$.
If $f$ is complex then at the next iteration
$U[\la_1] f = |W[\la_1] f|$ is
real so next stage wavelet transforms are computed only for
$\la_k \in \cLa_\infty$.
The scattering propagator of a complex function is
thus defined over ``positive'' paths
$p = (\la_1,\la_2,...,\la_m) \in \cLa_\infty^m$ and 
``negative'' paths denoted $-p = (-\la_1,\la_2,...,\la_m)$. 
This is analogous to the positive and negative frequencies of a
Fourier transform.
If $f$ is real then $W[-\la_1] f = W[\la_1] f^*$
so $U[-\la_1] f = U[\la_1] f$ and hence $U[-p] f = U[p] f$. 
To simplify explanations, all results are proved on real
functions with 
scattering propagators restricted to positive paths.
These results apply to
complex functions by including negative paths.

\begin{definition}
\label{definwindoscas}
Let $\cP_\infty$ be the set of all finite paths.
The scattering transform of $f \in \LU$ is defined for
any $p \in \cP_\infty$ by
\begin{equation}
\label{firsndfinsdfsdfon}
\overline S f(p) = \frac 1 {\mu_p}\, \int U[p] f (x)\, dx
\,\,\,\mbox{with}\,\,\,\mu_p = \int U[p] \delta (x)\, dx\,.
\end{equation}
\end{definition}

A scattering is a 
translation-invariant operator which transforms 
$f\in \LU$ into a function of the frequency path variable $p$.
The normalization factor 
$\mu_p$ results from a path measure 
introduced in Section \ref{infdnsdfonsdsec}. Conditions are given
so that $\mu_p$ does not vanish.
This transform is then well-defined for any $f \in \LU$ and any $p$
of finite length $m$. Indeed $\|\psi_\la \|_1 = \|\psi \|_1$ 
so (\ref{relastn62}) implies that $\|U[p] f\|_1 \leq \|f\|_1\, \|\psi\|_1^m$. 
We shall see that a scattering transform has similarities with
Fourier transform modulus, where the path $p$ plays the role of a
frequency variable. However, 
as opposed to a Fourier modulus, a scattering transform 
is stable to the action of
diffeomorphisms, because it is computed by iterating on wavelet transforms
and modulus operators, which are stable.
For complex-valued functions, $\oS f$ is also defined on negative paths,
and $\oS f(-p) = \oS f(p)$ if $f$ is real.

If $p \neq \nill$ then
$\oS f (p) $ is non-linear but it preserves amplitude factors:
\begin{equation}
\label{ampldifnpre}
\forall \mu \in \R~~,~~\oS (\mu  f) (p)  = |\mu|\, \oS f  (p).
\end{equation}
A scattering has similar scaling and
rotation covariance properties as a Fourier transform.
If $f$ is scaled and rotated, ${2^l g} \circ f (x) = f(2^l g x)$, then
(\ref{wavedfn98}) implies that
$U[\la] (2^l g \circ f) = 2^l g \circ U[{2^{-l} g \la}] f$ and cascading
this result shows that
\begin{equation}
\label{pathcomansdfus29df}
\forall p \in \cP_\infty~~,~~
U[p] (2^l g \circ f) = 2^l g \circ U[{2^{-l} g p}] f~.
\end{equation}
Inserting this result in the definition (\ref{firsndfinsdfsdfon}) proves that
\begin{equation}
\label{pathcomansdfus2}
\overline S (2^l g \circ f) (p) = 2^{-dl} \,\oS f (2^{-l} g\, p) ~.
\end{equation}
Rotating $f$ thus rotates identically 
its scattering, whereas if $f$ is scaled by $2^l$ then
the frequency paths $p$ is scaled by $2^{-l}$.
The extension of the scattering transform in $\LD$ is done
as a limit of windowed scattering transforms, that we now introduce.

\begin{definition}
\label{windscatdef}
Let $J \in \Z$ and $\cP_J$ be the set of finite paths
$p = (\la_1,...,\la_m)$ with $\la_k \in \cLa_J$ and
hence $|\la_k| = 2^{j_k} > 2^{-J}$. A windowed scattering
transform is defined for all $p \in \cP_J$ by
\begin{equation}
\label{relastn5}
S_J [p] f(x) = U[p] f \star \phi_{2^J} (x) = \int U[p] f(u) \, \phi_{2^J} (x-u)\, du~.
\end{equation}
\end{definition}

The convolution with
$\phi_{2^J}(x) = 2^{-dJ} \phi(2^{-J} x)$ 
localizes the scattering
transform over spatial domains of size proportional to $2^{J}$:
\[
S_J [p] f(x) = |~ |f \star \psi_{\la_1}| 
\star \psi_{\la_2}| \cdots |\star \psi_{\la_{m}}| \star \phi_{2^J}(x)~.
\]
It defines
an infinite family of functions indexed by $\cP_J$, denoted
\[
S_J [\cP_J] f = \{ S_J [p] f \}_{p \in \cP_J}~.
\]
For complex-valued functions, 
negative paths are also included in $\cP_J$,
and $S_J[-p] f = S_J[p] f$ if $f$ is real.

Section \ref{upsdfnwd} 
proves that for appropriate wavelets,
$\|f\|^2 = \sum_{p \in \cP_J} \|S_J [p] f\|^2$.
However, the signal energy is mostly
concentrated on a much smaller set of frequency-decreasing
paths $p = (\la_k)_{k \leq m}$ for which 
$|\la_{k+1}| \leq |\la_k|$. Indeed, the propagator $U[\la]$
progressively pushes the energy towards lower frequencies.
The main theorem of
Section \ref{invasnsf2} proves that a windowed scattering is Lipschitz 
continuous to the action of diffeomorphisms.

Since $\phi(x)$ is continuous at $0$, if $f \in \LU$ then its 
windowed scattering transform converges pointwise to its scattering transform
when the scale $2^J$ goes to $\infty$:
\begin{equation}
\label{firsndfinsdf}
\forall x \in \R^d~,~\lim_{J \rightarrow \infty} 2^{dJ} 
\,S_J [p] f(x) = \phi(0)\,\int U[p] f(u)\, du = \phi(0)\,\mu_p\,\overline S(p)~.
\end{equation}
However, when $J$ increases,
the path set $\cP_J$ also increases. 
Section \ref{infdnsdfonsdsec} shows that
$\{ \cP_J \}_{J \in \Z}$ 
defines a multiresolution path approximation of a much larger set
$\cPinf$ including paths of infinite length. 
This path set is not countable as opposed to each $\cP_J$, and
Section \ref{infdnsdfonsdsec} introduces
a measure $\mu$ and a metric on $\cPinf$. 

Section \ref{Renonscat} extends
the scattering transform $\oS f (p)$ to all $f \in \LD$ and to all 
$p \in \cPinf$, and proves
that $\oS f \in \Ld(\cPinf, d\mu)$.
A sufficient condition is given to guarantee 
a strong convergence of $S_J f$ to $\overline S f$,
and it is conjectured that it is valid on $\LD$. 
Numerical examples 
illustrate this convergence and show that a scattering transform has
strong similarities with a Fourier transforms modulus, when mapping the path $p$
to a frequency variable $\omega \in \R^d$.

\subsection{Scattering Propagation and Norm Preservation}
\label{upsdfnwd} 

We prove that a windowed scattering 
$S_J$ is nonexpansive, and preserves the $\LD$ norm. 
Family of operators indexed by
a path set $\Omega$ are written $S_J[\Omega] = \{S_J[p]\}_{p \in \Omega}$
and $U[\Omega] = \{U[p]\}_{p \in \Omega}$.

\begin{figure}[htb]
\includegraphics[width=\columnwidth]{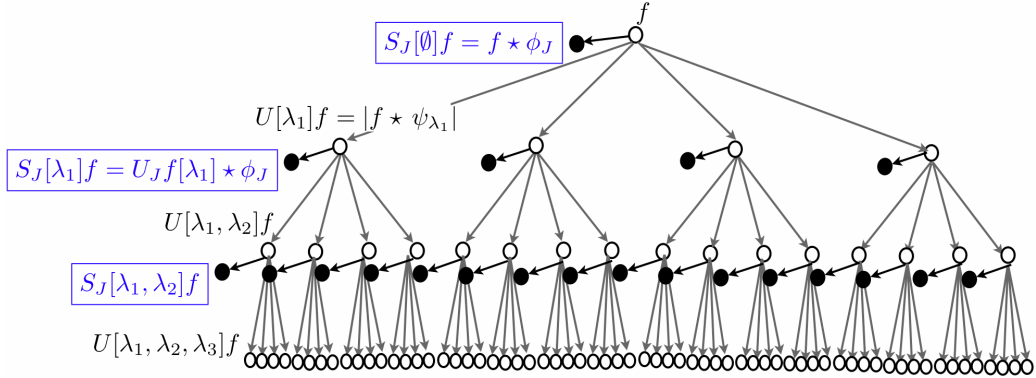}
	\caption{A scattering propagator $\oUJ$ applied to $f$ 
computes each
$U[\la_1] f = |f \star \psi_{\la_1}|$ and outputs 
$S_J [\emptyset] f = f \star \phi_{2^J}$.
Applying $\oUJ$ to each $U[\la_1] f$ computes all
$U[\la_1,\la_2] f$ and outputs 
$S_J[\la_1] = U[\la_1] \star \phi_{2^J}$.
Applying iteratively 
$\oUJ$ to each $U[p]f$ outputs $S_J[p]f = U[p]f \star \phi_{2^J}$
and computes the next path layer.}
	\label{scattering-cascade}
\end{figure}

A windowed scattering can be computed by iterating on the 
{\it one-step propagator} defined by
\[
\oU_J f = \{A_J f , ( U[\la] f )_{\la \in \cLa_J} \}\,,
\] 
with $A_J f = f \star \phi_{2^J}$ and $ U[\la] f = |f \star \psi_{\la}|$.
After calculating $U_J f$, 
applying again $\oU_J$ to each $U[\la] f$ 
yields a larger infinite family of functions. The decomposition is further 
iterated by recursively applying $\oU_J$ to each $U[p] f$. 
Since $U[\la] U[p] = U[p+\la]$ and $A_J U[p] = S_J [p]$, it results that
\begin{equation}
\label{propasndfs2}
\oU_J U[p] f = \{S_J [p] f \,,\,
( U[p+\la ] f )_{\la\in \cLa_J} \}~.
\end{equation}
Let $\cLa_J^m$ be the set of paths of length $m$, with
$\cLa_J^0 = \{ \nill \}$. It is propagated into
\begin{equation}
\label{propasndfs}
\oU_J U[\cLa_J^m] f  =  \{ S_J[\cLa_J^{m}] f \,,\, U[\cLa_J^{m+1}] f\}~.
\end{equation}
Since $\cP_J = \cup_{m \in \N} \cLa_J^{m}$, one can compute
$S_J[\cP_J] f$ from
$f = U [\nill] f$ by iteratively computing
$\oUJ U[\cLa_J^m] f$ for $m$ going from $0$ to $\infty$, as
illustrated in Figure \ref{scattering-cascade}.

Scattering calculations follow the general architecture of 
convolution neural-networks introduced by LeCun \cite{LeCun1}.
Convolution networks
cascade convolutions and a ``pooling'' non-linearity,
which is here the modulus
of a complex number. Convolution networks typically use
kernels that are not predefined functions such as wavelets, but which are 
learned with backpropagation algorithms.
Convolution network  architectures
have been successfully applied to number of recognition tasks
\cite{LeCun1} and are
studied as models for visual perception \cite{Poggio,poggio2}.
Relations between scattering operators and path formulations of quantum
field physics are also studied in \cite{Glinsky}.

The propagator $\oU_J f = \{A_J f , (|W[\la]f| )_{\la \in \cLa_J} \}$ 
is nonexpansive because
the wavelet transform $\oW_J$ is unitary and
a modulus is nonexpansive in the sense that 
$||a|-|b|| \leq |a-b|$ for any $(a,b) \in \C^2$.
This is valid whether $f$ is real or complex.
As a consequence
\begin{eqnarray}
\nonumber
\| \oUJ f -  \oUJ h \|^2 &=& \|A_J f - A_J h \|^2+ 
\sum_{\la \in \cLa_J}
\||W[\la] f| - |W[\la] h| \|^2 \\
& \leq & 
\label{inequasvdons}
\| \oWJ f -  \oWJ h \|^2 
\leq \|f - h \|^2  .
\end{eqnarray}
Since $\oWJ$ is unitary, setting $h = 0$ also proves
that $\|\oUJ f \| = \|f\|$, so $\oUJ$ preserves the norm.

For any path set $\Omega$ the norms of $S_J[\Omega] f$ and $U[\Omega] f$ are
\[
\|S_J[\Omega] f \|^2 = \sum_{p \in \Omega} \|S_J [p] f\|^2~~\mbox{and}~~
\|U[\Omega] f \|^2 = \sum_{p \in \Omega} \|U [p] f\|^2~.
\]
Since $S_J[\cP_J]$ iterates on $U_J$ which is nonexpansive, 
the following proposition derives that
$S_J[\cP_J]$ is also nonexpansive.

\begin{proposition}
\label{unidfnsdofnw}
The scattering transform is nonexpansive:
\begin{equation}
\label{onesnfsfd8}
\forall (f,h) \in \LD^2~~,~~\|S_J[\cP_J]f - S_J[\cP_J] h \| 
\leq \|f - h \|~.
\end{equation}
\end{proposition}.

{\it Proof:} 
Since $\oUJ$ is nonexpansive,
it results from (\ref{propasndfs}) that
\begin{eqnarray*}
\|U [\cLa^m_J]f -   U[\cLa^m_J] h \|^2 &\geq&
\| \oUJ U [\cLa^m_J]f -  \oUJ U[\cLa^m_J] h \|^2 \\
&=&\| S_J[ \cLa^{m}_J] f - S_J[\cLa^{m}_J] h\|^2 + 
\| U[ \cLa^{m+1}_J] f - U[\cLa^{m+1}_J] h\|^2 ~.
\end{eqnarray*}
Summing these equations for $m$ going from $0$ to $\infty$ proves that
\begin{equation}
\label{Oncfondfos8ufs}
\|S_J[\cP_J] f - S_J[\cP_J] h 
\|^2 = \sum_{m=0}^{\infty} \|S_J [\cLa^m_J]f -   S_J [\cLa^m_J] h \|^2 \leq 
\|f - h \|^2 ~. ~\Box
\end{equation}

Section \ref{fastalgo} explains that each
$U[\la] f = |f \star \psi_{\la}|$ captures the frequency energy of $f$
over a frequency band covered by $\hat \psi_\la$ and 
propagates this energy towards lower frequencies. 
The following theorem proves this result by showing that
the whole scattering energy ultimately reaches the minimum 
frequency $2^{-J}$ and is trapped by the low-pass filter $\phi_{2^J}$.
The propagated scattering energy thus
goes to zero as the path length increases, and the theorem derives that
$\|S_J[\cP_J] f \| = \|f\|$.
This result also applies to complex-valued functions by incorporating 
negative paths $(-\la_1,\la_2,...,\la_m)$ in $\cP_J$.

\begin{theorem}
\label{energydecth}
A scattering wavelet $\psi$ is said to be admissible
if there exists $\eta \in \R^d$ 
and $\rho  \geq 0$, with $|\hat \rho (\om)| \leq |\hat \phi (2 \om)|$ and
$\hat \rho(0) = 1$, such that the function
\begin{equation}
\label{conditionprogredf0}
\hat \Psi (\omega) = 
| \hat \rho (\om - \eta)|^2 - \sum_{k=1}^{+\infty} 
k\, \Big(1 - |\hat \rho (2^{-k} (\om - \eta))|^2\Big)~
\end{equation}
satisfies
\begin{equation}
\label{conditionprogre}
\alpha = \inf_{1 \leq |\omega|\leq 2}\sum_{j=-\infty}^{+\infty} 
\sum_{\ga \in \G}
\hat \Psi (2^{-j} \ga^{-1} \om) \,|\hat \psi (2^{-j} \ga^{-1} \om)|^2 
 > 0~.
\end{equation}
If the wavelet is admissible then 
\begin{equation}
\label{dnslnfs8302}
\forall f \in \LD~~,~~\lim_{m \rightarrow \infty} \| U[\cLa^m_J] f\|^2 = 
\lim_{m \rightarrow \infty} \sum_{n=m}^\infty \| S_J[\cLa^n_J] f\|^2 = 0
\end{equation}
and
\begin{equation}
\label{conditionprogredf}
\|S_J[\cP_J] f\| = \|f \|~.
\end{equation}
\end{theorem}

{\it Proof:} 
We first prove that 
$\lim_{m \rightarrow \infty} \| U[\cLa^m_J] f\| = 0$ is equivalent to having
$\lim_{m \rightarrow \infty} \sum_{n=m}^\infty \| S_J[\cLa^n_J] f\|^2 = 0$
and to $\|S_J[\cP_J] f \| = \|f\|$.
Since $\|\oUJ h \| = \|h \|$ for any $h \in \LD$ and $\oUJ U[\cLa_J^n] f  =  \{S_J[\cLa_J^{n}] f \,,\,U[\cLa_J^{n+1}]\}$,
\begin{equation}
\label{Unfidnfsdf0}
\| U [\cLa^n_J]f \|^2 = \| \oUJ U [\cLa^n_J]f \|^2 = 
\| S_J[\cLa^{n}_J] f\|^2 + 
\| U[\cLa^{n+1}_J] f \|^2 .
\end{equation}
Summing 
for $m \leq n < \infty$ proves that
$\lim_{m \rightarrow \infty} \| U[\cLa^m_J] f\| = 0$ is equivalent to
$\lim_{m \rightarrow \infty} \sum_{n=m}^\infty \| S_J[\cLa^n_J] f\|^2 = 0$.
Since $f = U[\cLa^0_J] f$, summing (\ref{Unfidnfsdf0}) for $0 \leq n < m$
also proves that
\begin{equation}
\label{Unfidnfsdf7}
\|f \|^2 = \sum_{n=0}^{m-1} \| S_J[\cLa^n_J] f\|^2 + \| U[\cLa^m_J] f\|^2 ~,
\end{equation}
so $\|S_J[\cP_J] f\|^2 = \sum_{n=0}^{\infty} \| S_J[\cLa^n_J] f\|^2 = \|f\|^2$
if and only if $\lim_{m \rightarrow \infty} \|U[\cLa_J^m] \| = 0$.

We now prove that condition (\ref{conditionprogredf0})  implies
that $\lim_{m \rightarrow \infty} \| U[\cLa_J^m] f\|^2 = 0 $.
It relies on the following lemma, which
gives a lower bound of $|f \star \psi_{\la}|$ convolved
with a positive function.

\begin{lemma}
\label{boundendf}
If $h \geq 0$ 
then for any $f \in \LD$
\begin{equation}
\label{dnslnfs}
 |f \star \psi_{\la}| \star h  \geq \sup_{\eta \in \R^d}
| f \star \psi_{\la} \star h_\eta  |~~\mbox{with}~~h_\eta (x)
= h(x)\,e^{i \eta x}~.
\end{equation}
\end{lemma}

The lemma is proved by computing
\begin{eqnarray*}
 |f \star \psi_{\la}| \star h (x)  &=& 
\int \left|\int f(v) \psi_{\la}(u-v) \, dv \right| \,h(x-u)\,du
\\
&=& \int \left|\int f(v) \psi_{\la}(u-v) \, e^{i\eta(x-u)}\,
h(x-u)\,dv \right| \,du
\\
&\geq & \left|\int \int f(v) \psi_{\la}(u-v) \, h(x-u)\,e^{i\eta(x-u)}\,dv \,du
\right|\\
&= & \left|\int f(v) \int \psi_{\la}(x-v-u') \, h(u')\,e^{i\eta u'}\,du' dv
\right|\\
& = &
\left|\int f(v) \psi_{\la} \star h_\eta (x-v) dv \right| 
 =   | f \star \psi_{\la} \star h_\eta | ~,
\end{eqnarray*}
which finishes the lemma's proof.

Appendix \ref{proofLemma8} uses this lemma to show that the
scattering energy propagates progressively towards lower frequencies,
and proves the following lemma.

\begin{lemma}
\label{lemma8}
If (\ref{conditionprogre}) is satisfied and
\begin{equation}
\label{condroref5}
\|f\|_w^2 = \sum_{j=0}^\infty \sum_{\ga \in \G^+} j\, \|W[{2^j\ga}] f \|^2 < \infty~
\end{equation}
then
\begin{equation}
\label{condroref}
\frac {\alpha} 2 \, \|U[\cP_J]f \|^2 \leq \max(J+1,1)\,
\|f \|^2 + \|f\|_w^2~.
\end{equation}
\end{lemma}



The class of function for which $\|f\|_w < \infty$ is a
logarithmic Sobolev class, corresponding
to functions having an average modulus of continuity in $\LD$.
Since
\[
\|U[\cP_J]f \|^2 = \sum_{m=0}^{+\infty} \|U[\cLa_J^m]f \|^2,
\]
if $\|f\|_w < \infty$ then  
(\ref{condroref}) implies that 
$\lim_{m \rightarrow \infty} \|U[\cLa_J^m]f \| = 0$.
This result is extended in $\LD$ by density.
Since $\phi \in \LU$ and $\hat \phi(0) = 1$, 
any $f \in \LD$ satisfies
$\lim_{n \rightarrow -\infty} \|f - f_n \| = 0$,
where $f_n = f \star \phi_{2^n}$ and $\phi_{2^n} (x) = 2^{-nd} \phi(2^{-n} x)$.
We prove that $\lim_{m \rightarrow \infty} \| U[\cLa_J^m] f_n\|^2  = 0$
by showing that $\|f_n \|_w < \infty$. 
Indeed  
\begin{eqnarray*}
\|W[2^j\ga] f_n \|^2 &=& \int |\hat f (\om)|^2\, |\hat \phi(2^{n} \om)|^2\,
|\hat \psi (2^{-j} \ga^{-1} \om)|^2 \, d\om \\
&\leq&
C\,2^{-2n-2j} \int |\hat f (\om)|^2\, d\om ,
\end{eqnarray*}
because $\psi$ has a vanishing moment so
$|\hat \psi (\om)| = O(|\om|)$, and 
the derivatives of $\phi$ are in $\LU$ so
$|\om|\, |\hat \phi(\om)|$ is bounded. It results 
that $\|f_n \|_w < \infty$. 

Since $U[\cLa^m]$ is nonexpansive
$\|U [\cLa^m_J]f -   U[\cLa^m_J] f_n \| \leq \|f - f_n \|$ so
\[
\| U[\cLa^m_J] f \|
\leq \|f - f_n \| + \| U[\cLa^m_J] f_n \|.
\]
Since $\lim_{n \rightarrow -\infty} \|f - f_n \| = 0$ and
$\lim_{m \rightarrow \infty} \| U[\cLa^m_J] f_n\|  = 0$
it results that 
$\lim_{m \rightarrow \infty} \| U[\cLa^m_J] f\|^2  = 0$
for any $f \in \LD$.
$\Box$

The proof shows that the scattering
energy propagates progressively towards lower frequencies.
The energy of $U [p] f$
is mostly concentrated along frequency-decreasing paths
$p = (\la_k)_{k \leq m}$ for which $|\la_{k+1}| < |\la_{k}|$. 
For example, if $f = \delta$ then paths of length $1$ have an energy
$\|U [2^j r] \delta \|^2 = \| \psi_{2^j r} \|^2 
= 2^{-dj} \|\psi\|^2$. This energy is 
then propagated among all paths $p \in \cP_J$.
For a cubic spline wavelet in dimension $d = 1$, over $99.5\%$ of this energy
is concentrated along frequency-decreasing paths.
Numerical implementations
of scattering transforms thus limits computations to these frequency
decreasing paths. 
The scattering transform of a signal of size $N$ is computed 
along all frequency-decreasing paths, with $O(N \log N)$ operations,
by using a filter bank implementation \cite{MallatEUSIPCO}.

The decay of 
$\sum_{n=m}^\infty \| S_J[\cLa^n_J] f\|^2$ 
implies that we can neglect all paths of length larger than 
some $m > 0$. 
The numerical decay of $\|S_J [\cLa^n_J] f\|^2$ appears to be exponential
in image and audio processing applications.
The path length is limited to $m = 3$ in classification 
applications \cite{Joakim,Bruna}.

Theorem \ref{energydecth} requires a unitary wavelet transform and
hence an admissible wavelet which satisfies
Littlewood-Paley condition 
$\beta \sum_{(j,r) \in \Z\times \G}
|\hat \psi (2^j \ga \om)|^2  = 1$.
There must also exist
$\rho \geq 0$ and $\eta \in \R^d$
with $|\hat \rho (\om)| \leq |\hat \phi (2 \om)|$
such that 
$\sum_{(j,r) \in \Z\times \G}
|\hat \psi (2^j \ga \om)|^2 |\hat \rho(2^j \ga \om-\eta)|^2 $ 
is sufficiently large so that $\alpha > 0$.
This can be obtained if according to
(\ref{wavedesin}), $\psi (x) = e ^{i \eta x} \theta(x)$ and hence
$\hat \psi(\om) = \hat \theta(\om - \eta)$, where $\hat \theta$ 
and $\hat \rho $ have their energy concentrated over nearly the same 
low frequency domains.
For example, an analytic cubic spline Battle-Lemari\'e 
wavelet is admissible in one dimension with $\eta = 3  \pi/2$.
This is verified by choosing 
$\rho$ to be a positive cubic box spline, in which case a
numerical evaluation of (\ref{conditionprogre}) gives
$\alpha = 0.2766 > 0$.

\subsection{Translation Invariance}
\label{invasnsf1}

We show that the scattering distance
$\|S_J [\oP_J] f - S_J [\oP_J] h\|$ 
is non-increasing when $J$ increases, and thus converges when $J$ goes
to $\infty$. It defines a limit distance which is proved to be translation
invariant. Section \ref{infdnsdfonsdsec} studies
the convergence of $S_J [\cP_J] f$ when $J$ goes to $\infty$,
to the translation invariant scattering transform $\oS f$.

\begin{proposition}
\label{Undfonsdf}
For all $(f,h) \in \LD^2$ and $J \in \Z$
\begin{equation}
\label{8dfv0s1}
\|S_{J+1}[\cP_{J\pl 1}] f - S_{J+1} [\cP_{J\pl 1}] h \|
\leq \|S_J[\cP_J] f - S_J[\cP_J] h \|~.
\end{equation}
\end{proposition}

{\it Proof:} 
Any 
$p' \in \cP_{J\pl 1}$ can uniquely be written as an extension of 
a path $p \in \cP_J$ where
$p $ is the longest prefix
of $p'$ which belongs to $\cP_J$,
and $p' = p + q$ for some $q \in \cP_{J\pl 1}$.
The set of all extensions of $p \in \cP_{J}$ in $\cP_{J+1}$ is
\begin{equation}
\label{extenset}
\cP_{J\pl 1}^p = \{p\} \cup \{p + 2^J\ga + p''\}_{\ga \in \G^+,\,
p'' \in \cP_{J\pl 1}} ~. 
\end{equation}
It defines 
a non-intersecting partition 
of $\cP_{J\pl 1} = \cup_{p \in \cP_J} \cP^p_{J\pl 1}$.
We shall prove
that such extensions are nonexpansive:
\begin{equation}
\label{pathconsdf9u8sd}
\sum_{p' \in \cP^p_{J\pl 1}} 
\|S_{J\pl 1}[p'] f - S_{J\pl 1}[p'] h\|^2 
\leq \|S_J[p] f - S_J[p] h\|^2 .
\end{equation}
To later prove Proposition \ref{bornsdfnsdf}, we
also verify that it preserves energy
\begin{equation}
\label{pathconsdf9u8sd2}
\sum_{p' \in \cP^p_{J\pl 1}} \|S_{J\pl 1}[p'] f\|^2 =
\|S_J[p] f\|^2  ~.
\end{equation}
Summing (\ref{pathconsdf9u8sd}) on all $p \in \cP_J$ proves (\ref{8dfv0s1}).

Appendix \ref{proofLemma8} proves in 
(\ref{conservq}) that for all $g\in \LD$
\[
\|g \star \phi_{2^{J\pl 1}} \|^2 + \sum_{\ga \in G^+} \|g \star \psi_{2^J\ga}  \|^2 
= \| g \star \phi_{2^J}\|^2 ~. 
\]
Applying it to $g = U[p] f - U[p] h$ 
together with $U[p] f \star \phi_{2^J} = S_J[p] f$ and
$U[p] f \star \psi_{2^J \ga} = U[p+ 2^J r] f$ gives
\begin{eqnarray}
\label{formsidfonsd8}
\|S_J[p] f - S_J[p] h\|^2 &=& \|S_{J\pl 1}[p] f - S_{J\pl 1}[p] h \|^2 \\
\nonumber
&&+ \sum_{r \in G^+} \|U[p + 2^J r]f-U[p + 2^J r]h  \|^2 
\end{eqnarray}
Since
$S_{J\pl 1}[\cP_{J\pl 1}] U[p+2^Jr ] f = \{S_{J\pl 1}[p+2^Jr+p'']\}_{p''\in \cP_{J\pl 1}}$,
and $S_{J+1}[\cP_{J\pl 1}]f$ is nonexpansive,
it implies
\begin{eqnarray*}
\lefteqn{
\|S_J[p] f - S_J[p] h\|^2 \geq \|S_{J\pl 1}[p] f - S_{J\pl 1}[p] h \|^2 }\\
&&+ \sum_{p'' \in \cP_{J\pl 1}} \sum_{r \in G^+} \|
S_{J\pl 1}[p+2^J r+p'']f- S_{J\pl 1}[p+2^J r + p'']h \|^2  ,
\end{eqnarray*}
which proves (\ref{pathconsdf9u8sd}).
Since $S_J[\cP_{J\pl 1}]f$ preserves the norm, setting
$h = 0$ in
(\ref{formsidfonsd8}) gives
\[
\| S_J[p] f\|^2 =
\|S_{J\pl 1}[p] f \|^2 + 
\sum_{p'' \in \cP_{J\pl 1}} \sum_{G^+} \|
S_{J\pl 1}[p+2^J r+p'']f \|^2  ~,
\]
which proves (\ref{pathconsdf9u8sd2}).
$\Box$

This proposition proves that $\|S_J[\cP_J] f - S_J[\cP_J] h \|$ 
is positive and non-increasing when $J$ increases, and thus
converges.
Since $S_J[\cP_J]$ is nonexpansive,
the limit metric is also nonexpansive
\[
\forall (f,h) \in \LD^2~~,~~
\lim_{J \rightarrow \infty} \|S_J[\cP_J] f - S_J[\cP_J] h \|
\leq \|f - h \|~.
\]
For admissible scattering wavelets
which satisfy (\ref{conditionprogre}), Theorem \ref{energydecth} proves
that $\|S_J[\cP_J]f \| = \|f\|$ so
$\lim_{J \rightarrow \infty} \|S_J[\cP_J] f\| = \|f\|$.
The following theorem proves that the limit metric is translation invariant.

\begin{theorem}
\label{theoInva}
For admissible scattering wavelets
\[
\forall f \in \LD~\forall c \in \R^d~,~
\lim_{J \rightarrow \infty} \|S_J[\cP_J] f - S_J[\cP_J] L_c f \| = 0~.
\]
\end{theorem}

{\it Proof:} 
Since $S_J[\cP_J] \,L_c = L_c\, S_J[\cP_J]$ and
$S_J[\cP_J] f = A_J\, U[\cP_J] f$
\begin{eqnarray}
\|S_J[\cP_J] L_c f - S_J[\cP_J]f \| &=& 
\|L_c A_J U[\cP_J] f - A_J U[\cP_J] f \| \nonumber \\
\label{indfnsdf8df}
&\leq &\|L_c A_J - A_J \| \, \|U[\cP_J] f \|\,.
\end{eqnarray}

\begin{lemma}
\label{lemma1}
There exists $C$ such that 
for all $\tau \in \bCd$ with $\|\nabla \tau \|_\infty \leq 1/2$ we have
\begin{equation}
\label{lemma1eq1}
\| L_\tau A_J f - A_J f \| 
\leq {C\,\|f\|\,2^{-J}~ \|\tau\|_\infty}~.
\end{equation}
\end{lemma}
This lemma is proved in Appendix \ref{proofLemma1}.
Applying it to $\tau = c$ and hence
$\|\tau \|_\infty = |c|$ proves that
\begin{equation}
\label{indfnsdf8df2}
\| L_c A_J - A_J \| 
\leq {C\,\,2^{-J}~ |c|}~.
\end{equation}
Inserting this in (\ref{indfnsdf8df}) gives
\begin{equation}
\label{indfnsdf8df31}
\|L_c S_J[\cP_J]f - S_J[\cP_J]f \|
\leq C\,\,2^{-J}\, |c|\,\|U[\cP_J] f \|\,.
\end{equation}

Since the admissibility condition (\ref{conditionprogre}) is satisfied,
Lemma \ref{lemma8} proves in (\ref{condroref}) that for $J > 1$
\begin{equation}
\label{foradfonsdf}
\frac{\alpha} 2 \, \|U[\cP_J]f \|^2 \leq 
(J+1)\, \|f\|^2 + \|f\|_w^2.
\end{equation}
If $\|f\|_w < \infty$ then 
it results from (\ref{indfnsdf8df31}) that
\[
\|L_c S_J[\cP_J]f - S_J[\cP_J]f \|^2 \leq ((J+1)\, \|f\|^2 + \|f\|_w^2) \,
{C^2\,2\,\alpha^{-1}\,2^{-2J}~ |c|^2}
\]
so $\lim_{J \rightarrow \infty} \|L_c S_J[\cP_J]f - S_J[\cP_J]f \| = 0$.

We then prove that  $\lim_{J \rightarrow \infty} \|L_c S_J[\cP_J]f - S_J[\cP_J]f \| = 0$ for all  $f \in \LD$, with a similar density argument as
in the proof of Theorem \ref{energydecth}. Any $f \in \LD$ can
be written as a limit of $\{f_n \}_{n \in \N}$ with $\|f_n \|_w < \infty$,
and since $S_J [\cP_J]$ is nonexpansive and $L_c$ unitary, one can verify that
\[
\|L_c S_J[\cP_J]f - S_J[\cP_J]f\| \leq
\|L_c S_J[\cP_J]f_n - S_J[\cP_J]f_n \| + 2\, \|f - f_n\|~.
\]
Letting $n$ go to $\infty$ proves that
$\lim_{J \rightarrow \infty} \|L_c S_J[\cP_J]f - S_J[\cP_J]f \| = 0$,
which finishes the proof.$\Box$

\subsection{Lipschitz Continuity to Actions of Diffeomorphisms}
\label{invasnsf2}

This section proves that a windowed 
scattering is Lipschitz continuous
to the action of diffeomorphisms.
A diffeomorphism of $\R^d$ sufficiently close to a translation maps
$x$ to $x - \tau (x)$ where $\tau (x)$ is a displacement field
such that $\|\nabla \tau \|_\infty < 1$. 
The diffeomorphism action on $f \in \LD$ is 
$L_\tau f (x) = f(x- \tau(x))$. 
The maximum increment of $\tau$ is denoted
$\|\Delta \tau \|_\infty = \sup_{(x,u) \in \R^{2d}} |\tau (x) - \tau (u)|$.
Let $S_J$ be a windowed scattering operator computed with an
admissible scattering wavelet which satisfies (\ref{conditionprogre}).
The following theorem computes an upper bound
of $\|S_J [\cP_J] L_\tau f - S_J [\cP_J] f \|$ as a function of
a mixed $(\lu,\LD)$ scattering norm: 
\begin{equation}
\label{err1099}
\|U [\cP_J] f \|_{1} 
= \sum_{m=0}^{+\infty} \|U [\cLa_J^m]f \|~.
\end{equation}
We denote $\cP_{J,m}$ the subset of $\cP_J$ of paths 
of length strictly smaller than $m$, and $(a \vee b) = \max(a,b)$.

\begin{theorem}
\label{DefomsProps}
There exists $C$ such that all $f \in \LD$ with
$\|U[\cP_J] f \|_{1} < \infty$ and all $\tau \in \bCd$ with
$\|\nabla \tau\|_\infty \leq 1/2$ satisfy
\begin{equation}
\label{err1090}
\| S_J[\cP_J] L_\tau f -  S_J[\cP_J] f   \| \leq  C\,\|U[\cP_J] f \|_{1} \,
K (\tau)
\end{equation}
with
\begin{equation}
\label{err10901}
K (\tau) = 
2^{-J} \|\tau \|_\infty + \|\nabla \tau \|_\infty (\log 
\frac{\|\Delta \tau\|_\infty}{\|\nabla \tau \|_\infty } \vee 1)
+ \|\rH \tau \|_\infty ~,
\end{equation}
and for all $m \geq 0$
\begin{equation}
\label{err10902}
\| S_J [\cP_{J,m}] L_\tau f -  S_J [\cP_{J,m}] f   \| 
\leq C\,m\,\|f\| \, K(\tau)~.
\end{equation}
\end{theorem}

{\it Proof:}
Let $[S_J[\cP_J],L_\tau] = S_J[\cP_J] \,L_\tau - L_\tau \, S_J[\cP_J]$,
\begin{equation}
\label{errnddfosnf001}
\| S_J[\cP_J]L_\tau f  -S_J[\cP_J]f \|  \leq 
\| L_\tau S_J[\cP_J]f -S_J[\cP_J]f \| + \| [S_J[\cP_J]\, ,\,L_\tau ] f\|~.
\end{equation}
Similarly to (\ref{indfnsdf8df})
the first term on the right satisfies
\begin{equation}
\label{err10df8s}
\| L_\tau S_J[\cP_J]f -S_J[\cP_J]f \| \leq \| L_\tau A_J  - A_J \|\,
\|U[\cP_J] f \|~.
\end{equation}
Since 
\[
\|U[\cP_J] f \| =
\left( \sum_{m=0}^{+\infty} \|U[\cLa_J^m] f \|^2\right)^{1/2}  \leq 
\sum_{m=0}^{+\infty} \|U[\cLa_J^m] f \|
\]
it results that
\begin{equation}
\label{err10df8s9}
\| L_\tau S_J[\cP_J]f -S_J[\cP_J]f \| 
\leq \| L_\tau A_J  - A_J \|\,\|U[\cP_J]f\|_{1}~.
\end{equation}
Since $S_J[\cP_J]$ iterates on $U_J$ which
is nonexpansive, Appendix \ref{Proofsdfn} proves the following upper
bound on scattering commutators.

\begin{lemma}
\label{LemmaComsnuta}
For any operator $L$ on $\LD$
\begin{equation}
\label{comsndfwnsf}
\| [S_J[\cP_J]\,,\, L] f \|
 \leq \|U[\cP_J] f \|_{1} \,\|[U_J\,,\,L]\|~.
\end{equation}
\end{lemma}

The operator $L = L_\tau$ also satisfies
\begin{equation}
\label{comsndfwnsf00}
\|[U_J,L_\tau] \| \leq \|[\oWJ,L_\tau] \| ~.
\end{equation}
Indeed, $\oUJ = M\,\oWJ$,
where $M \{h_J , (h_\la )_{\la \in \cLa_J}\} = \{h_J , (|h_\la| )_{\la \in \cLa_J}\}$
is a nonexpansive modulus operator. 
Since $M L_\tau = L_\tau M$ 
\begin{equation}
\label{err10dfi79}
\|[\oUJ , L_\tau] \|= \|M_J\,[\oWJ,L_\tau]\| \leq \|[\oWJ,L_\tau]\|~.
\end{equation}
Inserting (\ref{comsndfwnsf}) with (\ref{comsndfwnsf00}) 
and (\ref{err10df8s9}) in (\ref{errnddfosnf001}) gives
\begin{equation}
\label{err10}
\| S_{J} [\cP_J]L_\tau f -  S_{J} [\cP_J]f   \| \leq  \|U[\cP_J]f \|_{1} \Bigl(
\| L_\tau A_J - A_J \| + \|[ \oWJ,L_\tau] \| \Bigr)~.
\end{equation}
Lemma \ref{lemma1} proves that
$\| L_\tau A_J  - A_J \| \leq C\,2^{-J} \, {\|\tau \|_\infty}$.
This inequality and (\ref{err10}) imply that
\begin{equation}
\label{err10dfnsd2}
\| S_{J} [\cP_J]L_\tau f -  S_{J} [\cP_J]f   \| \leq  C\,\|U[\cP_J]f \|_{1} \Bigl(
2^{-J} \, {\|\tau\|_\infty} + \|[ \oWJ,L_\tau] \| \Bigr)~.
\end{equation}
To prove (\ref{err1090}), 
the main difficulty is to compute an upper bound of 
$\|[ \oWJ,L_\tau] \|$, and hence of
$\|[ \oWJ,L_\tau] \|^2 
= \|[ \oWJ , L_\tau]^*\,[ \oWJ , L_\tau] \|$,
where $A^*$ is the adjoint of an operator $A$. 
The wavelet commutator applied to $f$ is 
\[
[ \oWJ ,L_\tau] f = \{ [A_J,L_\tau] f ~,~([W[\la],L_\tau] f)_{\la \in \cLa_J} \}\,,
\]
whose norm is
\begin{equation}
\label{commnaf0}
\|[ \oWJ , L_\tau] f\|^2 = \|[A_J , L_\tau] f\|^2 + \sum_{\la \in \cLa_J}
\|[W[\la] , L_\tau] f\|^2 .
\end{equation}
It results that
\[
[ \oWJ , L_\tau]^*\,[ \oWJ , L_\tau] 
= [A_J , L_\tau]^*\,[A_J , L_\tau] +
\sum_{\la \in \cLa_J} [W[\la] , L_\tau]^*\,[W[\la] , L_\tau] .
\]
The operator 
$[ \oWJ , L_\tau]^*\,[ \oWJ , L_\tau]$ has a singular kernel along the diagonal 
but Appendix \ref{prooftheosdfn} proves that its norm is bounded.

\begin{lemma}
\label{thedw0}
There exists $C > 0$ 
such that all $J \in \Z$ and all $\tau \in \bCd$ with
$\|\nabla \tau\|_\infty\leq 1/2$ satisfy
\begin{equation}
\label{commnaf}
\|[ \oWJ , L_\tau] \| \leq 
C\, \Bigl(\|\nabla \tau \|_\infty (\log \frac{\| \Delta \tau\|_\infty}{\|\nabla \tau \|_\infty }
\vee 1)+ \|\rH \tau \|_\infty\Bigr) ~.
\end{equation}
\end{lemma}

Inserting the wavelet commutator bound (\ref{commnaf}) in
(\ref{err10dfnsd2}) proves the theorem inequality (\ref{err1090}).
One can verify that (\ref{err1090}) remains valid when replacing
$\cP_J$ by the subset of paths of length smaller than $m$:
$\cP_{J,m} = \cup_{n< m} \cLa_J^n$, if we replace 
$\|U[\cP_{J}] f \|_1$ by $\|U[\cP_{J,m}] f \|_1$. 
The inequality (\ref{err10902}) results from
\begin{equation}
\label{indwansford9}
\|U[\cP_{J,m}] f \|_1 = \sum_{n=0}^{m-1} \|U [\cLa_J^n]f \| \leq m\,\|f\|~.
\end{equation}
This is obtained by observing that
\begin{equation}
\label{indwansford}
\|U[\cLa_J^n]f \| \leq 
\|U[\cLa_J^{n-1}]f \| \leq \|f \| ~,
\end{equation}
because $U[\cLa_J^n]f$ is computed in (\ref{propasndfs2}) by applying 
the norm-preserving operator $\oUJ$ on $U[\cLa_J^{n-1}]f$.
$\Box$

The condition $\|\nabla \tau \|_\infty \leq 1/2$ can be 
replaced by $\|\nabla \tau \|_\infty < 1$ if
$C$ is replaced by $C\,(1-\|\nabla \tau \|_\infty)^{-d}$. 
Indeed $\|S_J[\cP_J] f \| = \|f \|$ and 
$\|S_J[\cP_J] L_\tau f \| \leq \|f \|(1-\|\nabla \tau \|_\infty)^{-d}$.
This remark applies to all subsequent theorems 
where the condition $\|\nabla \tau\|_\infty \leq 1/2$ appears.
The theorem proves that
the distance
$\| S_J[\cP_J] L_\tau f -  S_J[\cP_J] f   \|$ produced by the
diffeomorphism action $L_\tau$ is bounded by 
a translation term proportional to $2^{-J} \|\tau\|_\infty$ and
a deformation error proportional
to $\|\nabla \tau\|_\infty$. This deformation term results from the
wavelet transform commutator $[\oWJ,L_\tau]$.
The term
$\log (\| \Delta \tau\|_\infty/\|\nabla \tau \|_\infty )$ can also  
be replaced by $\max(J,1)$ in the proof of Theorem \ref{DefomsProps}.
For compactly supported functions $f$, Corollary \ref{Theoelas} 
replaces this term by the log of the support radius. 

If $f \in \LD$  has a weak form of regularity such as an 
average modulus of continuity in $\LD$ then
Lemma \ref{lemma8} proves that 
$\|U[\cP_J] f\|^2 = \sum_{n=0}^{\infty} \|U[\cLa_J^n] f\|^2$ is finite.
Numerical experiments indicate that 
$\|U[\cLa_J^n] f\|$ has exponential decay 
for a large class of functions,
but we do not characterize here
the class of functions for which 
$\|U[\cP_J] f\|_1 = \sum_{n=0}^{\infty} \|U[\cLa_J^n] f\|$ is finite.
In audio and image processing 
applications \cite{Joakim,Bruna}, 
the percentage of scattering energy becomes negligible
over paths of length larger than $3$ so (\ref{err10902}) 
is applied with $m=4$.

The following corollary derives from Theorem \ref{DefomsProps} that a
windowed scattering is Lipschitz continuous to the action of diffeomorphisms
over compactly supported functions.

\begin{corollary}
\label{Theoelas}
For any compact $\Omega \subset \R^d$ 
there exists $C$ such that for all $f \in \LD$ supported in $\Omega$
with $\|U[\cP_J] f \|_{1} < \infty$ 
and for all $\tau \in \bCd$ with $\|\nabla \tau\|_\infty \leq 1/2$, if
$2^J \geq \frac{\|\tau\|_\infty} {\|\nabla \tau \|_\infty }$ then 
\begin{equation}
\label{err2fsd2d2}
\| S_J[\cP_J] L_\tau f -  S_J[\cP_J] f   \| \leq  
C\,\|U[\cP_J] f \|_{1} \,
\Bigl(\| \nabla \tau \|_\infty + \| \rH \tau \|_\infty \Bigr)~.
\end{equation}
\end{corollary}

{\it Proof:} 
The inequality (\ref{err2fsd2d2}) is proved by applying
(\ref{err1090}) to a $\tilde \tau$ with
$L_{\tilde \tau} f = L_\tau f$, and showing that there exists $C'$ 
which only depends on $\Omega$ such that 
\begin{equation}
\label{adnfo8dsfbsdfdf8}
2^{-J} \|\tilde \tau \|_\infty + \|\nabla \tilde \tau \|_\infty (\log 
\frac{\|\Delta \tilde \tau\|_\infty}{\|\nabla \tilde \tau \|_\infty } \vee 1)
+ \|\rH \tilde \tau \|_\infty 
\leq C'\Bigl(\| \nabla \tau \|_\infty + \| \rH \tau \|_\infty 
\Bigr)~.
\end{equation}
Since $f$ has a support in $\Omega$, 
$L_{\tilde \tau} f = L_\tau f$ is equivalent to 
$\tilde \tau (x) = \tau (x)$ for all $x \in \Omega_\tau = 
\{x\,:\,x - \tau(x) \in \Omega\}$ and $\tilde \tau^{-1}(\Omega) = \Omega_\tau$.
If $\Omega$ has a radius $R$ then
the radius of $\Omega_\tau$
is smaller than $2R$, because $\|\nabla \tau \|_\infty \leq 1/2$.
We define $\tilde \tau$ as a regular extension of $\tau$ 
equal to $\tau(x)$ for $x \in \Omega_\tau$ and
to the constant $\min_{x \in \Omega_\tau} \tau(x)$
outside a compact $\widetilde \Omega_\tau$ of radius 
$(4R+2)$ including $\Omega_\tau$. It results that 
\begin{equation}
\label{adnfo8dsfbsdf}
{\|\Delta \tilde \tau\|_\infty} 
= \sup_{(x,u) \in \widetilde \Omega_\tau^2} |\tilde \tau (x) - \tilde \tau (u)| \leq
(4R+2)\,\|\nabla \tilde \tau \|_\infty~.
\end{equation}
The extension in $\widetilde \Omega_\tau - \Omega_\tau$
can be made regular in the sense that
$\|\nabla \tilde \tau\|_\infty + \|H \tilde \tau\|_\infty \leq 
\alpha\, (\|\nabla \tau\|_\infty \|+ \|H \tau\|_\infty) $
for some $\alpha > 0$ which depends on $\Omega$. This property 
together with (\ref{adnfo8dsfbsdf})
proves (\ref{adnfo8dsfbsdfdf8}). $\Box$

Similarly to Theorem \ref{DefomsProps},
if $\cP_J$ is replaced by the subset $\cP_{J,m}$ of 
paths of length smaller than $m$, then
$\|U[\cP_J] f \|_{1}$ is replaced by $m \, \|f\|$ 
in (\ref{err2fsd2d2}). If $L_{\tau} f(x) = f((1-s)\,x)$ with
$|\nabla \tau(x)| = |s| < 1$ then
the upper bound (\ref{err2fsd2d2})
is proportional to $m\,|s|\, \|f\|$. In this case, a lower bound
is simply obtained by observing that since $\|S_J[\cP_J] f\| = \|f\|$ 
and $\|S_J[\cP_J] L_\tau f\| = \|L_\tau f\| = (1-s)^{-1}\,\|f\|$ 
\[
\| S_J[\cP_J] L_\tau f -  S_J[\cP_J] f   \| 
\geq  |\| L_\tau f\| -  \| f \|  | 
> 2^{-1} \,s\, \|f\|. 
\]
Together with the upper bound (\ref{err2fsd2d2}), it
proves that if $\tau(x)= s x$ then the scattering distance of
$f$ and $L_\tau f$ is of the order of $\|\nabla \tau\|_\infty\, \|f\|$.

The next theorem reduces the translation error term $2^{-J} \|\tau\|_\infty$ in 
Theorem \ref{DefomsProps} to a second-order term $2^{-2J} \|\tau\|^2_\infty$, with 
first-order Taylor expansion of each $S_J [p] f$.
We denote $\nabla S_J[\cP_J] f(x) = \{ \nabla S_J [p] f(x) \}_{p \in \cP_J}$ and
$\tau(x) \cdot \nabla S_J[\cP_J] f(x) = \{ \tau(x) \cdot \nabla S_J [p] f (x)\}_{p \in \cP_J}$.

\begin{theorem}
\label{Theoelas2} 
There exists $C$ such that all $f \in \LD$
with $\|U[\cP_J] f \|_{1} < \infty$ and all $\tau \in \bCd$ with
$\|\nabla \tau\|_\infty \leq 1/2$ satisfy
\begin{equation}
\label{err1092}
\| S_J[\cP_J] L_\tau f -  S_J[\cP_J] f   
+ \tau \,\cdot\, \nabla  {S_J} [\cP_J] f\| \leq  C\,\|U[\cP_J] f \|_{1} \,
K (\tau)
\end{equation}
with
\begin{equation}
\label{err10920}
K (\tau) = 
2^{-2J} \|\tau \|^2_\infty + \|\nabla \tau \|_\infty 
(\log \frac{\| \Delta \tau\|_\infty}{\|\nabla \tau \|_\infty }\vee 1)
+ \|\rH \tau \|_\infty ~.
\end{equation}
\end{theorem}

{\it Proof:}
The proof proceeds as the proof of Theorem \ref{DefomsProps}.
Replacing 
$S_J[\cP_J] L_\tau  -  S_J[\cP_J]$ by 
$S_J[\cP_J]L_\tau  -  S_J[\cP_J]    + 
\tau \,.\, \nabla  {S_J}[\cP_J]$
in the derivation steps 
of the proof of Theorem \ref{DefomsProps} amounts
to replace $L_\tau A_J - A_J$ by 
$L_\tau A_J - A_J + \nabla A_J$. Equation (\ref{err10})
then becomes
\begin{eqnarray*}
\| S_J[\cP_J]L_\tau f -  S_J[\cP_J]f   + \tau \,.\, \nabla  {S_J} [\cP_J]
\| \leq  \|U[\cP_J] f \|_{1}&\Bigl(&
\|L_\tau A_J - A_J + \nabla A_J \| \\
& & + \|[\oWJ,L_\tau]\|\Bigr)~.
\end{eqnarray*}
Appendix \ref{proofLemma2} proves that
there exists $C > 0$ such that 
\begin{equation}
\label{err10dfi01}
\| L_\tau A_J f - A_J + \nabla A_J \| 
\leq C\,2^{-2J} \, {\|\tau \|^2_\infty}~.
\end{equation}
Inserting the upper bound (\ref{commnaf}) 
of $\|[\oWJ,L_\tau]\|$ proves (\ref{err1092}).
$\Box$

If $2^{J}  \gg \|\tau\|_\infty$ and $\|\nabla \tau \|_\infty 
+ \|\rH \tau \|_\infty \ll 1$ then 
$K(\tau)$ becomes negligible and
$\tau(x)$ can be estimated at each $x$ by solving 
the system of linear equations resulting from (\ref{err1092}):
\begin{equation}
\label{systenadfms}
\forall p \in \cP_J~~,~~ S_J [p] L_\tau f(x) -  
S_J [p] f(x) + \tau(x) \,.\, \nabla  S_J [p] f (x)
\approx 0.
\end{equation}
In dimension $d$, the displacement 
$\tau (x)$ has $d$ coordinates which can be computed
if the system (\ref{systenadfms}) has rank $d$. Estimating $\tau (x)$ 
has many applications. In image processing, the displacement field
$\tau (x)$ between two consecutive images 
of a video sequence is proportional to the optical flow velocity of 
image points.

\section{Normalized Scattering Transform}
\label{infdnsdfonsdsec}

To define the convergence of $S_J [\cP_J]$, all countable sets
$\cP_J$ are embedded in a non-countable set $\cPinf$.
Section \ref{measure} constructs a
measure $\mu$ and a metric in $\cPinf$. 
Section \ref{Renonscat} redefines the 
scattering transform $\oS f$
as limit of windowed scattering transforms
over $\cPinf$, with $\oS f \in \Ld (\cPinf,d\mu)$ for $f \in \LD$.
Numerical comparisons between $\oS f$ and $|\hat f|$ are given in
Section \ref{numerics}.

\subsection{Dirac Scattering Measure and Metric}
\label{measure}

A path $p \in \cP_J$ can be extended into an infinite set of paths 
in $\cP_{J\pl 1}$ which refine $p$. In that sense, $\cP_{J\pl 1}$ is
a set of higher resolution paths.
When $J$ increases to $\infty$,
these progressive extensions converge to paths of infinite length, which belong
to an uncountable path set $\cPinf$. 
A measure and a metric are defined on $\cPinf$.

A path $p = (\la_1,...,\la_m)$ of length $m$ belongs to 
the finite product set $\cLa_\infty^m$ with $\cLa_\infty = 2^\Z \times G^+$.
An infinite path $p$
is an infinite ordered string which belongs to 
the infinite product set $\cLa_\infty^\infty$. 
For complex-valued functions, adding negative paths
$(- \la_1,\la_2,...,\la_m)$ doubles the size of $\cLa_\infty^m$ 
and $\cLa_\infty^\infty$.
We concentrate on positive paths $(\la_1,\la_2,...,\la_m)$ and the
same construction applies to negative paths.
Since $\cLa_\infty = 2^\Z \times \G^+$ is
a discrete group,
its natural topology is the discrete topology where 
basic open sets are individual elements.
Open elements of the product topology of
$\cLa^\infty_\infty$ are cylinders 
defined for any $\la \in \cLa_\infty$ and $n \geq 0$ by
$C_n (\la) = \{q = \{ q_k \}_{k > 0} \in  \cLa^\infty_\infty~:~q_{n+1} = \la \}$
\cite{Willard}.
Cylinder sets are intersections of a finite number of open cylinders:
\[
C_n(\la_1,...,\la_m) = \{q \in \cLa^\infty_\infty~:~q_{n+1} = 
\la_1,...,q_{n+m} = \la_{m} \} = \cap_{i=1}^m C_{n+i} (\la_i)~.
\] 
As elements of the topology, cylinder sets are open sets but are also closed.
Indeed the complement of a cylinder set is a union of cylinders and is thus
closed. As a result, the topology is a sigma algebra, on which
a measure $\mu$ can be defined. The measure of a cylinder set $C$
is written $\mu(C)$. 

Let $\cP_\infty$ be the set of all finite paths including the $\nill$ path:
$\cP_\infty = \cup_{m \in \N} \cLa_\infty^m$. To any
$p = (\la_1,...,\la_m) \in \cP_\infty$, we associate a cylinder set:
\[
C(p) = C_0(p) = \{ q \in \cLa^\infty_\infty~:~q_1 = \la_1,...,q_m = \la_{m} \}~.
\]
This family of cylinder sets generates the same
sigma algebra as open cylinders since open cylinders can be written
$C_n (\la) = \cup_{(\la_1,...,\la_{n}) \in \cLa_\infty^n} C(\la_1,...,\la_n,\la)$.
The following proposition defines a measure on $\cLa^\infty_\infty$ 
from the scattering of a Dirac:
\[
U[p] \delta = |~||\psi_{\la_1}| \star \psi_{\la_2}|\star ...|\star \psi_{\la_m}|~.
\]

\begin{proposition}
\label{Disfdnsdfj}
There exists a unique 
$\sigma$-finite
Borel measure $\mu$, called Dirac scattering measure, 
such that $\mu(C(p)) = \|U[p] \delta \|^2$ for all $p \in \cP_\infty$.\\
For all $2^l g \in \Lambda_\infty$ and $p \in \cP_\infty$, 
$\mu(C(2^l g p)) = 2^{d l} \mu(C(p))$\\
If $|\hat \psi(\om)|+|\hat \psi(-\om)| \neq 0$ 
almost everywhere then
$\|U[p] \delta\| \neq 0$ for $p \in \cP_\infty$.
\end{proposition}

{\it Proof:} The Dirac scattering 
measure is defined as
a subdivision measure over the tree that generates all paths.
Each finite path $p$ corresponds to a node of the subdivision tree.
Its sons are the $\{p + \la \}_{\la \in \cLa_\infty}$,
and $C(p) = \cup_{\la \in \cLa_\infty} C(p+\la)$ is a non-intersecting partition. 
Since
\[
\|U[p] \delta \|^2 = \|\oUJ U[p] \delta \|^2 = \sum_{\la \in \cLa_\infty} \|U[p + \la] \delta \|^2,
\]
it results that $\mu(C(p)) = \sum_{\la \in \cLa_\infty} \mu(C(p+\la))$.
The sigma additivity of the Dirac measure over all cylinder sets results from
the tree structure, and the decomposition of the measure of
a node $\mu(C(p))$ as a sum of the measures $\mu(C(p+\la))$ of all its sons.
This subdivision measure is
uniquely extended to the Borel sigma algebra through
the sigma additivity. 
Since $\cLainf = \cup_{\la \in \cLa_\infty} C(\la)$ and 
$\mu(C(\la)) = \|U[\la] \delta \|^2 = \|\psi_{\la} \|^2$,
this measure is $\sigma$-finite.

We showed in (\ref{pathcomansdfus29df}) that
$U[p] (2^l g \circ f) = 2^l g \circ U[2^{-l} g p] f$.
Since $2^l g \circ \delta = 2^{-dl} \delta$ it results
$\|U[2^{-l} g p] \delta \|^2 = 2^{-dl} \|U[p] \delta \|^2$ and hence
$\mu(C(2^l g p)) = 2^{d l} \mu(C(p))$.

If the set of $\omega \in \R^d$ where $\hat \psi(\om)=0$ and $\hat \psi(-\om)=0$
is of measure $0$, let us prove by induction on the path
length that $U[p] f \neq 0$ if $f \in \LD \cup \LU$ or if $f = \delta$. 
We suppose that $U[p] f \neq 0$ and verify that $U[p+\la] f \neq 0$ for 
any $\la \in \cLa_\infty$
Since $U[p] f$ is real,
$|\widehat {U[p] f} (\om)| = |\widehat {U[p] f} (-\om)|$. But
$\hat \psi_{\la} (\om) = \hat \psi(\la^{-1} \om)$, so 
$\hat \psi_\la (\om)$ and $\hat \psi_\la (-\om)$ vanish simultaneously on a set
of measure $0$. It results that $\widehat {U[p+\la]} f = \widehat {U[p] f} \, 
\hat \psi_\la \neq 0$ if $\widehat {U[p] f} \neq 0$ so $U[p+\la] f$ is a non-zero
function. $\Box$

A topology and a metric can now be 
constructed on the path set $\La^\infty_\infty$.
Neighborhoods are defined with cylinder sets of frequency resolution $2^{J}$:
\begin{equation}
\label{dentadfnsdf}
C_{J} (p) = \cup_
{\la \in \La_\infty \atop |\la| \leq 2^{-J}}
C(p+\la) \subset C(p)~.
\end{equation}
Clearly  $C_{J+1}(p) \subset C_{J} (p)$. The following proposition 
proves that $\mu(C_J (p))$
decreases at least like $2^{-dJ}$ when $2^J$ increases, and it 
defines a distance from these measures.
The set $\La^\infty_\infty$ of infinite paths
is not complete with this metric. It is completed by embedding the
set $\cP_\infty$ of finite paths, and we denote
$\cPinf = \cP_\infty \cup \cLainf$ the completed set. This embedding
is defined by adding each finite path $p \in \cP_\infty$ 
to $C(p)$ and to each $C_J (p)$ for all $J \in \Z$, without modifying
their measure. 
We still denote $C_J (p)$ the resulting subsets of $\cPinf$.
For complex valued functions, the size of $\cPinf$ 
is doubled by adding finite and infinite negative paths 
$(-\la_1,\la_2,...,\la_m,...)$.

\begin{proposition}
\label{propodist}
If $p \in \cP_\infty$ is a path of length $m$ then 
\begin{equation}
\label{sizeneighboh}
\mu(C_J(p)) = \|S_J\delta[p] \|^2 \leq  2^{-dJ}\,\|\phi\|^2\,\|\psi\|_1^{2m}~.
\end{equation}
Suppose that $|\hat \psi(\om)|+|\hat \psi(-\om)| \neq0$ 
almost everywhere. For any $q \neq q' \in \cPinf$ 
\begin{equation}
\label{disdnfs0d98fns}
\bar d(q,q') = \inf_{(q,q') \in C_J(p)^2} \mu(C_J(p)) ~~\mbox{and}~~\bar d(q,q) = 0
\end{equation}
defines a distance on $\cPinf$, and $\cPinf$ is complete for this
metric.
\end{proposition}

{\it Proof:}
According to (\ref{dentadfnsdf})
\[
\mu(C_J(p)) = \sum_{\la \in \La_\infty \atop |\la| \leq 2^{-J}} \mu(C(p+ \la)) = 
\sum_{\la \in \La_\infty \atop |\la| \leq 2^{-J}} 
\|U[p+\la] \delta \|^2  .
\]
Since $U[p+\la] \delta  = U[p] \delta  \star \psi_{\la}$ and
$|\hat \phi_{2^{J}} ( \om)|^2 = 
\sum_{\la \in \La_\infty \atop |\la| \leq 2^{-J}} 
|\hat \psi_{\la} (\omega)|^2$, the Plancherel formula implies
\[
\mu(C_J(p)) = \sum_{\la \in \La_\infty \atop |\la| \leq 2^{-J}} 
\|U[p] \delta  \star \psi_{\la}\|^2 =
\|U[p] \delta \star \phi_{2^J} \|^2 = \|S_J[p] \delta \|^2~.
\]
Since $S_J[p] \delta = U[p] \delta \star \phi_{2^J}$, 
Young's inequality implies
$\|S_J[p] \delta\| \leq \|U[p] \delta\|_1 \, \|\phi_{2^J}\|$.
Moreover $\|U[\la] f \|_1 \leq \|\psi_\la\|_1 \|f\|_1$ with $\|\psi_\la \|_1 = \|\psi \|_1$, so we verify by induction that
$\|U[p] \delta \|_1 \leq \|\psi \|^m$. Inserting
$\|\phi_{2^J}\|^2 = 2^{-dJ} \|\phi \|^2$ proves
(\ref{sizeneighboh}).

Let us now prove that $\bar d$ defines a distance. 
If $q \neq q'$, we denote $\bar p \in \cP_\infty$ their common prefix of
longest size $m$, which may be $0$, and show that $\bar d(q,q') \neq 0$. 
Let $|q_{m+1}| = 2^{j_{m+1}}$ and $|q'_{m+1}| = 2^{j'_{m+1}}$ be the frequencies
of their first different coordinate.
If $2^{-J} = \max(|q_{m+1}|, |q'_{m+1}|)$ then $(q,q') \in C_J(\bar p)^2$ and 
it is the smallest set including both paths so 
$\bar d(q,q') = \mu(C_J (\bar p))$.
It results that $\bar d(q,q') \neq 0$ because 
$\mu(C_J (\bar p)) \geq \mu(C (\bar p + 2^J r))$ for $r \in \G^+$ and
Proposition \ref{Disfdnsdfj} proves that $\mu(C(p)) \neq 0$ 
for all $p \in \cP_\infty$, so $\bar d(q,q') \neq 0$.

The triangle inequality is proved by showing that
\begin{equation}
\label{inestrian}
\forall (q,q',q'') \in \overline {\mathcal P}_\infty^3~~,~~
\bar d (q',q'') \leq \max \Bigl(\bar d (q,q')\,,\,\bar d (q,q'')\Bigr)~.
\end{equation}
This is verified by writing
$\bar d(q,q') = \mu(C_J (\bar p))$, $\bar d(q',q'') = \mu(C_{J'} (\bar p'))$
and  $\bar d(q',q'') = \mu(C_{J''} (\bar p''))$.
Necessarily $\bar p$ is a substring of $\bar p'$ or vice versa, and
$\bar p''$ is larger then the smallest of the two.
If $\bar p''$ is strictly larger then the smallest say $\bar p$, 
then $\mu(C_{J''} (\bar p'')) \leq \mu(C (\bar p'')) \leq C_J (\bar p)$,
so (\ref{inestrian}) is satisfied. If $\bar p'' = \bar p = \bar p'$
then $2^{-J''} \leq \max(2^{-J}, 2^{-J'})$ and (\ref{inestrian}) is satisfied.
Otherwise $\bar p'' = \bar p$ is strictly smaller than $\bar p'$ and
necessarily $2^{J''} = 2^J$ so (\ref{inestrian}) is also satisfied.

To prove that $\cPinf$ is complete, consider 
a Cauchy sequence $\{q_j \}_{j \in \N}$ in $\cPinf$.
Let $p_k$ be the common prefix of maximum length $m_k$
among all $q_j$ for $j \geq k$. It is a growing string which
either converges to a finite string $p \in \cP_\infty$ if $m_k$ is bounded
or to an infinite string $p \in \cLainf$.
Among all paths
$\{q_j \}_{j \geq k}$ whose maximum common prefix with $p$ has a length $m_k$,
let $q_{j_k}$ be a path whose next element $\la_{m_k+1}$
has a maximum frequency amplitude $|\la_{m_k+1}|$.
One can verify that
\[
\sup_{j,j' \geq k} \bar d (q_j,q_{j'}) = \bar d(q_{j_k},p) = 
\sup_{j \geq k} \bar d (p,q_j)~.
\]
Since $\sup_{j,j' \geq k} \bar d (q_j,q_{j'})$ 
converges to $0$ as $k$
increases, it implies that
$\sup_{j \geq k} \bar d (p,q_j)$ also converges to $0$ and hence 
that $\{q_j \}_{j \in \N}$ converges to $p$.
$\Box$

\subsection{Scattering  Convergence}
\label{Renonscat}

For $h \in \Ld(\cPinf,d\mu)$, we
denote $\|h\|^2_\cPinf = \int_\cPinf |h(q)|^2\,d\mu(q)$,
where $\mu$ is the Dirac scattering measure.
This section redefines the scattering transform
$ \oS f$ as a limit of windowed scattering
transforms, and proves that $\oS f \in \Ld(\cPinf,d\mu)$
for all $f \in \LD$. 
We suppose that $\psi$ is an admissible scattering wavelet,
and that
$|\hat \psi(\om)|+|\hat \psi(-\om)| \neq 0$ almost everywhere.

Let $1_{C_J(p)} (q)$ be the indicator function of $C_J(p)$
in $\cPinf$.
A windowed wavelet
scattering  $S_J[\cP_J] f(x) = \{S_J[p] f(x)\}_{p \in \cP_J}$ 
is first extended into a normalized function of $(q,x) \in \cPinf \times \R^d$ 
\begin{equation}
\label{joinsdfpnas}
S_J f(q,x) = \sum_{p \in \cP_J} 
\frac {S_J[p] f(x)} {\|S_J [p] \delta \|} \, {1_{C_J(p)}(q)}~.
\end{equation}
It satisfies
$S_Jf(p,x) = S_J [p] f(x)/\|S_J [p] \delta \|$ for $p \in \cP_\infty$.
Since $\mu(C_J(p)) = \|S_J [p] \delta \|^2$,
for all $(f,h) \in \LD^2$
\[
\int_{\cPinf} \int_{\R^d} |S_J f(q,x) - S_J h(q,x)|^2\, d\mu(q)\,dx = 
\|S_J[\cP_J] f - S_J[\cP_J] h\|^2 \leq \|f - h \|^2~,
\]
\[
\int_{\cPinf} \int_{\R^d} |S_J f(q,x)|^2\, d\mu(q)\,dx = 
\|S_J[\cP_J] f\|^2 = \|f\|^2~,
\]
so $S_J f(q,x)$ can be interpreted as a
scattering energy density in $\cPinf \times \R^d$.

The windowed scattering $S_J f(q,x)$ 
has a spatial resolution $2^{-J}$ along $x$ 
and a resolution $2^{J}$ along the frequency path $q$.
When $J$ goes to $\infty$, $S_J f(q,x)$ loses its spatial localization, and
Theorem \ref{theoInva} proves that the asymptotic metric on $S_J[\cP_J]f$ 
and hence on $S_J f(q,x)$ is
translation invariant. The convergence of $S_J f(q,x)$ to 
a function which depends only on $q \in \cPinf$ 
is studied by introducing the marginal 
$\LD$ norm of $S_J f(q,x)$ along $x$ for $q$ fixed:
\begin{equation}
\label{indfods08hwef0}
\forall q \in \cPinf~~,~~\overline S_J f(q) = \int |S_J f(q,x)|^2\,dx 
= \sum_{p \in \cP_J} \frac{\| S_J [p]f\|} {\|S_J [p] \delta \|}\,
1_{C_J (p)} (q)\,.
\end{equation}
It is a piecewise constant function of the path variable $q$,
whose resolution increases with $J$.
Since $\mu(C_J(p)) = \|S_J [p] \delta \|^2$,
\begin{equation}
\label{indfods08hwef0d}
\|\oS_J f - \oS_J h \|^2_\cPinf = \int_\cPinf |\oS_J f(q) - \oS_J h(q) |^2\,d\mu(q) 
= \sum_{p \in \cP_J} \Big| \|S_J[p] f\| - \|S_J[p] h\| \Big|^2~.
\end{equation}
The following proposition proves that $\oS_J$ is a nonexpansive operator
which preserves the norm.

\begin{proposition}
\label{bornsdfnsdf}
For all $(f,h) \in \LD^2$ and $J \in \Z$
\begin{equation}
\label{compardisng30}
\|\oS_J f - \oS_J h \|_\cPinf \leq \|\oS_{J+1} f - \oS_{J+1} h \|_\cPinf~,
\end{equation}
\begin{equation}
\label{compardisng3}
\|\oS_J f - \oS_J h \|_\cPinf \leq 
\|S_J[\cP_J] f - S_J[\cP_J] h \|  \leq \|f - h\|~,
\end{equation}
\begin{equation}
\label{indfods08hwef}
\|\oS_J f\|_\cPinf = \|f\|~.
\end{equation}
\end{proposition}

{\it Proof:}  We proved in (\ref{pathconsdf9u8sd2}) that
\begin{equation}
\label{pathconsdf9u8sd222}
\|S_J[p] f\|^2  = \sum_{p' \in \cP^p_{J+1}} \|S_{J\pl 1}[p'] f\|^2 , 
\end{equation}
where $\cP_{J\pl 1} = \cup_{p \in \cP_J} \cP^p_{J\pl 1}$ is a disjoint partition.
Applying this to $f$ and $h$ implies 
\[
\Big| \|S_J[p] f\| - \|S_J[p] h\| \Big|^2 \leq
\sum_{p' \in \cP^p_{J\pl 1}} \Big| \|S_{J\pl 1}[p'] f\| - \|S_{J\pl 1}[p'] h\| \Big|^2 .
\]
Summing over $p \in \cP_{J}$ and inserting (\ref{indfods08hwef0d})
proves (\ref{compardisng30}).

Since $\Big|\|S_J[p] f\| - \|S_J[p] h \| \Big| \leq 
\|S_J[p] f - S_J[p] h \|$, 
summing this inequality over $p \in \cP_{J}$ and inserting (\ref{indfods08hwef0d})
proves the first inequality of (\ref{compardisng3}).
The second inequality is obtained because $S_J[\cP_J]$ is nonexpansive.
Setting $h = 0$ proves that $\|\oS_J \|_\cPinf = \|S_J[\cP_J] f\|$ and
Theorem \ref{energydecth} proves $\|S_J[\cP_J] f\| = \|f\|$, which gives
(\ref{indfods08hwef}).
$\Box$

Since $\|\oS_J f - \oS_J h \|_\cPinf$ is non-decreasing and bounded
when $J$ increases,
it converges to a limit which is smaller than the limit of the non-increasing
sequence $\|S_J[\cP_J] f - S_J [\cP_J] h \|$.
The following proposition proves that $\oS_J f$ converges pointwise
to the scattering transform on $\cP_\infty$ introduced in
Definition \ref{definwindoscas}.

\begin{proposition}
If $f \in \LU$ then 
\begin{equation}
\label{limsindfods08hwef}
\forall p \in \cP_\infty~~,~~
\lim_{J \rightarrow \infty} \oS_J f(p) = \oS f(p) = 
\frac 1 {\mu_p}~\int U[p] f (x)\,dx
\end{equation}
with $\mu_p = \int U[p] \delta (x)\,dx$.
\end{proposition}

{\it Proof:}
If $p \in \cP_\infty$ then for $J$ sufficiently large
$\oS_J f(p) = {\| S_J[p]f\|}/ {\|S_J[p] \delta \|}$. 
Let us prove that 
\begin{equation}
\label{firsndfinsdf2}
\lim_{J \rightarrow \infty} 2^{dJ/2} \| S_J[p] f \|= 
\|\phi\|\,\int U[p] f (x)\, dx~,
\end{equation}
and that this equality also holds for $f = \delta$.
Since $S_J[p] f = U[p]f \star \phi_{2^J}$, 
the Plancherel formula implies
\begin{equation}
\label{firsndfinsdf3}
2^{dJ}\,\| S_J[p] f \|^2 = 2^{dJ} \,(2 \pi)^{-d}\,\int |\widehat {U[p]f} (\omega)|^2\,
|\hat \phi(2^{J} \omega)|^2\, d\omega~.
\end{equation}
Since derivatives of $\phi$ are in $\LU$, we have
$\hat \phi(\omega) = O((1+|\omega|)^{-1})$
and hence $(2 \pi)^{-d} 2^{dJ}\,|\hat \phi(2^{J} \omega)|^2$ converges to
$\|\phi\|^2\, \delta(\omega)$. Moreover, if $f \in \LU$ then 
$U[p] f \in \LU$ so 
$\widehat {U[p]f} (\omega)$ is continuous at $\omega = 0$.
It results from (\ref{firsndfinsdf3}) that
$\lim_{J \rightarrow \infty} 2^{dJ} \| S_J[p] f \|^2= |\widehat {U[p]f} (0)|^2\,
\|\phi\|^2$ which proves (\ref{firsndfinsdf2}). The same 
derivations hold to prove this result for $f = \delta$.

Since $|\hat \psi(\om)| + |\hat \psi(-\om)| \neq 0$ almost everywhere,
Proposition \ref{Disfdnsdfj} proves that $U[p] \delta \neq 0$.
Since it is positive, it has a non-zero integral.
It results from (\ref{firsndfinsdf2}) that
$\lim_{J \rightarrow \infty} {\| S_J[p]f\|}/ {\|S_J[p] \delta \|} =
\int U[p] f(x) dx / \int U[p] \delta (x) dx$ which proves
(\ref{limsindfods08hwef}). $\Box$

The scattering transform $\oS f$ can now
be extended to $\cPinf$ as a windowed scattering limit:
\[
\forall q \in \cPinf~~,~~\oS f(q) = \liminf_{J \rightarrow \infty} \oS_J f(q)~.
\]
Proposition \ref{bornsdfnsdf} proves that
$\|\oS_J f \|_\cPinf = \|f\|$ so
Fatou's lemma  implies that
$\oS f \in \Ld(\cPinf,d \mu)$.
The following theorem gives a sufficient condition so that
$\oS_J f$ converges strongly to $\oS f$, which then preserves
the $\LD$ norm of $f$.

\begin{theorem}
\label{proconsdf}
If for $f \in \LD$ 
there exists $ \Omega_J^f \subset \cP_J$ with
\begin{equation}
\label{convadfndsoihfdsf}
\lim_{J\rightarrow  \infty}   \|S_J [\Omega_J^f] f\|^2 = 0
~~\mbox{and}~~
\lim_{J\rightarrow  \infty} 
\sup_{p \in \cP_J -  \Omega_J^f} \left\|\frac {S_J [p] f}  {\|S_J [p] f\|} - 
\frac {S_J [p] \delta}  {\|S_J [p] \delta\|} \right\| = 0
\end{equation}
then $\oS_J f$ converges in norm to $\oS f$ with 
$\|\oS f\|_\cPinf = \|f\|$  and
\begin{equation}
\label{sdfn089s}
\forall p \in \cP_\infty~~,~~
\int_{C(p)} | \oS_J f (q)|^2\, d \mu(q) = \|U[p] f\|^2~.
\end{equation} 
If $(f,h) \in \LD^2$ satisfy (\ref{convadfndsoihfdsf}) then
\begin{equation}
\label{convadfndsoihfdsf9}
\lim_{J\rightarrow \infty} \|S_J[\cP_J] f - S_J[\cP_J] h \| =
\| \oS  f -  \oS  h \|_\cPinf~.
\end{equation}
If (\ref{convadfndsoihfdsf}) is satisfied
in a dense subset of $\LD$ then 
$S_J f$ converges strongly to $\oS f$ for all $f \in \LD$ and
both (\ref{sdfn089s}) and (\ref{convadfndsoihfdsf9}) are satisfied
in $\LD$.
\end{theorem}

{\it Proof:}
The following lemma proves that $ \{ \oS_J f \}_{J \in \N}$ is 
Cauchy and hence converges in norm 
to $\oS f \in \Ld(\cPinf,d\mu)$.
The proof is in Appendix \ref{Cauchsdfproof}. 

\begin{lemma}
\label{Cauchsdf}
If $f \in \LD$ satisfies 
(\ref{convadfndsoihfdsf}) then
$\{\oS_J f \}_{J \in \N}$ is a Cauchy sequence in $\Ld(\cPinf,d\mu)$.
\end{lemma}

Since $\Ld(\cPinf,d\mu)$ is complete,
$\oS_J f(q)$ converges in norm to its limit inf $\oS f$.
Since $\|\oS_J f \| = \|f\|$, it also implies that $\|\oS f \|_\cPinf = \|f\|$.
Moreover, $U[p+q] = U[q] U[p]$ so
$\|\oS_J U[p] f \|_\cPinf^2 = \int_{C(p)} |\oS_J f (q)|^2 \, d\mu(q)$.
Since $\|\oS_J U[p] f \|_\cPinf^2 = \|U[p] f \|^2$ taking the limit when $J$
goes to $\infty$ proves  (\ref{sdfn089s}).

The windowed scattering 
convergence (\ref{convadfndsoihfdsf9}) relies on the
following lemma.

\begin{lemma}
\label{Rensdofnsdf}
If $(f,h) \in \LD^2$ satisfy 
(\ref{convadfndsoihfdsf}) then
\begin{equation}
\label{convadfndsoihf90}
\lim_{J\rightarrow \infty} \|S_J[\cP_J] f - S_J[\cP_J] h \| = 
\lim_{J\rightarrow \infty} 
\|\oS_J f - \oS_J h \|_\cPinf~.
\end{equation}
\end{lemma}

Since (\ref{convadfndsoihfdsf}) implies that 
$\oS_J f$ and $\oS_J h$ respectively converge in norm to $\oS f$ and $\oS h$,
the convergence (\ref{convadfndsoihfdsf9}) results from 
(\ref{convadfndsoihf90}). Proving 
(\ref{convadfndsoihf90}) is equivalent to proving that
$\lim_{J \rightarrow \infty} \sum_{p \in \cP_J} I_J(f,h)[p] = 0$ for
\[
I_J(f,h)[p] = \|S_J[p] f - S_J[p] h \|^2 - \Big|\|S_J[p] f\| - \|S_J[p] h\|
\Big|^2~.
\]
Observe that
\begin{eqnarray}
\nonumber
I_J (f,h)[p] &=& \|S_J[p] f\|\, \|S_J[p] h\| 
\left\|
\frac {S_J[p] f} {\|S_J[p] f\|} - \frac {S_J[p] h} {\|S_J[p] h\|} \right\|^2.\\
\label{089dfsndsofi3e}
&\leq & 2 \|S_J[p] f\|\, \|S_J[p] h\| \left(
\left\|\frac {S_J[p] f} {\|S_J[p] f\|} - \frac {S_J[p] \delta} {\|S_J[p] \delta\|} 
\right\|^2 
\right. \\
\nonumber
&&~~~~~~~~~~~~~~~~~~~~~~~~~+ \left.
\left\|\frac {S_J[p] h} {\|S_J[p] h\|} - \frac {S_J[p] \delta} {\|S_J[p] \delta\|} 
\right\|^2 \right). 
\end{eqnarray}
When summing over $p \in \cP_J$, we separate
$ \Omega_J^f \cup  \Omega^h_J$ from its complement in $\cP_J$. 
Since $\lim_{J\rightarrow \infty} \|{S_J [ \Omega_J^f ] f}\|^2 = 0$,
$\|{S_J [ \cP_J ] f}\|^2 = \|f\|^2$,
$\lim_{J\rightarrow \infty} \|{S_J [ \Omega_J^h] h}\|^2 = 0$,
and $\|{S_J [ \cP_J ] h}\|^2 = \|h\|^2$, dividing the
sum over $\Omega_J^f $ and $\Omega_J^h $ and applying
Cauchy-Schwartz proves that
\[
\lim_{J \rightarrow \infty} \sum_{p \in  \Omega_J^f \cup  \Omega_J^h}
\|S_J[p] f\|\, \|S_J[p] h\| = 0~,
\]
and $ \sum_{p \in \cP_J}
\|S_J[p] f\|\, \|S_J[p] h\| \leq \|f\|\, \|h\|$.
The hypothesis (\ref{convadfndsoihfdsf}) applied to $f$ and $h$ gives
\[
\lim_{J \rightarrow \infty} \sup_{p \in \cP_J -  \Omega_J^f \cup  \Omega_J^h}
\left(
 \left\|\frac {S_J[p] f} {\|S_J[p] f\|} - \frac {S_J[p] \delta} {\|S_J[p] \delta\|} 
\right\|^2 +
\left\|\frac {S_J[p] h} {\|S_J[p] h\|} - \frac {S_J[p] \delta} {\|S_J[p] \delta\|} 
\right\|^2 \right) = 0 ~
\]
so (\ref{089dfsndsofi3e})
implies that $\lim_{J \rightarrow \infty} \sum_{p \in \cP_J} I_J(f,h)[p] = 0$, which
finishes the Lemma proof.

Suppose that (\ref{convadfndsoihfdsf}) is satisfied
in a dense subset of $\LD$. 
Any $f \in \LD$ is the
limit of $\{f_n\}_{n > 0}$ in this dense set. Since $\oS$
and $\oS_J$ are nonexpansive
\[
\|\oS f - \oS_J f \|_\cPinf \leq 2\, \|f - f_n \| + \|\oS f_n - 
\oS_J f_n \|_\cPinf~.
\]
Since $f_n$ satisfies (\ref{convadfndsoihfdsf}), 
we proved that $\oS_J f_n$ converges in norm to $\oS f_n$. Letting $n$ go
to $\infty$ implies that $\oS_J f$ converges in norm to $\oS$. 
The previous derivations then implies that
both (\ref{sdfn089s}) and (\ref{convadfndsoihfdsf9}) are satisfied 
in $\LD$. $\Box$

If $f \in \LU$ and $p \in \cP_\infty$,
since $S_J [p] f(x) = U [p] f \star \phi_{2^J}$ and $\|U [p] f \|_1 < \infty$,
applying the Plancherel formula proves that
\begin{equation}
\label{asdf089dfwefs}
\lim_{J\rightarrow \infty} 
 \Big\|\frac {S_J[p] f}  {\|S_J[p] f\|} -
\frac{S_J [p] \delta } {\|S_J [p]\delta  \|}  \Big\|^2 = 0~.
\end{equation}
This is however not sufficient to prove (\ref{convadfndsoihfdsf}) because
the sup is taken over all $p \in \cP_J - \Omega_J^f$ which grows when $J$ 
increases. For $f \in \LU$, one can find paths $p_J \in \cP_J$, 
which are not frequency-decreasing, 
where ${S_J[p_J] f}/  {\|S_J[p_J] f\|}$ does not converge
to ${S_J[p_J] \delta}/  {\|S_J[p_J] \delta\|}$. The main difficulty is
to prove that over the set $\Omega_J^f$ of all such paths, a
windowed scattering transform has a norm
$ \|S_J[\Omega_J^f] f\|$ which converges to zero.
Numerical experiments indicate that this property could be valid for all
$f \in \LU$. It also seems that if $f \in \LU$ then
$\oS f(q)$ is a continuous function of the path
$q$, relatively to the Dirac scattering metric.
This is analogous to the Fourier transform continuity when $f \in \LU$.

\begin{conjecture} 
\label{conjectsdf}
Condition (\ref{convadfndsoihfdsf}) 
holds for all $f \in \LU$. 
Moreover, if $f \in \LU$ then $\oS f(q)$ is
continuous in $\cPinf$ relatively to the Dirac scattering metric.
\end{conjecture}

If this conjecture is valid, since $\LU$ is dense in $\LD$, then
Theorem \ref{proconsdf} proves
that $\oS_J$ converges strongly to
$ \oS f$ for all $f \in \LD$, and $\|\oS f \|_\cPinf = \|f\|$.
Property (\ref{convadfndsoihfdsf9}) also proves that 
$\|S_J [\cP_J] f - S_J [\cP_J] h\|$ converges to
$\|\oS f - \oS h\|_\cPinf$ as $J$ goes to $\infty$. 
Through this limit,
the Lipschitz continuity of 
$S_J$ under the action of diffeomorphisms can then be extended to 
the scattering transform $\oS$.

\subsection{Numerical Comparisons with Fourier}
\label{numerics}

Let $\R^{d+}$ be the half frequency space of all $\om = (\om_1,...,\om_d) \in \R^d$ with
$\om_1 \geq 0$ and $\om_k \in \R$ for $k > 1$.  
To display numerical examples for real functions, 
the following proposition 
constructs a function from
$\R^{d+}$ to $\cPinf$ 
which maps the Lebesgue measure of $\R^{d+}$
into the Dirac scattering measure.
It provides a representation of $\oS f$ over  $\R^{d+}$.
We assume that $\psi$ is an admissible scattering wavelet,
and that
$|\hat \psi(\om)|+|\hat \psi(-\om)| \neq 0$ almost everywhere.

\begin{proposition}
\label{isodfmpsdfon}
There exists a surjective function $q(\omega)$
from $\R^{d+}$ onto $\cPinf$ 
such that for all measurable sets $\Omega \subset \cPinf$
\begin{equation}
\label{additnasfd}
\mu(\Omega) = \int_{q^{-1}(\Omega)}  d \omega ~. 
\end{equation}
\end{proposition}

{\it Proof:} The proof first
constructs the inverse $q^{-1}$ 
by mapping each cylinder $C(p)$ for 
$p \in \cP_\infty$ into 
a set $q^{-1}(C(p)) \subset \R^{d+}$ satisfying the following properties:
$\mu(C(p)) = \int_{q^{-1}(C(p))} d \om$, and
$q^{-1}(C(p)) \cap q^{-1}(C(p')) = \emptyset$ if
$C(p) \cap C(p') = \emptyset$, and
$q^{-1}(C(p)) \subset q^{-1}(C(p'))$ if $C(p) \subset C(p')$. 
Let $\overline{q^{-1}(C(p))}$ be the closure of $q^{-1}(C(p))$ in $\R^{d+}$. 
For all $p \neq \emptyset$, we also impose that 
the frontier of $q^{-1}(C(p))$ is a set
of measure $0$ in $\R^{d+}$, and that
$\overline {q^{-1}(C(p+\la))} \subset q^{-1}(C(p))$ for all $\la \in \cLa_\infty$.
The cylinders $C(p)$ generate the sigma algebra on which
the measure $\mu$ is defined. A measurable set $\Omega$ can be
approximated by sets $\Omega_k$ which 
are union of disjoint cylinder sets $C(p)$ with
$\lim_{k \rightarrow \infty} \mu(\Omega-\Omega_k) = 0$.
The properties of $q^{-1}$ on the cylinders $C(p)$ imply that 
$\int_{q^{-1} (\Omega_k)} d\om = \mu(\Omega_k)$ and when $k$ goes to $\infty$
we get (\ref{additnasfd}).

Once all $q^{-1}(C(p))$ are constructed, the inverse
$q(\om)$ is uniquely defined for all $\om \in \R^{d+}$, as follow.
Let $p_m$ be the prefix of
$\bar q \in \cPinf$ of length $m$. We define 
$q^{-1}(\bar q) = \cap_{m \in \N} q^{-1}(C(p_m))$.
Since $\overline {q^{-1}(C(p+\la))} \subset q^{-1}(C(p))$ for all $\la \in \cLa_\infty$,
it results that 
$\cap_{m \in \N} q^{-1}(C(p_m)) = \cap_{m \in \N} \overline{q^{-1}(C(p_m))}$.
It is a closed non-empty set
because $\overline{q^{-1}(C(p_m))} \subset 
\overline{q^{-1}(C(p_{m-1}))}$ is a 
non-empty set of measure $\|U[p_m] \delta \| \neq 0$.
We verify that
$q(\om) = \bar q$ for all $\om \in q^{-1}(\bar q)$
defines a surjective function on $\R^{d+}$ by showing that 
$\cup_{\bar q \in \cPinf} q^{-1}(\bar q)$ is a partition
of $\R^{d+}$.
If $\cP_m$ is the set of all path of length $m$ then
$\cup_{p \in \cP_m} C(p)$ is a partition of $\cPinf$, so
the recursive construction of $q^{-1}$ implies that
$\cup_{p \in \cP_m} q^{-1}(C(p))$ is a partition of 
$\R^{d+}$. Letting $m$ go to infinity proves
that $\cup_{\bar q \in \cPinf} q^{-1}(\bar q)$ is a partition
of $\R^{d+}$.

The sets $q^{-1}(C(p))$ satisfying the previously mentioned properties
are defined recursively on the path length, 
with a subdivision procedure.
In dimension $d = 1$, 
each $q^{-1}(C(p))$ is recursively defined as an interval of $\R^+$.
We begin with paths $p = 2^j$ of length $1$
by defining $q^{-1}(C(2^j )) = [2^{j} \|\psi\|^2 , 2^{j+1} \|\psi\|^2)$,
whose width is $2^{j} \|\psi\|^2  = \mu(C(2^j))$.
Suppose now that $q^{-1}(C(p))$ is an interval of width equal to $\mu(C(p))$.
All $q^{-1}(C(p+2^j))$ for $j \in \Z$ are 
defined as consecutive intervals $[a_j,a_{j-1})$, which
define a partition of $q^{-1}(C(p)) = \cup_{j \in \Z} [a_j,a_{j-1})$ with
$a_{j-1} - a_j = \|U[p+2^j]\delta \|^2 = \mu(C(p+2^j))$.
One can verify that this recursive construction defines
intervals $q^{-1}(C(p))$ which satisfy all mentioned properties.
Moreover, in this case the resulting function 
$q(\om)$ is bijective from $\R^+$ to $\cPinf$.

In higher dimensions $d \geq 1$, this construction is extended as follow.
All cylinders $C(\la)$ for all paths $p = \la = 2^j r$ of length $1$ 
are mapped to non-intersecting 
hyper-rectangles $q^{-1}(C(\la))$ of measure
\[
\int_{q^{-1}(C(2^j r))} d\om = \mu(C(2^j r)) = \|U[2^j r] \delta \|^2 =  
2^{dj}\,\|\psi\|^2 ~.
\]
These hyper-rectangles are chosen to 
define a partition of $\R^{d+}$, and hence $\R^{d+}
 = \cup_{\la \in \cLa_\infty } q^{-1} (C(\la))$ with 
$q^{-1} (C(\la)) \cap q^{-1} (C(\la')) = \emptyset$ for $\la \neq \la'$.
Suppose now that $q^{-1}(C(p))$, with $\int_{q^{-1}(C(p))} dx = \|U[p] \delta \|^2$,
is defined for all paths $p$ of length $m$. 
Since $U$ preserves the norm
$\sum_{\la \in \La_\infty} \|U[p+\la] \delta \|^2 = \|U[p] \delta \|^2$.
We can thus partition $q^{-1}(C(p))$ into subsets
$\{q^{-1}(C(p+\la)) \}_{\la \in \cLa_\infty}$ with
$\int_{q^{-1}(C(p+\la)} d \omega = \|U [p+\la] \delta\|^2$,
whose frontiers are piecewise hyperplanes of dimensions $d-1$ and hence
have a zero measure.

The property
$\overline {q^{-1}(C(p+\la))} \subset q^{-1}(C(p))$
for all $\la \in \cLa_\infty$ is obtained with
a progressive packing strategy. We first construct
$q^{-1}(C(p+\la))$ for all $\la = 2^j r$ with $j \geq 0$, by defining
a partition of a closed subset of 
$q^{-1}(C(p))$ of measure 
$\sum_{\la \in \La_\infty,|\la| \geq 1} \|U[p+\la] \delta \|^2$.
The remaining $q^{-1} (C(p+\la))$ 
are then progressively constructed
for $\la = 2^j r$ and $j$ going from $-1$ to $-\infty$,
within the remaining closed subset of $q^{-1}(C(p))$ not already
allocated. This is possible since we guarantee that the frontier
of each $q^{-1}(C(p))$ has a zero measure. 
$\Box$

The function $q(\omega)$ maps the Lebesgue measure into the Dirac 
scattering measure, but it is discontinuous at all 
$\omega \in \R^{d+}$ such that $q(\omega) \in \cP_\infty$.
Indeed these $\omega$ are then at a boundary of the subdivision procedure
used to construct $q(\omega)$. As a result,
if $\omega$ and $\omega'$ are
on opposite sides of a subdivision boundary then they
are mapped to paths $q(\omega)$ and $q(\omega')$ whose distance
$\bar d(q(\om),q(\om'))$ does not converge to $0$ as $|\omega - \omega'|$
goes to $0$.

Measure preservation (\ref{additnasfd}) implies that $q(\omega)$
defines a
scattering function $\oS f(q(\omega)) \in \Ld(\R^{d+})$ with
\[
\|\oS f(q(\omega))\|_{\R^{d+}}^2 = 
\int_{\R^{d+}} |\oS f(q(\omega))|^2 \, d\om = 
\int_{\cPinf} |\oS f(q)|^2 \, d\mu(q) = 
\|\oS f\|_\cPinf^2~.
\]
If $f$ is a complex-valued function,
then $\cPinf$ is a union of positive paths 
$q= (\la_1,\la_2,\la_3...)$ and negative paths $-q = (-\la_1,\la_2,\la_3...)$. 
Setting $q(-\om) = -q(\om)$  
defines a surjective function from $\R^d$ to $\cPinf$ which
satisfies (\ref{additnasfd}). It results that
$ \oS f (q(-\omega)) = \oS f (-q(\omega))$ for all $\om \in \R^d$,
and $\oS f(q(\om)) \in \LD$ with
$\|\oS f(q(\om))\| = \|\oS f\|_\cPinf$.

If $f$ satisfies (\ref{convadfndsoihfdsf}) then
$ \oS f (q(\omega))$ and $|\hat f (\om)|$
have an equivalent decay over dyadic frequency bands, because their
norm is equal over these frequency bands.
Indeed, for a frequency band $\lambda = 2^j r$ of 
radius proportional to $|\lambda| = 2^j$,
the measure preservation 
(\ref{additnasfd}) together with (\ref{sdfn089s}) prove that
$\|U[\la] f\| = \|f \star \psi_{\la}\|$ satisfies
\begin{equation}
\label{dfosinsdfs}
\int_{q^{-1}(C(\la))}
| \oS f (q(\omega))|^2 \,d\omega = \|U[\la] f\|^2 = 
\frac 1 {2 \pi} \int |\hat f(\omega)|^2\, |\hat \psi (\la^{-1} \omega)|^2\,d \omega~.
\end{equation}
If Conjecture \ref{conjectsdf} is valid then this is true for all $f \in \LD$.
In dimension $d = 1$,
$q^{-1}(C(2^j)) = [\|\psi\|^2 2^{j} , \|\psi\|^2 2^{j+1})$ and
$|\hat \psi(2^j \om)|$ is non-negligible on a similar dyadic frequency
interval.  Hence
$\oS f (q(\omega))$ and $|\hat f(\om)|$
have equivalent energy over dyadic frequency intervals.

\begin{figure}[htb]
\includegraphics[width=\columnwidth]{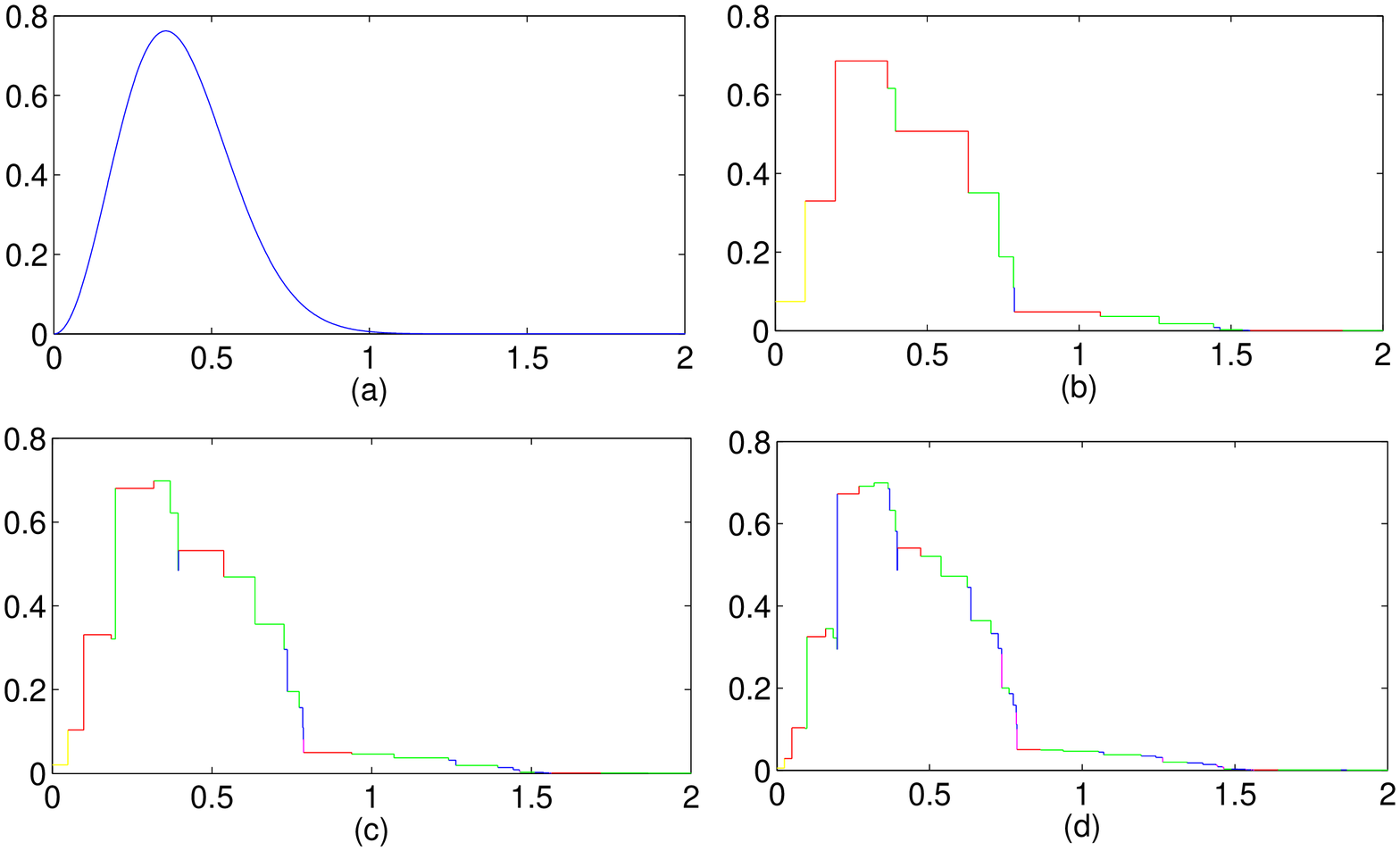} 
	\caption{(a): Fourier modulus $|\hat f(\om)|$
of a Gaussian second derivative, as a function of $\om \in [0,2]$. 
(b,c,d): Piecewise constant
graphs of $\oS_J f(q(\omega))$, as a function of $\omega \in [0,2]$. 
The color specifies the length of each path
$q(\omega)$: $0$ is yellow, $1$ red, $2$ green, $3$ blue, $4$ violet.
The frequency resolution $2^{J}$ increases from (b) to (c) to (d), 
and $\oS_J f(q(\omega))$ converges to a limit function $\oS f(q(\om))$.
}
	\label{examplesconv}
\end{figure}

Figure \ref{examplesconv}(c,d,e) illustrates the
convergence of the windowed scattering transform
$\oS_J f (q(\omega)) $ when $J$ increases, for
a Gaussian second derivative $f$. 
$\oS_J f (q(\omega))$ is constant if
$q(\omega) = p$ is constant and hence if
$\omega \in q^{-1}(C_J(p))$. The
frequency interval $q^{-1}(C_J(p))$ has
a width $\mu(C_J(p)) = \|S_J \delta [p] \|^2$,
which goes to zero as $J$ goes to $\infty$ as shown by 
(\ref{sizeneighboh}). When $J$ increases, each 
$q^{-1}(C_J(p))$ is subdivided into smaller intervals $q^{-1} (C_{J+1}(p'))$ 
corresponding to paths $p$ which are prolongations of $p$.
For each $\omega$, the graph color specifies the length of the
path $p = q(\omega)$.
At low frequencies, $q(\omega) = \emptyset$ is shown as
a yellow interval.
Paths $q(\omega)$ of length $1$ to $4$ are respectively coded in
red, green, blue and violet. 

In these numerical examples, the total energy of
$\oS_J f(q(\om))$ on frequency-decreasing paths $q(\omega)$ is about 
$10^{5}$ times larger than the energy of scattering coefficients on all
other paths. We thus only compute $\oS_J f(q(\om))$ for
frequency-decreasing paths, with an $O(N \log N)$ filter
bank algorithm described in \cite{MallatEUSIPCO}. It is implemented with
the complex cubic spline Battle-Lemari\'e wavelet $\psi$.
As expected from (\ref{dfosinsdfs}), 
$\oS_J f (q(\omega)) f$ has an amplitude and a frequency
localization which is similar to the Fourier modulus $|\hat f(\omega)|$ shown
in Figure \ref{examplesconv}(a).
The discontinuities of 
$ \oS (q(\omega)) f$ along $\omega$ are produced by the discontinuities of
the mapping $q(\omega)$, as opposed to discontinuities of 
$ \oS (q) f$ relatively to the scattering metric in $\cPinf$.

\begin{figure}[htb]
\begin{tabular}{c}
	\includegraphics[width=\linewidth]{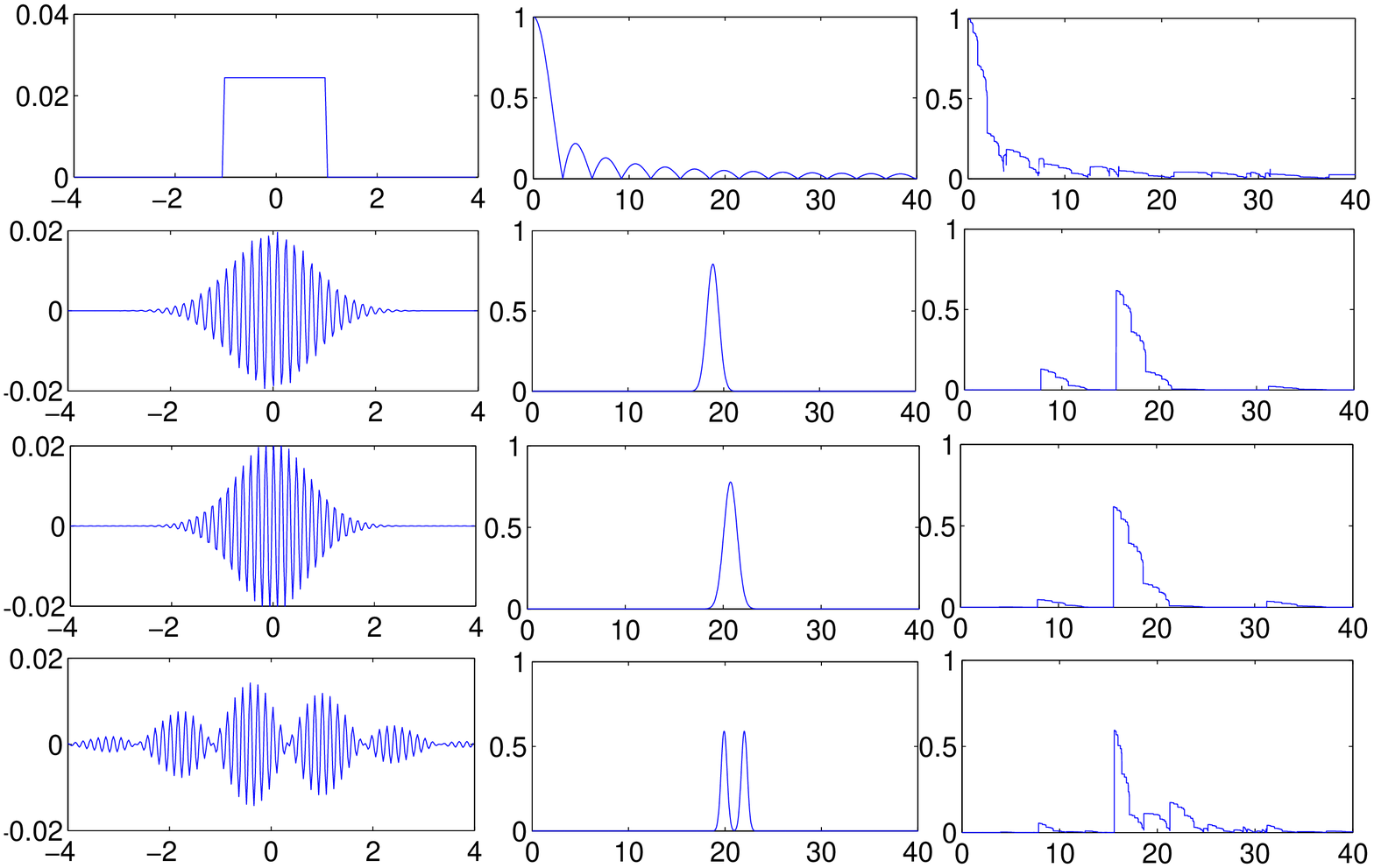}\\
{
(a)\hspace{3.5cm}
(b)\hspace{3.5cm}
(c)}
\end{tabular}
	\caption{(a): Each row $1 \leq i \leq 4$
gives an example of function $f_i (x)$. 
(b): Graphs of the Fourier modulus 
$|\hat f_i (\omega)|$, as a function of $\omega$.
(c): Graphs of the scattering $\oS f_i (q(\omega))$, as a function of $\omega$.}
	\label{examples}
\end{figure}

Figure \ref{examples} compares $ \oS (q(\omega)) f_i$
and $|\hat f_i(\omega)|$ for four functions $f_i$ with $1 \leq i \leq 4$. 
For $f_1 = 1_{[0,1]}$,
the first row of Figure \ref{examples} shows that
$|\hat f_1 (\omega)| = O((1+|\omega|)^{-1})$ has
the same decay in $\omega$ as $ \oS f_1 (q(\omega)) $.
The second row corresponds to 
a Gabor function
$f_2 (x) = e^{i \xi x}\, e^{-x^2/2}$ and the third row shows a small
scaling $f_3 (x) = f_2 ((1-s) x)$ with $s = -0.1$.
The support
of $\hat f_3(\omega) = (1-s)^{-1}\,\hat f_2 ((1-s)^{-1} \omega)$ 
is shifted towards
higher frequencies relatively to the support of $\hat f_2$. A numerical
computation gives
$\||\hat f_2 | - |\hat f_3 |\| = C\, |s|\, \|f_2 \|$
with $C = 13.5$. 
As shown by (\ref{instablsndf}), the constant
$C$ grows proportionally to
the center frequency $\xi$ of $\hat f_2$. It illustrates the
instability of the Fourier modulus under the action of diffeomorphisms.
On the contrary,
the scattering distance remains stable. We numerically obtain
$\| \oS f_2  -  \oS f_3 \| = C\,|s|\, \|f_2 \|$, 
with $C = 1.5$, and this constant does not grow with $\xi$.
It illustrates the Lipschitz continuity of a scattering 
relatively to deformations. 
In the fourth row, $f_4$ is a sum of two high-frequency
Gabor functions, and $|\hat f_4 (\om)|$ includes two narrow
peaks localized within the support of $\hat f_3$. 
The wavelet transform has a bad frequency localization at
such high frequencies, and 
can not discriminate the two frequency peaks
of $\hat f_4$ from $\hat f_3$. However, these two frequency
peaks create low frequency interferences, which appear in the
graph of $f_4$, and which are captured by second order scattering coefficients.
As a result, $\oS f_4$ is very different from $\oS f_3$, which illustrates
the high frequency resolution of a scattering transform obtained 
through interferences.

\section{Scattering Stationary Processes}
\label{Stationary}

A scattering defines a representation of 
stationary processes in $\ld(\cP_\infty)$,  
having different properties than a Fourier power spectrum.
The Fourier power spectrum depends only on second-order moments.
A scattering transform incorporates higher-order moments
that can discriminate processes having same second-order moments. Section \ref{consdfessec} shows that it is 
Lipschitz continuous to random deformations, up to a log term.

\subsection{Expected Scattering}
\label{SecDefoProc}

The properties of a scattering transform 
in $\LD$ are extended to 
stationary processes $X(x)$ with finite second-order moments.
The r\^ole of $\LD$ norm on functions is replaced by the 
mean square norm $E(|X(x)|^2)^{1/2}$ on stationary stochastic processes, 
which does not depend upon $x$ and is thus denoted
$E(|X|^2)^{1/2}$.
Convolutions as well as a modulus preserve stationarity.
If $X(x)$ is stationary, it results that
$U[p] X(x)$ is also stationary and its expected value thus does not
depend upon $x$.

\begin{definition}
The expected scattering transform of a sationary process $X$ is
defined for all $p = (\la_1,...,\la_m) \in \cP_\infty$ by
\[
\overline S X(p) = E (U[p] X)= E(|~|X \star \psi_{\la_1}| \star ... | \star \psi_{\la_m}|).
\]
\end{definition}

This definition replaces the normalized integral of the 
scattering transform (\ref{firsndfinsdfsdfon}) 
by an expected value.
The expected scattering distance between
two stationary processes $X$ and $Y$ is 
\[
\|\oS X  -\oS  Y \|^2 = \sum_{p \in \cP_\infty} |\oS X(p) - \oS Y(p)|^2~.
\]

Scattering coefficients depend upon
normalized high order moments of $X$. This is shown by
decomposing
\[
|U[p]X(x)|^2 = E (|U[p] X |^2)\, (1 + \epsilon(x))~.
\]
A first-order approximation assumes that
$|\epsilon| \ll 1$. Since $\int \psi_{\la} (x)\,dx = 0$,
and $U[p]X(x) = \sqrt{|U[p]X(x)|^2}$,
computing $U[p+\la] X  = |U [p] X \star \psi_{\la}|$ 
with $\sqrt{1+\epsilon} \approx 1+\epsilon/2$ gives
\begin{equation}
\label{resd0f9unsf}
U[p+\la] X  \approx \frac{ ||U [p] X|^2 \star \psi_{\la}|} {2\,E (|U[p] X|^2)^{1/2}}~.
\end{equation}
Iterating on (\ref{resd0f9unsf}) 
proves that  $\oS X(p) = E (U [p] X)$ for $p = (\la_1,...,\la_m)$
depends on normalized moments of $X$
of order $2^{m}$, successively filtered by the wavelets
$\psi_{\la_k}$ for $1 \leq k \leq m$. 

The expected scattering transform is estimated by computing 
a windowed scattering transform of a realization $X(x)$:
\[
S_J[\cP_J] X = \{S_J[p] X \}_{p \in \cP_J}~~\mbox{with}~~S_J [p] X = U[p] X \star \phi_{2^J}.
\]
Since $\int \phi_{2^J} (x)\,dx= 1$, it results that
$E (S_J[p] X)  = E (U[p] X) = \oS X(p)$. So
$S_J[\cP_J] X$ is an unbiased estimator of $\{\oS X(p)\}_{p \in \cP_J}$.

The autocovariance of a real stationary process $X$ is denoted
\[
RX(\tau) = E \Bigl((X(x) - E (X))\,(X(x-\tau) - E (X))\Bigr)~.
\]
Its Fourier transform $\widehat RX (\omega)$ is the power spectrum of $X$.
The mean-square norm of $S_J[\cP_J] X = \{S_J [p] X \}_{p \in \cP_J}$ is written
\[
E (\|{S_J[\cP_J]\, X} \|^2)
= \sum_{p \in \cP_J}  
E ( |{S_J[p] X}|^2).
\]
The following proposition proves that
$S_J[\cP_J]X$ and 
$\oS X$ are nonexpansive and 
that  $\oS X \in \ld(\cP_\infty)$.
The wavelet $\psi$ is assumed to satisfy
the Littlewood-Paley condition (\ref{consenas}).

\begin{proposition}
\label{confppfoidnsf}
If $X$ and $Y$ are finite second-order stationary processes then
\begin{equation}
\label{relastn1050}
E (\|S_J[\cP_J] X - S_J[\cP_J] Y\|^2)
\leq E (|X - Y|^2)~,
\end{equation}
\begin{equation}
\label{relastn105}
\| \oS X - \oS Y\|^2
\leq E (|X - Y|^2)~
\end{equation}
and
\begin{equation}
\label{relastn1059}
\|\oS X \|^2 \leq E (|X|^2)~.
\end{equation}
\end{proposition}

{\it Proof:} We first show that the wavelet transform 
$W_J X = \{A_J X , (W[\la] X)_{\la \in \cLa_J}\}$ 
is unitary over stationary processes.
Let us denote
\[
E (\|{W_J} X\|^2) = 
E (| {A_J} X|^2) + \sum_{\la \in \cLa_J} E (| W[\la] X|^2)~.
\]
Both $A_J X = X \star \phi_{2^J}$ and
$W[\la] X = X \star \psi_\la$ are stationary.
Since $\int \phi_{2^J}(x)\,dx=1$ and $\int \psi_\la(x)\,dx=0$
it results that $E (A_J X) = E (X)$ and $E (W[\la] X) = 0$. 
Since the power spectrum of $A_J X$ and $W[\la] X$ is
respectively $\widehat RX (\om)\, |\hat \phi(2^J \om)|^2$ and
$\widehat RX (\om)\, |\hat \psi_\la( \om)|^2$, we get
\[
E (|A_J X  |^2) = \int \widehat RX (\om)\, |\hat \phi(2^J \om)|^2 
\,d\om + E (X)^2~
\]
and
\[
E (|W[\la] X |^2) = \int \widehat RX (\om)\, |\hat \psi_\la( \om)|^2\,d \om.
\]
Since $E (|X|^2) = \int \widehat RX (\om)\,d\om + E (X)^2$,
the same proof as in Proposition \ref{Littlewood} shows that
the wavelet condition
(\ref{consenas}) implies that $E (\| W_J X\|^2)  = E (|X |^2)$.

The propagator 
${U_J} X = \{ A_J X , (|W[\la] X| )_{\la \in \cLa_J} \}$ satisfies
\[
E (\|{U_J} X - {U_J} Y\|^2) 
\leq E (\|{W_J} X - W_J Y\|^2) = E (|X - Y|^2)
\]
and is thus nonexpansive on stationary processes. 
We verify as in (\ref{propasndfs}) that
\[
\oUJ \,U[\cLa_J^m] X  = 
\{S_J[\cLa_J^{m}] X\,,\, U[\cLa_J^{m+1}] X \}~.
\]
Since $\cP_J = \cup_{m=0}^{+\infty} \cLa_J^m$, one can compute
$S_J[\cP_J] X$ by iteratively applying the nonexpansive
operator $U_J$. The nonexpansive property (\ref{relastn1050}) 
is derived from the fact that $U_J$ is nonexpansive,
as in Proposition \ref{unidfnsdofnw}.

Let us prove (\ref{relastn105}). Since $\oS X(p) = E (S_J [p] X)$ and
$\oS Y(p) = E (S_J [p] Y)$ 
\[
\sum_{p \in \cP_J} | \oS X(p) - \oS Y(p)|^2 \leq E (\|S_J[\cP_J] X - S_J[\cP_J] Y\|^2)
\leq E (|X - Y|^2)~.
\]
Letting $J$ go to $\infty$ proves (\ref{relastn105}).
The last inequality (\ref{relastn1059}) is obtained by setting $Y = 0$. 
$\Box$

Paralleling the scattering norm preservation
in $\LD$, the following theorem proves
that $S_J[\cP_J]$ preserves the mean-square norm of stationary processes.

\begin{theorem}
\label{theenconfstat}
If the wavelet satisfies the admissibility
condition (\ref{conditionprogre})
and if $X$ is stationary with $E(|X|^2) < \infty$ then
\begin{equation}
\label{relastn1023}
E (\|S_J[\cP_J]X \|^2) = E (|X|^2)~.
\end{equation}
\end{theorem}

{\it Proof:}
The proof of (\ref{relastn1023}) is almost identical to the proof of 
(\ref{dnslnfs8302}) in Theorem \ref{energydecth}, if we replace 
$f$ by $X$, $|\hat f(\omega)|^2$ by the power
spectrum $\widehat RX (\omega)$ and $\|f\|^2$ by $E (|X|^2)$. 
We proved that $E (\|W_J X\|^2)  = E (|X |^2)$ so we also have
$E (\|U_J X\|^2)  = E (|X |^2)$.
In the derivations of Lemma \ref{lemma8}, 
replacing $f_p= U[p] f$ by $X_p= U[p] X$,
and $|\hat f_p (\om)|^2$ by
$\widehat R{X_p} (\om)$, proves that
\[
\frac {\alpha} 2 \, E (\|U[\cP_J]X \|^2) \leq \max(J+1,1) \, E(|X|^2) + 
\sum_{j > 0}\sum_{r \in G^+} j\, E (|X \star \psi_{2^j\ga} |^2)~.
\]
Since $\cP_J = \cup_{m \in \N} \cLa_J^m$, if the right hand-side term is finite then 
\begin{equation}
\label{relastn102}
\lim_{m \rightarrow \infty} E (\|U[\cLa_J^m]X \|^2)  = 0~.
\end{equation}
The same density argument as in the proof of 
Theorem \ref{energydecth} proves that (\ref{relastn102}) also holds 
if $E (|X|^2)< \infty$ because $\widehat R{X} (\om)$ is integrable.

Since $E (\| U_J X\|^2)  =  E (|X |^2)$
and $\oUJ \,U[\cLa_J^m] X  = \{ S_J[\cLa_J^{m}] X \,,\, U[\cLa_J^{m+1}] X\}$, 
iterating $m$ times on $U_J$ proves as in (\ref{Unfidnfsdf7}) that
\[
E (|X |^2) = \sum_{n=0}^{m-1} E (\| S_J[\cLa^n_J] X\|^2) + 
E(\| U[\cLa^m_J] X \|^2) ~.
\]
When $m$ goes to $\infty$, (\ref{relastn102}) implies
(\ref{relastn1023}).
$\Box$

A windowed scattering $S_J [p] = U[p] X \star \phi_J$
averages $U[p] X$ over a domain whose size is proportional to $2^J$.
If $U[p] X$ is ergodic, 
it thus converges to $\overline SX(p) = E (U[p] X)$ when $J$ goes to $\infty$.
The windowed transformed scattering $S_J[\cP_J] X$ is said to be
a {\it mean-square consistent} estimator of $\oS X$ if its
total variance over all paths converges to zero:
\[
\lim_{J \rightarrow \infty}  E ( \|S_J[\cP_J] X - \oS_J X \|^2) = 
\lim_{J \rightarrow \infty}  \sum_{p \in \cP_J} E ( |S_J[p] X - \oS X(p) |^2) = 
0~.
\]
Mean-square convergence implies convergence in probability and hence
that $S_J[\cP_J] X$ converges to $\oS X$ with
probability $1$. 

\begin{figure}[htb]
	\includegraphics[width=\columnwidth]{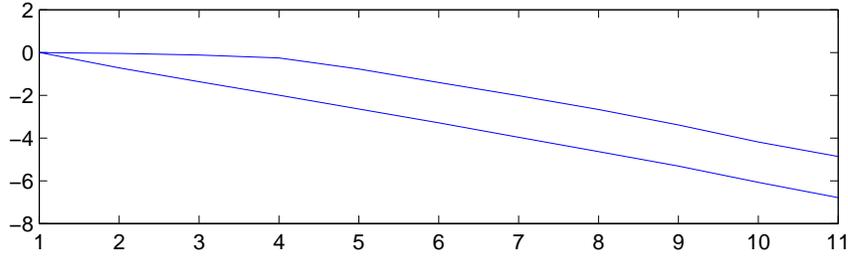}

	\caption{Decay of $\log_2 E ( \|S_J[\cP_J] X - \oS X\|^2)$ 
as a function of $J$ for a Gaussian white noise X (bottom line) and
a moving average Gaussian process (top line), along frequency-decreasing
paths.}
	\label{decayfigu}
\end{figure}

For a large class of ergodic processes $X$, including Gaussian processes,
mean-square convergence is observed numerically,
with  $E ( |S_J[\cP_J] X - \oS_J X) |^2) \leq C \, 2^{-\alpha J}$
for $C > 0$ and $\alpha > 0$. When $J$ increases,
the global variance of $S_J[\cP_J] X$ decreases
despite the path subdivision into new paths 
because each modulus reduces the variance by removing random phase variations.
The variance of $S_J[p]X$ thus decreases when the path length increases, 
and it is concentrated over a small number of frequency-decreasing paths.
For a Gaussian white noise and a moving average Gaussian process of unit
variance, Figure \ref{decayfigu} shows that
$\log E ( \|S_J[\cP_J] X - \oS X \|^2)$, computed over all frequency-decreasing paths, decays linearly as a function $J$.
For the correlated Gaussian process, the decay begins for $2^J \geq 2^{4}$,
which is the correlation length of this process. Indeed,
the averaging by $\phi_{2^J}$ effectively reduces the estimator variance
when $2^{J}$ is bigger than the correlation length.

\begin{conjecture}
If $X$ is a Gaussian stationary process
with $\|RX \|_1 < \infty$ then $S_J[\cP_J X]$ 
is a mean-square consistent estimator of  $\oS X$.
\end{conjecture}

The following corollary of
Theorem \ref{theenconfstat} proves that 
mean-square consistency implies an expected 
scattering energy conservation.

\begin{corollary}
\label{theenconfstatcor}
For an admissible scattering wavelet which satisfies condition
(\ref{conditionprogre}),
$S_J[\cP_J]X$ is mean-square consistent if and only if
\begin{equation}
\label{consisnten20}
\|\oS X \|^2 = E ( |X|^2)~,
\end{equation}
and mean-square consistency implies that for all $\lambda \in \Lambda_\infty$
\begin{equation}
\label{consisnten2}
\sum_{p \in \cP_\infty} |\oS X (\lambda + p )|^2 = 
E ( |X \star \psi_{\lambda}|^2)~.
\end{equation}
\end{corollary}

{\it Proof:}
It results from Theorem \ref{theenconfstat} that
$E (\|S_J[\cP_J ]X \|^2) = E (|X|^2)$. Since
\[
E (\|S_J[\cP_J] X \|^2)  = \sum_{p \in \cP_J} E (S_J[p] X)^2 + 
E (| S_J[\cP_J] X - E (S_J[\cP_J]X) |^2),
\]
and $E (S_J[p] X) = \oS X(p)$, 
we derive that
$\lim_{J \rightarrow \infty} E (\| S_J[\cP_J] X - E(S_J[\cP_J] X) \|^2) = 0$
if and only if $\|\oS X \|^2 = E ( |X|^2)$.
Moreover, for all $\la \in \Lambda_\infty$,
since $U[p] U[\la] X = U[\la+p] X$,
applying (\ref{consisnten20}) to $U[\la] X$ instead of 
$X$ proves (\ref{consisnten2}).
$\Box$

The expected scattering can 
be represented by a singular scattering spectrum in $\cPinf$.
Similarly to
Section \ref{Renonscat}, we associate to 
$\oS X(p) = E (U[p] X)$ 
a function that is piecewise constant in $\cPinf$:
\begin{equation}
\label{convsdfn8usdf0}
\forall q \in \cPinf~,~
P_J X(q) 
= \sum_{p \in \cP_J} {\oS X(p)^2} \,
\frac {1_{C_J (p)} (q)} {\|S_J[p] \delta \|^2}\, .
\end{equation}
The following proposition proves that $P_J$ converges to a singular measure,
called a {\it scattering power spectrum}.

\begin{proposition}
$P_J X(q)$ converges 
in the sense of distributions to a 
Radon measure in $\cPinf$, supported in $\cP_\infty$:
\begin{equation}
\label{convsdfn8usdf}
PX(q) = 
\lim_{J \rightarrow \infty} P_J X(q) = 
\sum_{p \in \cP_\infty} \oS X(p)^2 \, \delta(q-p)~.
\end{equation}
\end{proposition}

{\it Proof:} For any $p \in \cP_\infty$, the Dirac
$\delta(p-q)$ is defined as a linear form satisfying
$\int_{\cPinf} f(q) \, \delta(p-q) d\mu(q) = f(p)$ for 
all continuous functions $f(q)$ of $\cPinf$ relatively to the
scattering metric.
For all $J \in \Z$,
$\mu(C_J(p)) = {\|S_J[p] \delta \|^2}$, $p \in C_J(p)$.
and $\lim_{J \rightarrow \infty} \mu(C_J(p)) = 0$.
We thus obtain the following convergence 
in the sense of distributions:
\[
\lim_{J \rightarrow \infty } \frac{1_{C_J(p)} (q)} 
{\|S_J[p] \delta \|^2} = \delta (q-p) . 
\]
Letting $J$ go to $\infty$ in (\ref{convsdfn8usdf0}) 
proves (\ref{convsdfn8usdf}). $\Box$

If $S_J[\cP_J] X$ is mean-square consistent then 
(\ref{consisnten2}) implies that the scattering spectrum $PX(q)$ is
related to the Fourier power spectrum $\widehat RX (\omega)$ by
\begin{equation}
\label{powerspcvndsf}
\int_{C(\la)}
PX (q) \,d \mu(q) = 
E (|X \star \psi_{\la}|^2) = 
\frac 1 {2 \pi} \int \widehat RX (\omega)\, |\hat \psi (\la^{-1} \omega)|^2\,d \omega~.
\end{equation}
Let $q(\om)$ be the function of
Proposition \ref{isodfmpsdfon},
which maps the Lebesgue
measure of $\R^{d+}$ into the Dirac scattering measure of $\cPinf$.
It defines a scattering power spectrum 
$PX(q(\omega))$ over the half frequency space $\omega \in \R^{d+}$.
In dimension $d=1$, $q^{-1}(C(2^j)) = [\|\psi\|^2 2^{j}\,,\,\|\psi\|^2 2^{j+1})$,
so (\ref{powerspcvndsf}) implies
\[
\int_{\|\psi\|^2 2^{j}}^{\|\psi\|^2 2^{j+1}}
PX (q(\omega)) \,d \omega = 
\frac 1 {2 \pi} \int \widehat RX (\omega)\, |\hat \psi (2^j \omega)|^2\,d \omega~.
\]
Although $PX(q(\omega))$ and  $\widehat RX (\omega)$ have the same integral over
dyadic frequency intervals, they have
very different distributions within each of these intervals.
Indeed, (\ref{resd0f9unsf}) shows that if $p$ is of length $m$ then
$E (U[p] X)$ depends upon 
normalized moments of $X$ of order $2^{m}$. It results that
$PX(q(\omega))$ depends upon arbitrarily high order moments of $X$
where as $\widehat RX(\omega)$ only depends upon moments of order $2$. 
Hence, $PX (q)$ can discriminate different
stationary processes having same Fourier power spectrum
and thus same second-order moments.

\begin{figure}[htb]
\begin{tabular}{c}
	\includegraphics[width=\columnwidth]{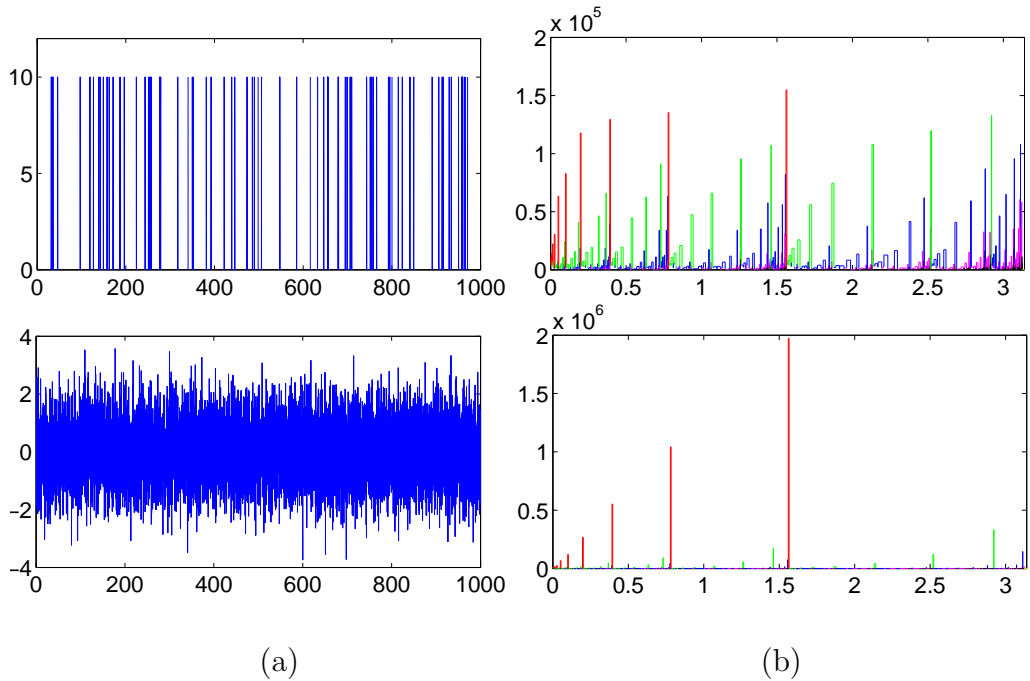}\\
{
(a)\hspace{6cm}
(b)}
\end{tabular}
	\caption{(a): Realization of a 
Bernoulli process $X_1(x)$ at the top and a 
Gaussian white noise $X_2(x)$ at the bottom, both having a unit variance.
(b): Scattering power spectrum 
$PX_i(q(\omega))$ of each process,
as a function of $\omega \in [0,\pi]$. 
The values of $PX_i(q(\omega))$ are displayed respectively
in red, green, blue and violet, for paths $q(\omega)$ of
length $1$, $2$, $3$ and $4$.}
	\label{randexamples}
\end{figure}

Figure \ref{randexamples} gives the scattering power
spectrum of a Gaussian
white noise $X_2$ and of a Bernoulli process $X_1$ in dimension $d = 1$, 
estimated from a realization sampled over $N = 10^4$ integer points. 
Both processes have a constant Fourier power spectrum $\widehat RX_i (\om) = 1$
but very different scattering spectrum.
Their scattering spectrum
$P X_i (q(\omega))$ is estimated by $P_J X_i (q(\omega))$ in
(\ref{convsdfn8usdf0}) at the maximum scale $2^J = N$. 
It is a sum of spikes in Figure \ref{randexamples}(b), which  
converges to a Radon measure supported in $\cP_\infty$ when increasing
$2^J = N$.
A Gaussian white noise $X_2$ has a scattering spectrum mostly concentrated
on paths $q(\omega) = (2^j)$ of length $1$. 
These scattering coefficients appear as 
large amplitude red spikes at dyadic positions, in the bottom graph of
Figure \ref{randexamples}(b). Their amplitude is proportional
to $\oS X_2 (2^j)^2 \sim 2^{j}$. 
Other spikes in green, correspond to paths $q(\omega) = (\la_1,\la_2)$
of length two. They have a much smaller amplitude.
Scattering coefficients for paths of length $3$ and $4$,
in blue and violet, are so small that they are not visible. 
The top of Figure \ref{randexamples}(b) shows 
the scattering spectrum $P_J X_1 (q(\omega))$ 
of a Bernouilli process $X_1$. It has
a maximum amplitude for paths $q(\omega)$ of length $1$ (in red), 
but longer paths
shown in green, blue and violet also produce large scattering
coefficients, as opposed to a Gaussian white noise scattering. 
Scattering coefficients for paths $p$ of length
$m$ depend upon the moments of $X$ up to the order $2^m$. 
For $m > 1$, large scattering coefficients
indicate a strongly non-Gaussian behavior of high order moments.

\subsection{Random Deformations}
\label{consdfessec} 

We now show that the scattering transform is nearly Lipschitz
continuous to the action of random deformations.
If $\tau$ is a random process with 
$\|\nabla \tau \|_\infty = |\nabla \tau (x)| < 1$ then $x - \tau (x)$ is a random diffeomorphism.
If $X(x)$ and $\tau (x)$ are 
independent stationary processes then the action of this random diffeomorphism
on $X(x)$ defines a randomly deformed process
$L_\tau X(x) = X(x-\tau(x))$ which remains stationary.

The following theorem adapts the result of Theorem \ref{DefomsProps}
by proving that the scattering distance produced by a random
deformation is dominated by a first-order term proportional to 
$E (\|\nabla \tau\|^2_\infty )$. Let us denote
\[
E ( \|U[\cP_J] X\|_{1}) = 
\sum_{m=0}^{+\infty} \left(\sum_{p \in \cLa_J^m} E (|U[p] X|^2) \right)^{1/2}~
\]
where $\cLa_J^m$ is the set of paths $p = (\la_k)_{k \leq m}$ 
of length $m$ with $|\la_k| < 2^{J}$.

\begin{theorem}
\label{Theoelas6} 
There exists $C$ such that for all independent stationary processes
$\tau$ and $X$ 
satisfying $\|\nabla \tau\|_\infty \leq 1/2$  with probability $1$,
if $E (\|U[\cP_J] X \|_1) < \infty$ then
\begin{equation}
\label{err109213}
E (\|S_J[\cP_J] L_\tau X - S_J[\cP_J]X \|^2)
\leq  C\,E (\|U[\cP_J]X\|_1)^2\,K (\tau)
\end{equation}
with
\begin{equation}
\label{err109215}
K(\tau) = 
 E\left(\Bigl(2^{-J} \|\tau \|_\infty + \|\nabla \tau \|_\infty 
(\log \frac{\|\Delta \tau\|_\infty}{\|\nabla \tau \|_\infty }
\vee 1)
+ \|\rH \tau \|_\infty \Bigr)^2\right) . 
\end{equation}
Over the subset $\cP_{J,m}$ of path in $\cP_J$ of length strictly smaller than $m$
\begin{equation}
\label{err109218}
E (\|S_J[\cP_{J,m}] L_\tau X - S_J[\cP_{J,m}]X \|^2)
\leq  C\,m\, E (|X|^2)\,K (\tau) ~.
\end{equation}
\end{theorem}

{\it Proof:} Similarly to the proof of
Theorem \ref{DefomsProps}, we decompose
\begin{eqnarray*}
E (\|S_J[\cP_J] L_\tau X - S_J[\cP_J]X \|^2) &\leq& 2\,
E(\|L_\tau S_J[\cP_J] X - S_J[\cP_J]X \|^2) \\
& & + 2 \,E(\|[S_J[\cP_J], L_\tau] X\|^2).
\end{eqnarray*}
Appendix \ref{theornadofns} proves that
$E(\|[S_J[\cP_J]\,,\, L_\tau] X \|^2) \leq E(\|U[\cP_J] X \|_1 )^2\,B(\tau)$
with
\begin{equation}
\label{commnaf87}
B(\tau) = 
C^2\,E \left( \Bigl(\|\nabla \tau \|_\infty (\log \frac{\|\Delta \tau\|_\infty}{\|\nabla \tau \|_\infty }
\vee 1)+ \|\rH \tau \|_\infty \Bigr)^2\right) ,
\end{equation}
and since
\begin{equation}
\label{commnaf87cnfs}
E(\|L_\tau S_J[\cP_J] X - S_J[\cP_J]X \|^2) \leq C^2\,E(\|U[\cP_J] X \|_1 )^2\,
E( 2^{-J} \|\tau\|_\infty^2)~,
\end{equation}
we get (\ref{err109213}). The commutator
$[S_J[\cP_J],L_\tau]$ and $L_\tau S_J[\cP_J]  - S_J[\cP_J]$
are random operators since $\tau$ is a random process.
The key argument of the proof
is provided by the following lemma which relates the 
expected $\LD$ sup norm of a random operator
to its norm on stationary processes.
This lemma is proved in Appendix \ref{prooflemmarandom}. 

\begin{lemma}
\label{lemmarandom}
Let $K_\tau$ be an integral operator with a kernel $k_\tau (x,u)$ 
which depends upon a random process $\tau$. 
If the following two conditions are satisfied
\[
E\Bigl(k_\tau (x,u)\,k^*_\tau (x,u') \Bigr) = \bar k_\tau (x-u,x-u')~~\mbox{and}~~
\iint |\bar k_\tau (v,v')| \,|v-v'|\,dv\,dv' < \infty~,
\]
then for any stationary process $Y$ independent of $\tau$, 
$E(|K_\tau Y (x)|^2)$ does not depend upon $x$ and  
\begin{equation}
\label{indqualityopdfns}
E(|K_\tau Y |^2) \leq E(\|K_\tau \|^2)\, E(|Y|^2)~,
\end{equation}
where $\|K_\tau\|$ is the operator norm in $\LD$ for each realization of $\tau$. 
\end{lemma}

This result remains valid when
replacing $S_J[\cP_J]$ by $S_J[\cP_{J,m}]$ and
$U[\cP_J]$ by $U[\cP_{J,m}]$. With the same argument as in the proof
of (\ref{indwansford9}), we verify that
\[
E(\|U[\cP_{J,m}]X \|_1) \leq m\,E(|X |^2)^{1/2}
\]
which proves (\ref{err109218}).
$\Box$

Small stationary deformations of stationary processes result in small
modifications of the scattering distance, which is important to characterize
deformed stationary processes as in image textures \cite{Bruna}. 
The following corollary proves that the expected
scattering transform
is almost Lipschitz continuous in the size of
the stochastic deformation gradient $\nabla \tau$, up to a log term. 

\begin{corollary}
\label{Theoelas6df} 
There exists $C$ such that for all independent stationary processes
$\tau$  and $X$ satisfying
$\|\nabla \tau\|_\infty \leq 1/2$  with probability $1$, if
$E(\|U[\cP_\infty]X\|_1) < \infty$ then
\begin{equation}
\label{err109213df}
\|\oS L_\tau X - \oS X \|^2 
\leq  C\,E(\|U[\cP_\infty]\|_1) ~ E(|X|^2)\,K (\tau)
\end{equation}
with
\begin{equation}
\label{err109215r}
K(\tau) = 
 E\Bigl\{\Bigl(\|\nabla \tau \|_\infty (\log \frac{\|\Delta \tau\|_\infty}{\|\nabla \tau \|_\infty }
\vee 1)
+ \|\rH \tau \|_\infty \Bigr)^2\Bigr\} . 
\end{equation}
\end{corollary}

{\it Proof:} $E (\|S_J[\cP_J] L_\tau X - S_J[\cP_J]X \|^2) \leq
\|E (S_J[\cP_J] L_\tau X) - E(S_J[\cP_J]X) \|^2$, so letting $J$ go
to $\infty$ in (\ref{err109213}) 
proves (\ref{err109213df}). $\Box$

\section{Invariance to Actions of Compact Lie Groups}
\label{MultiGroup}

Invariant scattering are extended to the action of compact Lie Groups $G$.
Section \ref{RotatGroup} builds scattering operators in $\Ld(G)$ which
are invariant to the action of $G$.
Section \ref{Combined} defines a translation and rotation invariant
operators on $\LD$ by combining a scattering operator on $\LD$ and a
scattering operator on $\Ld(SO(d))$.

\subsection{Compact Lie Group Scattering}
\label{RotatGroup} 

Let $\G$ be a compact Lie group and
$\Ld (\G)$ be the space of measurable functions $f(\ga)$ such that
$\|f\|^2 = \int_\G |f(\ga)|^2\,d\ga < \infty$, where
$d \ga$ is the Haar measure of $\G$. 
The left action of $g \in G$ on $f \in \Ld(G)$ is defined by
$L_g f (r) = f(g^{-1} r)$.
This section introduces a scattering transform on $\Ld(G)$, which
is invariant to the action of $G$. It is obtained with a 
scattering propagator which
cascades the modulus of wavelet transforms defined on $\Ld (\G)$.

The construction of Littlewood-Paley decompositions on compact manifolds
and in particular on compact Lie groups was developed
by Stein \cite{Stein1}. Different wavelet constructions have been 
proposed over manifolds \cite{Donoho2}.
Geller and Pesenson \cite{Geller1} have built
unitary wavelet transforms on compact Lie groups, 
which can be viewed as analogs of unitary wavelet transforms on
the circle in $\R^2$. In place of sinusoids, they use the eigenvectors
of the Laplace-Beltrami operator of an invariant metric defined on the group.
Similarly to Meyer wavelets \cite{Meyerbook},
these basis elements
are regrouped into dyadic subbands with appropriate windowing. 
For any $2^L \leq 1$, it defines 
a scaling function $\widetilde \phi_{2^L}(r)$
and a family of wavelets $\{ \widetilde \psi_{2^j}(r) \}_{-L< j< \infty}$
which are in $\Ld(\G)$ \cite{Geller1}. The wavelet coefficients of 
$f \in \Ld(\G)$ are computed with left convolutions on the group $G$
for each $\tilde \la = 2^j$:
\begin{equation}
\label{LIewavel}
\widetilde W[\tilde \la] f (r) = f \star \widetilde \psi_{\tilde \la} (r) = \int_\G 
f(g)\, \widetilde \psi_{\tilde \la} (g^{-1}\,r) \,d g~
\end{equation}
and the scaling function performs an averaging on $G$:
\begin{equation}
\label{LIewave2}
\widetilde A_Lf  (r) = f \star \widetilde \phi_{2^L} (r) = 
\int_\G f(g)\, \widetilde \phi_{2^L} (g^{-1}\,r) \,d g~.
\end{equation}
The resulting wavelet transform of $f \in \Ld(\G)$ is
\[
{\widetilde W_L} f = \{ \widetilde A_L f \,,\,
(\widetilde W[\tilde \la] f)_{\tilde \la \in \tilde \cLa_L}\}~~\mbox{with}~~ 
\widetilde \cLa_L = \{ \tilde \la = 2^j~:~j > -L \}~.
\]
At the maximum scale $2^L = 1$, since $\widetilde \phi_1 (r) = \left(\int_\G ~d g \right)^{-1} = |\G|^{-1}$ is constant, the operator $\widetilde A_0$ performs
an integration on the group:
\begin{equation}
\label{A0fns}
\widetilde A_{0} f (r) = |\G|^{-1}\,\int_\G f(g)\, d g = cst.
\end{equation}
Wavelets are constructed to obtain a unitary operator \cite{Geller1}
\begin{equation}
\label{grounsdf09s}
\| \widetilde W_L f \| = \|f\|~,
\end{equation}
with 
\[
\| \widetilde W_L f \|^2 = \|\widetilde A_L f\|^2 + \sum_{\tilde \la \in 
\tilde \cLa_L}
\| \widetilde W[\tilde \la] f \|^2~.
\]

The Abelian group $\G = SO(2)$ of rotations in $\R^2$ is a simple example
parametrized by an angle in $[0,2 \pi]$. 
The space $\Ld(\G)$ is thus equivalent to $\Ld[0,2\pi]$. 
Wavelets in $\Ld(\G)$ are the well-known periodic 
wavelets in $\Ld[0,2\pi]$ \cite{Meyerbook}. 
They are obtained by periodizing
a scaling function $\phi_{2^L}(x) = 2^{-L} \phi(2^{-L} x)$
and 
wavelets $\psi_{2^j} (x) = 2^{j} \psi(2^{j} x)$ with $(\phi,\psi) \in \Ld(\R)^2$:
\begin{equation}
\label{periwavsdfs}
\widetilde \phi_{2^L} (x)= \sum_{m \in \Z} \phi_{2^L} (x - 2 \pi m)\,\,\,\mbox{and}
\,\,\,
\widetilde \psi_{2^j} (x)= \sum_{m \in \Z} \psi_{2^j} (x - 2 \pi m)~.
\end{equation}
We suppose that $\hat \phi(0) = 1$ and
$\hat \phi ( 2k \pi) = 0$ for $k \in \Z-\{0\}$. 
The Poisson formula implies that
$\widetilde \phi (x) = \sum_{n \in \Z} \phi (x - n) = 2 \pi$.
Convolutions (\ref{LIewavel}) and (\ref{LIewave2}) on the rotation 
group are circular convolutions of periodic functions in
$\Ld[0,2\pi]$. With
the Poisson formula, one can prove that the periodic
wavelet transform $\widetilde W_L$ is unitary 
if and only if $(\phi,\psi) \in \Ld(\R)^2$
satisfy the Littlewood-Paley equalities (\ref{consenas}).

For a general compact Lie group $G$, we define
a wavelet modulus operator by
$\widetilde U[\tilde \la] f = |\widetilde W_L[\tilde \la] f|$, and the 
resulting one-step propagator is
\[
\widetilde U_L f = \{ \widetilde A_L f \,,\,
( \widetilde U[\tilde \la] f| )_{\tilde \la \in \tilde \cLa_L} \}~.
\]
Since $\widetilde W_L$ is unitary, we verify as in (\ref{inequasvdons})
that $\widetilde U_L$ is nonexpansive
and preserves the norm in $\Ld(\G)$. 

A scattering operator on $\Ld(\G)$
applies $\widetilde U_L$ iteratively.
Let $\tP_L$ denote the set of all finite paths
$\tilde p = \{\tilde \la_1 , ... , \tilde \la_{m}\}$ of length $m$,
where $\tilde \la_k = 2^{j_k} \in \widetilde \cLa_L$. 
Following Definition \ref{propdefin},
a scattering propagator on $\Ld(G)$
is a path ordered product of non-commutative wavelet modulus
operators 
\[
\widetilde U[\tilde p] = \widetilde U[{\tilde \la_m}]\,...\,\widetilde U[{\tilde \la_1}]~,
\]
with $\widetilde U[\nill]  = Id$. 

Following Definition \ref{windscatdef},
a windowed scattering is defined by averaging
$\widetilde U[\tilde p]f$
through a group convolution with $\widetilde \phi_{2^L}$ 
\begin{equation}
\label{relastn0898sdfn}
\widetilde S_L[\tilde p ] f(r) = \widetilde A_L \widetilde U [\tilde p] f (r) = \int_\G U[\tilde p] f(g)\, \widetilde \phi_{2^L} (g^{-1}\,r) \,d g~.
\end{equation}
It yields an infinite family of functions
$\widetilde S_L[\tP_L] f = \{\widetilde S_L[\tilde p] f \}_{\tilde p \in \tP_L}$,
whose norm is
\[
\|\widetilde S_L[\tP_L] f \|^2 = 
\sum_{ \tilde p \in \tP_L}
\| \widetilde S_L[ \tilde p] f \|^2
~~\mbox{with}~~
\| \widetilde S_L [\tilde p] f\|^2 = \int_\G
| \widetilde S_L [\tilde p] f(r)|^2\,dr~.
\]

Since $\widetilde U[\tP_L]$ 
is obtained by cascading the nonexpansive operator $\widetilde U_L$,
the same proof as in Proposition \ref{unidfnsdofnw} shows that 
it is nonexpansive:
\[
\forall (f,h) \in  \Ld(G)^2~~,~~
\| \widetilde S_L [\tP_L] f - 
\widetilde S_L [\tP_L] h \| \leq \|f - h \|~.
\]
Since $\widetilde U[\tP_L]$ preserves the norm in $\Ld(G)$, to also
prove as in Theorem \ref{energydecth}
that $\|\widetilde S_L [\tP_L] f \| = \|f \|$, it is necessary
to verify that 
$\lim_{m \rightarrow \infty} \|\widetilde U[\widetilde \Lambda_L^m]\|^2 = 0$.
For the translation group where $G = \R^d$, Theorem \ref{energydecth}
proves this result by imposing a condition on the Fourier transform of
the wavelet. The extension of this result is not straightforward on
$\Ld(G)$ for general compact Lie groups $G$, but it remains valid
for the rotation group $G = SO(2)$ in $\R^2$. Indeed, periodic
wavelets $\tilde \psi_{\tilde \lambda} \in \Ld(SO(2)) = \Ld[0,2 \pi]$
are obtained by periodizing wavelets 
$\psi_{\tilde \lambda} \in \Ld(\R)$ in (\ref{periwavsdfs}), 
which is equivalent to subsample
uniformly their Fourier transform. If $\psi$ satisfies the admissibility
condition of Theorem \ref{energydecth} then by replacing convolutions
with circular convolutions in the proof, we verify that the periodized 
wavelets $\tilde \psi_{\tilde \lambda}$ define a scattering transform 
of $\Ld[0,2\pi]$ which preserves the norm 
$\|\widetilde S_L [\tP_L] f \| = \|f \|$. 

When $2^L = 1$, $\widetilde A_0$ is the integration
operator (\ref{A0fns}) on the group, so $\widetilde S_0 [\tilde p] f(r)$ does
not depend on $r$. Following Definition \ref{definwindoscas}, it defines
a scattering transform which maps any
$f \in \Ld(G)$ into a function of the path variable $\tilde p$:
\begin{equation}
\label{comvinsdf0ds}
\forall \tilde p \in \tP_0~~,~~\widetilde S_0 [\tilde p] f  = |\G|^{-1}\,
\int_\G U[\tilde p] f(g)\,dg~.
\end{equation}
Over a compact Lie group, the scattering transform
$\widetilde S_0[\tP_0] f = \{\widetilde S_0[\tilde p] f \}_{\tilde p \in \tP_0}$
is a discrete sequence in $\ld(\tP_0)$. 
The following proposition proves that it is invariant to the action 
$L_g f(r) = f(g^{-1} r)$ of $g \in G$ on $f \in \Ld(G)$.

\begin{proposition}
\label{groupinfaprop}
For any $f \in \Ld(G)$ and $g \in G$ 
\begin{equation}
\label{dispanfdson0}
\widetilde S_0[\tP_0] \,L_g f =  \widetilde S_0[\tP_0] f\,.
\end{equation}
\end{proposition}
 
{\it Proof:} Since $\widetilde A_0$ and $\widetilde W[\tilde \la] f$ 
are computed with
left convolutions on $\G$, they commute with $L_g$.
It results that $\widetilde U[\tilde \la]$ 
and hence $\widetilde S_0[\tP_L]$ also commutes with $L_g$.
If $\tilde p \in  \tP_0$,
since $\widetilde S_0 [ \tilde p] f(r)$ is constant in $r$,
$\widetilde S_0[\tilde p] \,L_g f = L_g \, \widetilde S_0 [ \tilde p] f
= \widetilde S_0 [\tilde p] f$,
which proves (\ref{dispanfdson0}). $\Box$

As in the translation case, 
the Lipchitz continuity of $\widetilde S_L$ 
to the action of diffeomorphisms 
relies on the Lipschitz continuity of the wavelet transform
$\widetilde W_L$.
The action of a small diffeomorphism on $f \in \Ld(G)$ can be written 
$L_\tau f(r) = f(\tau(r)^{-1} r)$ with $\tau (r) \in G$. The proof of
Theorem \ref{DefomsProps} on Lipschitz continuity 
applies to any compact Lie groups $G$. The main difficulty
is to prove Lemma \ref{thedw0}, which proves the Lipschitz continuity
of $\widetilde W_L$ by computing an upper bound of the
commutator norm $\|[\widetilde W_L , L_\tau ] \|$. 
The proof of this lemma can still be carried by applying Cotlar's lemma,
but integration by parts and the resulting bounds 
require appropriate hypothesis on the
regularity and decay of $\tilde \psi_{\tilde \lambda}$. If $G = SO(2)$
then the proof can be directly adapted from the proof on the translation
group, by replacing convolutions with circular convolutions. 
It proves that $\widetilde S_L$ is Lipschitz continuous to the action 
of diffeomorphisms on $\Ld(SO(2))$.

\subsection{Combined Translation and Rotation Scattering}
\label{Combined} 

We construct a scattering operator on $\LD$ which is invariant to 
translations and rotations, by combining a 
translation invariant scattering operator on
$\LD$ and a rotation invariant scattering operator on $\Ld(SO(d))$.

Let $G$ be a rotation subgroup of the general linear group of $\R^d$,
which also includes
the reflection $-\One$ defined by $-\One x = -x$.
According to (\ref{motherandfosdf}), the wavelet
transform in $\LD$ is defined for any $\lambda = 2^j r \in 2^\Z \times G$ by
$W[\lambda] f = f \star \psi_{\lambda}$, where
$\psi_{\lambda} (x)= 2^{dj} \psi(2^{j} r^{-1} x)$.
Section \ref{scalittwnwavel} considers the case of a finite group $G$,
which is a subgroup of $SO(d)$ if $d$ is even or which is 
a subgroup of $O(d)$ if $d$
is odd, while including $-\One$.
The extension to a compact subgroup potentially equal to $SO(d)$
or $O(d)$ is straightforward.
We still denote $G^+$ the quotient of $G$ by $\{-\One,\One\}$.
The wavelet transform of a complex valued functions is defined
over all $\lambda \in 2^\Z \times G$ but it is restricted to
$2^\Z \times G^+$ if $f$ is real. 
All discrete sums on $G$ and $G^+$ are replaced by integrals with
the Haar measure $dr$. The group is compact and thus has a finite measure
$|G| = \int_G dr$. It results that these integrals behave as finite
sums in all derivations of this paper. The theorems
proved for a finite group $G$ remains valid for a compact group $G$.
In the following we concentrate
on real valued functions.

Let $\cP_J$ be the countable set of all finite paths
$p = (\la_1,...,\la_m)$ with 
$\la_k  \in \cLa_J = \{ \la = 2^j r~:~j > -J \,,\,r \in G^+ \}$.
The windowed scattering $S_J [\cP_J] f = \{S_J[p] f \}_{p \in \cP_J}$
is defined in Definition \ref{windscatdef}, but
$\cLa_J$ and $\cP_J$ are not discrete sets anymore.
The scattering norm is defined by summing the $\LD$ norms of all
$S_J[p] f$ for all $p =(2^{j_1} r_1,...,2^{j_m} r_m) \in \oP_J$, with
the Haar measure:
\[
\|S_J[\cP_J] f \|^2 = 
\sum_{m=0}^{\infty} \sum_{j_1 > -J,...,j_m > -J} \int_{\G^{+m}} 
\|S_J[2^{j_1} r_1,...,2^{j_m} r_m] f \|^2 \, d r_1...dr_m\,~,
\]
which is written
\[
\|S_J [\cP_J] f \|^2 = 
\int_{\cP_J} \|S_J[p] f \|^2\, dp~.
\]
One can verify that $S_J[\cP_J]$ is nonexpansive as in the case where $G$
is a finite group. For an admissible scattering wavelet satisfying 
(\ref{conditionprogre}), Theorem \ref{energydecth} remains valid and
$\|S_J [\cP_J] f \| = \|f\|$.

The scattering $S_J$ is covariant to rotations.
Invariance to rotations
in $\G$ is obtained by applying
the scattering operator $\widetilde S_L$ defined on $\Ld(G)$
by (\ref{relastn0898sdfn}).
Any $p = (\la_1,...,\la_m) \in \cP_J$ with $\la_1 = r 2^{j_1}$
can be written as a rotation $p = r\, \bar p$ of a
normalized path $\bar p = (\bar \la_1,...,\bar \la_m)$, where
$\bar \la_k = r^{-1} \la_k$ and hence where $\bar \la_1 = 2^{j_1}$ is a scaling
without rotation. It results that
\[
S_J [p] f(x) = S_J [r\, \bar p] f(x) ~.
\]
For each $x$ and $\bar p$ fixed, 
$S_J[r \,\bar p] f(x)$ is
a function of $r$ which belongs to $\Ld(G)$. 
We can thus apply the scattering operator 
$\widetilde S_L [\tilde p]$ to
this function of $r$. The result is denoted
$\widetilde S_L [\tilde p] S_J[r\, \bar p] f(x)$ for all 
$\tilde p \in \tP_L$. This output can be indexed by
the original path variable $p = r\, \bar p$, and we denote the 
combined scattering:
\begin{equation}
\label{combinesnf}
\widetilde S_L [\tilde p]\, S_J[p] f(x) = \widetilde S_L [\tilde p] 
\,S_J[r\, \bar p] f(x)~.
\end{equation}
This combined scattering
cascades wavelet transforms and hence convolutions along the spatial
variable $x$ and along the rotation $r$, 
which is factorized from each path.
In $d = 2$ dimensions then $\Ld(SO(2)) = \Ld [0,2\pi]$. The wavelet
tranform along rotations is implemented by circular convolutions along
the rotation angle variable in $[0,2\pi]$, with the periodic wavelets
(\ref{periwavsdfs}).

A combined scattering transform computes
\[
\widetilde S_L [\tP_L]  S_J[\cP_J]f = \{ \widetilde S_L [\tilde p] S_J[p] f\}_{p \in \cP_J, \tilde p \in \tP_L}~.
\]
Its norm is computed by summing the $\LD$ norms
$\|\widetilde S_L [\tilde p] S_J[p] f\|^2$:
\begin{equation}
\label{scanofnsd22}
\|\widetilde S_L [\tP_L ] S_J[\cP_J] f \|^2 = 
\sum_{\tilde p \in \tP_L} \int_{\cP_J} \|\widetilde S_L[\tilde p] S_J[p] f \|^2\, dp~.
\end{equation}
Since $\widetilde S_L[\tP_L]$ and $S_J[\cP_J]$ 
are nonexpansive their cascade is also nonexpansive:
\[
\forall (f,h)\in \LD^2~~,~~
\|\widetilde S_L [\tP_L] S_J[\cP_J] f - 
\widetilde S_L [\tP_L] S_J[\cP_J] h\| \leq \|f - h \|~.
\]
If $S_J$ is computed with an admissible scattering wavelet then
$\|S_J[\cP_J] f \| = \|f\|$. In $d = 2$ dimensions, if $\widetilde S_L$
is computed with periodic wavelets derived in (\ref{periwavsdfs}) 
from a one-dimensional admissible scattering wavelet, 
then the combined scattering preserves the norm:
\[
\|\widetilde S_L [\tP_L] S_J[\cP_J] f \| = 
\| S_J[\cP_J] f \| = \|f\|~.
\]

By setting $L = 0$ and letting $J$ go to $\infty$, the following
proposition proves that the resulting combined scattering is
invariant to translations and rotations.
Such combined
scattering representations are used for rotation invariant classification
of image textures \cite{Sifre}.
For any $(c,g) \in \R^d \times SO(d)$, 
we denote $L_{c,g} f(x) = f(g^{-1}(x - c))$.

\begin{proposition}
\label{groupinfapropinvarosn}
For all $(c,g) \in \R^d \times SO(d)$ and all $f \in \LD$
\begin{equation}
\label{dispanfdson0sdf8}
\forall (\tilde p, p) \in \widetilde \cP_0 \times \cPinf\,\,,\,\, 
\widetilde S_0[\tilde p] \oS(p) L_{c,g} f =  
\widetilde S_0[\tilde p] \oS(p) f ~.
\end{equation}
\end{proposition}

{\it Proof:} The scattering transform in $\LD$ 
is translation invariant and covariant to rotations:
$\oS L_{c,g} f (p) = \oS f(g^{-1} p)$ for all $p \in \cPinf$.  
Since $g^{-1}$ acts as a rotation on the path $p$, 
Proposition \ref{groupinfaprop} proves that 
$\widetilde S_0[\tilde p] \oS f (g^{-1}\,p) = 
\widetilde S_0[\tilde p] \oS f (p)$, which gives (\ref{dispanfdson0sdf8}).
$\Box$

\appendix

\section{Proof of Lemma \protect \ref{lemma8}}
\label{proofLemma8}

The proof of (\ref{condroref}) shows that
the scattering energy propagates towards lower frequencies.
It computes the average arrival log frequency of the
scattering energy $\|U [p] f \|$ for paths of length $m$, 
and shows that it increases when $m$ increases.
The arrival log frequency of
$p = \{ \la_k = 2^{j_k} r_k \}_{k \leq m}$ 
is the log frequency index $\log_2 |\la_m| = j_m$ of the last path element.

Let us denote $e_m = \| U[\cLa_J^m] f\|^2$ and
$\overline e_{m} = \| S_J [\cLa_J^m] f\|^2$.
The average arrival log frequency
among paths of length $m$ is
\begin{equation}
\label{dnslnfs0fd8}
{a_m} = e_m^{-1}
\,\sum_{p \in \cLa^m_J} {j_{m}}\, \|U[p] f\|^2 \geq -J~.
\end{equation}
The following lemma
shows that when $m$ increases by $1$ then $a_m$ decreases by nearly $\alpha/2$,
where $\alpha$ is defined in (\ref{conditionprogre}).

\begin{lemma}
\label{wndsofnsw0}
If (\ref{conditionprogre}) is satisfied then
\begin{equation}
\label{dnslnfs0}
\forall m > 0~~,~~
\frac {\alpha} 2 
\, e_{m-1} \leq (a_m+J) e_m  - (a_{m+1}+J)  e_{m+1} + 
e_{m-1} - e_{m}~.
\end{equation}
\end{lemma}

We first show that (\ref{dnslnfs0}) implies
(\ref{condroref}) and then prove this lemma.
Summing over (\ref{dnslnfs0}) gives
\begin{equation}
\label{infds8sdf}
\frac{\alpha} 2 \sum_{k=0}^{m-1} e_k \leq (a_1+J) e_1 - (a_{m+1}+J) 
e_{m+1} +e_{0} - e_{m} \leq e_{0} + (a_1+J) e_1~.
\end{equation}
For $m=1$, $p=2^j\ga$ so
$a_1\,e_1 = \sum_{j>-J} \sum_{\ga \in {\G^+}}  j\,\|W[{2^j\ga}] f \|^2 $.
Moreover, $e_0 = \|f\|^2$, so
\[
e_{0} + (a_1+J) e_1 = \|f  \|^2 +
\sum_{j>-J} \sum_{\ga \in {\G^+}} (j+J)\,\|W[{2^j\ga}] f \|^2 ~.
\]
Inserting this in (\ref{infds8sdf}) for $m = \infty$
proves (\ref{condroref}).

Lemma \ref{wndsofnsw0} is proved by calculating 
the evolution of $a_m$ as $m$ increases.
We consider the advancement of a path $p$ 
of length $m-1$ with
two steps $p+2^j\ga + 2^lr'$, and denote
$f_p = U[p] f$.  
The average arrival log frequency $a_m$ can be written
as the average arrival log frequency
of $\|U[p+2^j\ga] f\|^2$ over
all $2^j\ga$ and all $p$ of length $m-1$:
\begin{equation}
\label{jmsfn}
a_m\, e_m
= \sum_{p\in \cLa_J^{m-1}}\sum_{j > -J} \sum_{\ga \in {\G^+}} j\,
\|f_p \star \psi_{2^j\ga} \|^2 ~.
\end{equation}
After the second step, the average arrival log frequency of
$\|U[p+2^j\ga+2^lr']f\|^2$ overall $p \in \cLa^{m-1}_J$, $2^j\ga$ and 
$2^lr'$ is $a_{m+1}$:
\[
a_{m+1}\, e_{m+1}
= \sum_{p\in \cLa^{m-1}_J}\sum_{j > -J} \sum_{\ga \in {\G^+}} 
\sum_{l > -J} \sum_{\ga' \in {\G^+}} 
l\, \||f_p \star \psi_{2^j\ga}| \star\psi_{2^l\ga'}  \|^2 \,.
\]
The wavelet transform is unitary and hence for any $h \in \LD$
\[
\|h\|^2 =  \sum_{l > -J} \sum_{\ga' \in {\G^+}}  \|h \star\psi_{2^l\ga'}  \|^2 + \|h \star \phi_{2^J}\|^2 ~. 
\]
Applied to each $h = f_p \star \psi_{2^j\ga}$ in (\ref{jmsfn}) this relations,
together with 
\[
\overline e_{m} = \sum_{p \in \cLa^{m-1}_J}\sum_{j > -J} \sum_{\ga \in \G^+} 
\||f_p \star \psi_{2^j\ga}| \star\phi_{2^J}  \|^2 dr, 
\]
shows that $I = a_m \,e_m - a_{m+1}\, e_{m+1}  + J\,  \overline e_{m} $ satisfies
\begin{eqnarray*}
I&=& \sum_{p \in \cLa^{m-1}_J}
\sum_{j > -J} \sum_{\ga' \in {\G^+}}
\Bigl(\sum_{l > -J} \sum_{\ga \in {\G^+}}
(j-l) \, \||f_p \star \psi_{2^j\ga}| \star\psi_{2^l\ga'}  \|^2 
\\
& &  ~~~~~~~~~~~~~+ ( j+J) \,
  \||f_p \star \psi_{2^j\ga}| \star \phi_{2^J}  \|^2 \Bigr) ~.
\end{eqnarray*}
A lower bound of $I$ is calculated by dividing the sum on $l$ for
$l \geq j$ and $l < j$. In the $j+J-1$ term
for $l < j$, the index $l$ is replaced by $j-1$ and the convolution with
$\phi_{2^J}$ is incorporated in the sum:
\begin{eqnarray}
I
&\geq & \sum_{p \in \cLa^{m-1}_J} 
\sum_{j > -J} \sum_{\ga' \in {\G^+}} \left[
\sum_{-J < l<j} \left( \sum_{\ga \in {\G^+}}
\| |f_p \star \psi_{2^j\ga}|  \star \psi_{2^l\ga'} \|^2  \right) \right.
\nonumber \\
\label{convpdsfjsdf8sd}
& & +  \| |f_p \star \psi_{2^j\ga}| \star \phi_{2^J}\|^2  - \left. 
\sum_{l>j} \sum_{\ga \in {\G^+}}
(l-j)  \, \| |f_p \star \psi_{2^j\ga}|  \star \psi_{2^l\ga'} d \ga'\|^2 
\right] d\ga~.
\end{eqnarray}

Since wavelets satisfy the unitary property (\ref{consenas}),
for all real functions $f \in \LD$ and all $q \in \Z$,
\begin{equation}
\label{conservq}
\sum_{-q \geq l > -J} \sum_{\ga \in \G^+}
\|f \star \psi_{2^l\ga} \|^2 + \|f \star \phi_{2^J} \|^2 = \|f \star \phi_{2^q} \|^2~.
\end{equation}
Indeed (\ref{consenas}) implies that 
\begin{equation}
\label{conservq98dsf}
|\hat \phi(2^{J} \om)|^2 + \frac 1 2 
\sum_{-q \geq l > -J} \sum_{\ga \in \G}
|\hat \psi (2^{-l} \ga^{-1} \om)|^2
= |\hat \phi(2^{q} \om)|^2~.
\end{equation}
If $f$ is real, then $\|f \star \psi_{2^j r}\| = \|f \star \psi_{-2^j r}\|$.
Multiplying (\ref{conservq98dsf})  by $|\hat f(\om)|^2$ and
integrating in $\omega$ proves (\ref{conservq}).
Inserting (\ref{conservq}) in (\ref{convpdsfjsdf8sd}) gives
\begin{eqnarray*}
I &\geq& \sum_{p \in \cLa^{m-1}_J} 
\sum_{j > -J} \sum_{\ga \in \G^+} \left(
\||f_p \star \psi_{2^j\ga}|  \star \phi_{2^{-j+1}} \|^2 
\right. \\
& & - 
\sum_{l>j} 
(l-j) \,\Bigl( \| |f_p \star \psi_{2^j\ga}|  \star \phi_{2^{-l}} \|^2 -
 \| |f_p \star \psi_{2^j\ga}|  \star \phi_{2^{-l+1}} \|^2 \Bigr) 
\Bigr) \, 
~.
\end{eqnarray*}
If $\rho \geq 0$ satisfies $|\hat \rho (\om)| \leq |\hat \phi (2 \om)|$,
then for any $f \in \LD$ and any $l \in \Z$,
\[
 \| f  \star \phi_{2^{-l+1}} \|^2 \geq  \| f  \star \rho_{2^l \ga} \|^2 ~~
\mbox{with}~~\rho_{2^l\ga} (x) = 2^{dl} \rho (2^{l} \ga^{-1} x)~.
\]
It results that
\begin{eqnarray*}
I &\geq &\sum_{p \in \cLa^{m-1}_J} 
\sum_{j > -J} \sum_{\ga \in \G^+}  \left(
\||f_p \star \psi_{2^j\ga}|  \star \rho_{2^j\ga} \|^2 \right. \\
& & - 
\sum_{l>j} 
(l-j) \,\Bigl( \| f_p \star \psi_{2^j\ga} \|^2 -
 \| |f_p \star \psi_{2^j\ga}|  \star \rho_{2^l \ga} \|^2 \Bigr) \Bigr) ~.
\end{eqnarray*}
Applying Lemma \ref{boundendf} for $h = \rho_{2^l \ga}$ and a frequency
$2^{j} \ga \eta$ proves that
\[
 \| |f_p \star \psi_{2^j\ga}|  \star \rho_{2^l \ga} \| \geq
 \| f_p \star \psi_{2^j\ga}  \star \rho_{2^l \ga,2^j} \| ~~
\mbox{with}~~\rho_{2^l \ga,2^j} (x) = \rho_{2^l \ga} (x) \, e^{i 2^{j} \ga \eta . x }~~
\]
and $\hat \rho_{2^l \ga,2^j} (\om) = \hat \rho (2^{-l} \ga^{-1} \om - 2^{j-l} \eta)$.
It results that
\begin{eqnarray*}
I &\geq &\sum_{p \in \cLa^{m-1}_J} 
\sum_{j > -J} \sum_{\ga \in \G^+}  \left(
\|f_p \star \psi_{2^j\ga}  \star \rho_{2^j\ga,2^j} \|^2 \right. \\
& & - 
\sum_{l>j} 
(l-j) \,\Bigl( \| f_p \star \psi_{2^j\ga} \|^2 -
 \| f_p \star \psi_{2^j\ga}  \star \rho_{2^l\ga,2^j} \|^2 \Bigr) \Bigr)\,~.
\end{eqnarray*}

We shall now rewrite this equation in the Fourier domain.
Since $f_p(x) \in \R$, $|\hat f_p (\om)| = |\hat f_p (-\om)|$,
applying Plancherel gives
\begin{eqnarray*}
I &\geq &\frac 1 2 \sum_{p \in \cLa^{m-1}_J} 
\int |\hat f_p (\om)|^2 \sum_{\ga \in \G} \sum_{j > -J} \left( 
|\hat \psi (2^{-j} \ga^{-1}\om)|^2\,|\hat \rho (2^{-j} \ga^{-1} \om - \eta)|^2   \right. \\
& & ~~~~~~~~~ 
- \sum_{l>j} (l-j) \,
|\hat \psi (2^{-j} \ga^{-1} \om)|^2\,(1 - |\hat \rho (2^{-l} \ga^{-1} \om - 2^{j-l}\eta)|^2)
\, \Bigr)\,d\om~.
\end{eqnarray*}
Inserting $\hat \Psi$ defined in (\ref{conditionprogredf0}) by
\[
\hat \Psi (\omega) = 
| \hat \rho (\om-\eta)|^2 - \sum_{k=1}^{+\infty} 
k\, (1 - |\hat \rho (2^{-k} (\om - \eta))|^2)~
\]
with $k = l-j$, gives
\[
I \geq \frac 1 2 
\sum_{p \in \cLa^{m-1}_J} 
\int |\hat f_p (\om)|^2
\sum_{j > -J} b (2^{-j} \omega)\,d\omega
\]
with $b( \omega) =  \sum_{\ga \in \G}
\hat \Psi( \ga^{-1} \om ) \,|\hat \psi ( \ga^{-1} \om)|^2 \,$.
Let us add to $I$
\[
\overline e_{m-1} = \sum_{p \in \cLa^{m-1}_J} 
\|f_p \star  \phi_{2^J}  \|^2  = 
\sum_{p \in \cLa^{m-1}_J} \int |\hat f_p (\om)|^2\,|\hat \phi(2^{J}\om)|^2\,d\om~.
\]
Since $\rho \geq 0$, $|\hat \rho (\om)| \leq \hat \rho(0) = 1$ 
and hence $\hat \Psi(\om) \leq 1$.
The wavelet unitary property (\ref{consenas}) 
together with $\hat \Psi (\om) \leq 1$ implies that
\[
|\hat \phi(2^{J} \om)|^2 = \frac 1 2 \,\sum_{j \leq -J} \sum_{\ga \in \G}
|\hat \psi (2^{-j} \ga^{-1} \om)|^2   \geq \frac 1 2\,
\sum_{j \leq -J} b (2^{-j} \omega)
\]
so
\[
I + \overline e_{m-1} \geq \,\frac 1 2 \,\sum_{p \in \cLa^{m-1}_J} 
\int |\hat f_p (\om)|^2 
\sum_{j =-\infty}^{+\infty} b (2^{-j} \omega)\,d\om~.
\]

If $\alpha = \inf_{1\leq |\omega| < 2} \sum_{j} b (2^{-j} \omega)$
then $\sum_{j} b (2^{-j} \omega) \geq \alpha$ for all
$\omega \neq 0$.
If the hypothesis (\ref{conditionprogre}) is satisfied 
and hence $\alpha > 0$ then
\begin{eqnarray*}
I + \overline e_{m-1}&\geq& \frac {\alpha} 2 \, \sum_{p \in \cLa^{m-1}_J} 
\int |\hat f_p (\om)|^2 \, d\om = \frac \alpha 2 \sum_{p \in \cLa^{m-1}_J} 
\|f_p \|^2\\
& = &\frac \alpha 2 \, \sum_{p \in \cLa^{m-1}_J} \|U[p] f \|^2  = \frac \alpha 2 \,e_{m-1} ~.
\end{eqnarray*}
Inserting $I = a_m \,e_m - a_{m+1}\, e_{m+1}  + J\,  \overline e_{m} $ 
proves that
\begin{equation}
\label{dnslnfs0fe}
a_m \,e_m - a_{m+1}\, e_{m+1}  + 
  J\,  \overline e_{m} + \overline e_{m-1} \geq
\frac \alpha 2 \,e_{m-1}~.
\end{equation}

Since $U_J$ preserves the norm,
$e_m = e_{m+1} + \overline e_{m}$,
indeed (\ref{propasndfs}) proves that
$\oUJ U[\cLa_J^m] f  =  \{U[\cLa_J^{m+1}]f\,,\,S_J[\cLa_J^{m}] f\}$.
Inserting $\overline e_{m} = e_{m} - e_{m+1}$ and $\overline e_{m-1} = e_{m-1} - e_{m}$ 
in (\ref{dnslnfs0fe})  gives
\[
\frac \alpha 2 \, e_{m-1} \leq (a_m+J) e_m  - (a_{m+1}+J)  e_{m+1} + e_{m-1} - e_{m}~,
\]
which finishes the proof of Lemma \ref{wndsofnsw0}.

\section{Proof of Lemma \protect \ref{lemma1}}
\label{proofLemma1}

Lemma \ref{lemma1} as well as all other upper bounds on 
operator norms are computed with Schur's lemma. 
For any
operator $K f (x) = \int f(u)\,k(x,u)\,du$, Schur's lemma proves that
\begin{equation}
\label{schs1eq2}
\int |k(x,u)|\,dx \leq C~~\mbox{and}~~\int |k(x,u)|\,du \leq C~~
{\Longrightarrow}~~\|K \| \leq C~,
\end{equation}
where $\|K\|$ is the $\LD$ norm of $K$.

The operator norm of $k_J = L_\tau A_J - A_J$ is computed by applying
Schur's lemma on its kernel
\begin{equation}
\label{kernappA}
k_J (x,u) = \phi_{2^J} (x-\tau(x) - u) - \phi_{2^J} (x - u)~.
\end{equation}
A first-order Taylor expansion proves that
\[
|k_J (x,u)| \leq 
 |\int_0^1 \nabla \phi_{2^J} (x-u- t\,\tau(x))\,.\,\tau(x)\, dt| \leq
\|\tau\|_\infty \, \int_0^1 
| \nabla \phi_{2^J} (x-u- t\,\tau(x))|\, dt 
\]
so
\begin{equation}
\label{asdfnwds}
\int |k_J(x,u)| \,du \leq \|\tau\|_\infty \, \int_0^1 \int 
| \nabla \phi_{2^J} (x-u- t\,\tau(x))|\, du\, dt ~.
\end{equation}
Since $\nabla \phi_{2^J} (x) = 2^{-dJ-J} \nabla \phi (2^{-J} x)$, it results that
\begin{equation}
\label{asdfnwds8df7}
\int |k_J(x,u)| \,du \leq 
\|\tau\|_\infty \, 2^{-dJ-J} 
\int | \nabla \phi (2^{-J} u')|\, du' = 
2^{-J}\,\|\tau\|_\infty\,\|\nabla \phi \|_1~.
\end{equation}
Similarly to (\ref{asdfnwds}) we prove that
\[
\int |k_J(x,u)| \,dx 
\leq \|\tau\|_\infty \, \int_0^1 \int 
| \nabla \phi_{2^J} (x-u- t\,\tau(x))| \,dx\, dt ~.
\]
The Jacobian of the
change of variable $v = x - t \,\tau (x)$ is $\Id - t \nabla \tau (x)$ whose
determinant is larger than $(1 -  \|\nabla \tau \|_\infty)^d \geq 2^{-d}$ 
so
\begin{eqnarray*}
\int |k_J(x,u)| \,dx 
&\leq& \|\tau\|_\infty \, 2^{d}\,
\int_0^1 \int 
| \nabla \phi_{2^J} (v-u)| \,dv\, dt\\
& =& 2^{-J}\,\|\tau\|_\infty\,
\|\nabla \phi \|_1
\, 2^{d}\,.
\end{eqnarray*}
Schur's lemma (\ref{schs1eq2}) applied to this upper bound
and (\ref{asdfnwds8df7}) proves the lemma result:
\[
\|L_\tau A_J - A_J \| \leq 2^{-J+d} \,\|\nabla \phi \|_1\,
\,\|\tau\|_\infty~.
\]

\section{Proof of  (\ref{err10dfi01})}
\label{proofLemma2}

We prove that
\begin{equation}
\label{lemma1eq3}
\|L_\tau A_J  f - A_J f  + \tau \, . \,\nabla A_J f\| \leq  C\,\|f\|\,
2^{-2J}~\, \|\tau\|^2_\infty ~
\end{equation}
by applying Schur's lemma (\ref{schs1eq2}) on the
kernel of $k_J = L_\tau A_J  - A_J  + \tau \cdot \nabla A_J$:
\[
k_J (x,u) = \phi_{2^J} (x-\tau(x) - u) - \phi_{2^J} (x - u) + 
\nabla \phi_{2^J}  (x-u) \cdot \tau(x) ~.
\]

Let $\rH f(x)$ the Hessian matrix of a function $f$ at $x$ and
$|\rH f(x)|$ the sup matrix norm of this Hessian matrix. 
A Taylor expansion gives
\begin{eqnarray}
|k_J(x,u)| &=&
\left|\int_0^1 t \tau(x) \cdot \rH \phi_{2^J} (u-x-(1-t)\tau(u)) \cdot \tau (x) \, dt
\right|\nonumber\\
\label{asdfnwds2}
&\leq&  \|\tau \|^2_\infty \, \int_0^1 |t| 
| \rH\phi_{2^J} (u-x- (1-t)\,\tau(x))|\, dt ~.
\end{eqnarray}
Since $\phi_{2^J}(x) = 2^{-dJ} \phi(2^{-J} x)$, 
$\rH \phi_{2^J} (x) = 2^{-Jd-2J} \rH \phi (2^{-J} x)$. With a change
of variable, (\ref{asdfnwds2}) gives
\begin{equation}
\label{fndosfinwe}
\int |k_J(x,u)| \,du \leq 
\|\tau \|^2_\infty \, 2^{-dJ-2J} 
\int | \rH \phi (2^{-J} u')|\, du' = 
2^{-2J}\,\|\tau \|^2_\infty\,\|\rH \phi \|_1~,
\end{equation}
where $\|\rH \phi \|_1= \int |\rH \phi(u)|\,du$ is bounded.
Indeed all second-order derivatives of $\phi$ at $u$ are $O((1+|u|)^{-d-1})$.

The upper bound (\ref{asdfnwds2}) also implies that
\[
\int |k_J(x,u)| \,dx 
\leq \|\tau \|^2_\infty \, \int_0^1 |t| \int 
|\rH \phi_{2^J} (u-x- (1-t)\tau(x))| \,du\, dt ~.
\]
The Jacobian of the
change of variable $v = x - (1-t) \,\tau (x)$ is $\Id - (1-t) \nabla \tau (x)$ whose determinant is larger than $(1 - \|\nabla \tau \|_\infty)^d$ so
\begin{eqnarray}
\int |k_J(x,u)| \,dx 
&\leq& \|\tau \|^2_\infty \, (1 -  \|\nabla \tau \|_\infty)^{-d}\,
\int_0^1 \int 
|\rH \phi_{2^J} (v-u)| \,dv\, dt \nonumber\\
\label{fndosfinwe2}
& =& 2^{-2J}\,\|\tau \|^2_\infty\,
\|\rH \phi \|_1
\, 2^{d}\,.
\end{eqnarray}
The upper bounds (\ref{fndosfinwe}) and (\ref{fndosfinwe2})
with Schur's lemma (\ref{schs1eq2}) proves (\ref{lemma1eq3}).

\section{Proof of Lemma \protect\ref{LemmaComsnuta}}
\label{Proofsdfn}

This appendix proves that for any operator $L$ 
and any $f \in \LD$
\begin{equation}
\label{infds3}
\| [S_J[\cP_J]\,,\, L] f \|
 \leq \|[ U_J\,,\,L]\|\, \|U[\cP_J] f \|_{1} =
\|[ U_J\,,\,L]\|\, 
\sum_{n=0}^{\infty} \|U[\cLa^n_J] f \|.
\end{equation}
If $A$ and $B$ are two operators, we denote
$\{A,B\}$ the operator defined by $\{A,B\}f = \{A f, Bf\}$. 
We introduce a wavelet modulus operator without averaging:
\begin{equation}
\label{onesnfs23}
V_J f = \{|W[{\lambda}] f| = |f \star \psi_{\lambda}|\}_{\la \in \cLa_J}~~
\mbox{with}~~\cLa_J = \{2^j\ga~:~j > -J~,\ga \in \G^+ \},
\end{equation}
and $U_J = \{A_J \,,\,V_J \}$.
The propagator $V_J$ creates all paths
$V_J U [\cLa^n_J] f = U[\cLa^{n+1}_J] f$
for any $n \geq 0$. 
Since $U[\cLa^0_J] = Id$, it results that $V_J^n = U[\cLa^n_J]$. 
Let $\cP_{J,m}$ be the subset of $\cP_J$ of paths $p$ of 
length smaller than $m$.
To verify (\ref{infds3}), we shall prove that 
\begin{equation}
\label{infds38}
[S_J[\cP_{J,m}], L]   =\sum_{n=0}^m K_{m-n} V_J^n ~,
\end{equation}
where $K_n = \{[A_J,L] \,,\,S_J[\cP_{J,n-1}] [V_J,L]  \} $
satisfies 
\begin{equation}
\label{infds39}
\|K_n  \| \leq \|[ U_J , L ] \| \,.
\end{equation}
Since $V_J^n f = U[\cLa^n_J] f$,  
it implies that for any $f \in \LD$
\[
\|[S_J[\cP_{J,m}], L] f \| \leq \sum_{n=0}^m \|K_{m-n} \|\, \| V_J^n f \|
\leq \|[ U_J , L ] \| \, \sum_{n=0}^{m-1} \|U[\cLa^n_J] f \|~,
\]
and letting $m$ tend to $\infty$ proves (\ref{infds3}).

Property (\ref{infds38}) is proved by first showing that
\begin{equation}
\label{infds309}
S_J[\cP_{J,m}] L  =\{ L A_J\,,\,S_J[\cP_{J,m-1}] L V_J\} + K_m~,
\end{equation}
where $K_m = \{[A_J,L] \,,\,S_J[\cP_{J,m-1}] \,[V_J,L]  \} $.
Indeed, since $V_J^n = U[\cLa^n_J]$, we have
$A_J\,V_J^n  = S_J[\cLa^{n}_J]$ and $\cP_{J,m} = \cup_{n=0}^{m-1} \cLa^n_J$
yields
$S_J[\cP_{J,m}]  = \{A_J\, V_J^n   \}_{0 \leq n < m}$.
It results that
\begin{eqnarray*}
S_J[\cP_{J,m}] L  &=& \{A_J\, V_J^n L   \}_{0 \leq n < m} \\
&=& \{L A_J  + [A_J,L]\,,\, A_J V_J^{n-1} L V_J  + A_J V_J^{n-1} [V_J,L] \}_{1 \leq n < m} \\
&=&\{L A_J \,,\,S_J[\cP_{J,m-1}] L V_J  \} +
\{[A_J,L] \,,\,S_J[\cP_{J,m-1}]\, [V_J,L]  \} \\
&=&\{L A_J \,,\,S_J[\cP_{J,m-1}] L V_J  \} + K_m~,
\end{eqnarray*}
which proves (\ref{infds309}).

A substitution of $S_J[\cP_{J,m-1}] L$ in (\ref{infds309}) by
the expression derived by this same formula  gives
\[
S_J[\cP_{J,m}] L  =\{ L A_J\,,\,L A_J V_J\,,\,
S_J[\cP_{J,m-2}] L V_J^2\} + K_{m-1} V_J + K_m~.
\]
With $m$ substitions, we obtain
\[
S_J[\cP_{J,m}] L  =\{ L A_J V_J^n\}_{0 \leq n < m} + \sum_{n=0}^m K_{m-n} V_J^n =
L S_J[\cP_{J,m}]  + \sum_{n=0}^m K_{m-n} V_J^n ~
\]
which proves (\ref{infds38}).

Let us now prove (\ref{infds39}) on
$K_m = \{[A_J,L] \,,\,S_J[\cP_{J,m-1}]\, [V_J,L]  \} $.
Since $S_J[\cP_J]$ is  nonexpansive, its restriction
$S_J[\cP_{J,m}]$ is also nonexpansive. Given that
$U_J = \{A_J \,,\,V_J \}$ we get
\begin{eqnarray*}
\|K_m f \|^2 &=& \|[A_J,L] f \|^2 + \|S_J[\cP_{J,m-1}]\, [V_J,L]f\|^2\\
& \leq & \|[A_J,L] f \|^2 + \|[V_J,L]f\|^2 = \|[ U_J , L]f\|^2 \leq
\|[ U_J , L]\|^2 \,\|f\|^2
\end{eqnarray*}
which proves (\ref{infds39}).

\section{Proof of Lemma \protect\ref{thedw0}}
\label{prooftheosdfn} 

This section computes an upper bound of $\|[W_J , L_\tau]\|$ by
considering
\[
[W_J , L_\tau]^*\,[W_J , L_\tau] =
\sum_{\ga \in \G^+} \sum_{j=-J+1}^{\infty}
[W[2^j\ga] , L_\tau]^*[W[2^j\ga] , L_\tau]
+ [A_J , L_\tau]^*\,[A_J , L_\tau].
\]
Since $\|[W_J , L_\tau]\| = \|
[W_J , L_\tau]^*\,[W_J , L_\tau]\|^{1/2}$, 
\begin{equation}
\label{cnvondf09s}
\|[W_J , L_\tau]\| \leq
\sum_{\ga \in \G^+} 
\Big\|\sum_{j=-J+1}^{\infty} [W[2^j\ga] , L_\tau]^*[W[2^j\ga] , L_\tau]\Big\|^{1/2}\, 
+ \|[A_J , L_\tau]^*\,[A_J , L_\tau]\|^{1/2}~.
\end{equation}
To prove the upper bound (\ref{commnaf}) of Lemma \ref{thedw0}, 
we compute an upper bound for each term on the right under the
integral and the last term, which is done
by the following lemma.

\begin{lemma}
\label{ndfondflemma}
Suppose that $h(x)$, as well as all
its first and second-order derivatives have a decay in
$O((1 + |x|)^{-d-2})$. Let $Z_j f = f \star h_j$ 
with $h_j (x) = 2^{dj} h(2^{j} x)$. 
There exists $C > 0$ such that if $\|\nabla \tau \|_\infty \leq$
then
\begin{equation}
\label{firsnfidfnw}
\|[Z_j , L_\tau ]\| \leq C\, \|\nabla \tau \|_\infty
\end{equation}
and if $\int h(x)\,dx = 0$ then 
\begin{equation}
\label{firsnfidfnw2}
\Big\|\sum_{j=-\infty}^{+\infty} [Z_{j} , L_\tau]^*[Z_{j} , L_\tau]\Big\|^{1/2} \leq
C\,\Bigl( \max\Bigl(\log \frac {\|\Delta \tau\|_\infty}{\|\nabla \tau\|_\infty}\,,\, 1\Bigr)\,\|\nabla \tau \|_\infty + \|\rH \tau \|_\infty \Bigr)~.
\end{equation}
\end{lemma}

The inequality (\ref{firsnfidfnw2}) clearly remains valid if the summation
is limited to $-J$ instead of $-\infty$ since 
$[Z_{j} , L_\tau]^*[Z_{j} , L_\tau]$ is a positive operator. Inserting 
in (\ref{cnvondf09s}) both
(\ref{firsnfidfnw}) with $h = \phi$ and 
(\ref{firsnfidfnw2}) with $h(x) = \psi(\ga^{-1} x)$ for each $\ga \in \G^+$, 
and replacing $-\infty$ by $-J$
proves the upper bound (\ref{commnaf}) of Lemma \ref{thedw0}. 

To prove Lemma \ref{ndfondflemma}, 
we factorize
\[
[Z_j,L_\tau]
 = K_j\, L_\tau~~\mbox{with}~~K_j = Z_j - L_\tau \,Z_j\,L_\tau^{-1}~.
\]
Observe that
\begin{equation}
\label{adsfnwonfs0}
\|[Z_j,L_\tau]^*[Z_j,L_\tau] \|^{1/2} = \|L_\tau^* K_j^* K_j L_\tau \|^{1/2} \leq 
\|L_\tau \|\,\| K_j^*\,K_j\|^{1/2} ,
\end{equation}
and that
\begin{equation}
\label{adsfnwonfs}
\Big\|\sum_{j=-\infty}^{+\infty} [Z_{j} , L_\tau]^*[Z_{j} , L_\tau]\Big\|^{1/2} \leq
\|L_\tau \|\,\|\sum_{j=-\infty}^{+\infty} K_j^*\,K_j\|^{1/2}~,
\end{equation}
with $\|L_\tau \| \leq (1 - \|\nabla \tau\|_\infty)^{-d} $.
Since $L_\tau^{-1} f(x) = f(\xi(x))$ with
$\xi(x-\tau(x)) = x$,
the kernel of $K_j = Z_j - L_\tau \,Z_j\,L_\tau^{-1}$ is 
\begin{equation}
\label{epsilonfdnwd0}
k_j (x,u) = h_j (x-u) - h_j (x - \tau(x) - u+\tau(u)) \,
\det(\Id - \nabla \tau(u))~.
\end{equation}

The lemma is proved by computing upper bounds
of $\|K_j\|$ and $\|\sum_{j=-\infty}^{+\infty} K_j^*\,K_j\|$.
The sum over $j$ is divided in three parts
\begin{equation}
\label{threpasfns}
\|\sum_{j=-\infty}^{+\infty} K_j^*\,K_j\|^{1/2} \leq
\|\sum_{j=-\infty}^{-\gamma-1} K_j^*\,K_j\|^{1/2} +
\|\sum_{j=-\gamma}^{-1} K_j^*\,K_j\|^{1/2} +
\|\sum_{j=0}^{\infty} K_j^*\,K_j\|^{1/2} ~.
\end{equation}
We shall first prove that
\begin{equation}
\label{toprove1}
\|\sum_{j=-\infty}^{-\gamma} K_j^* K_j\|^{1/2} \leq
C\,  \Bigl( \|\nabla \tau \|_\infty + 2^{-\gamma}\,\|\Delta \tau\|_\infty  
+ 2^{-\gamma/2}\,\|\Delta \tau \|^{1/2}_\infty  \,\|\nabla \tau \|_\infty^{1/2} \Bigr)~.
\end{equation}
Then we verify that $\|K_j\| \leq C\,\|\nabla \tau\|_\infty$ and hence that
\begin{equation}
\label{toprove2}
\|\sum_{j=-\gamma}^{-1} K_j^*\,K_j\|^{1/2} \leq \gamma\, \|K_j\| \leq 
C\,\gamma \, \|\nabla \tau\|_\infty ~.
\end{equation}
The last term carries the singular part and we prove that
\begin{equation}
\label{toprove3}
\|\sum_{j=0}^{\infty} K_j^*\,K_j\|^{1/2} \leq
C\, (\|\nabla \tau\|_\infty + \|\rH \tau \|_\infty)~.
\end{equation}
Choosing $\gamma = \max(\log \frac{\|\Delta \tau \|_\infty} {\|\nabla \tau\|_\infty}  , 1)$ 
yields
\[
\|\sum_{j=-\infty}^{+\infty} K_j^* K_j\|^{1/2} \leq
C\, \Bigl(\max\Bigl(\log \frac {\|\Delta \tau\|_\infty}{\|\nabla \tau\|_\infty}\,,\, 1\Bigr)\,\|\nabla \tau \|_\infty + \|\rH \tau \|_\infty \Bigr)~.
\]
Inserting this result in (\ref{adsfnwonfs}) will
prove the second lemma result (\ref{firsnfidfnw2}) .
In the proof, $C$ is a generic constant which depends only on 
$h$ but which evolves along the calculations. 

The proof of (\ref{toprove1}) is done by decomposing
$K_j = \tilde K_{j,1} + \tilde K_{j,2}$, with a first kernel
\begin{equation}
\label{epsilonfdnwd011}
\tilde k_{j,1} (x,u) = 
a(u)\, h_j (x-u)~~\mbox{with}~~a(u) = (1- \det(\Id - \nabla \tau(u)))\,,
\end{equation}
and a second kernel 
\begin{equation}
\label{epsilonfdnwd012}
\tilde k_{j,2} (x,u) = \det(\Id - \nabla \tau(u)) \Bigl(
h_j (x-u) - h_j (x - \tau(x) - u+\tau(u)) \Bigr)~.
\end{equation}
This kernel 
has a similar form as the kernel (\ref{kernappA}) in Appendix \ref{proofLemma1}
by $\tau(x)$ is replaced here by $\tau(x) - \tau(u)$.
The same proof shows that 
\begin{equation}
\label{upperboundd40dfd9}
\| \tilde K_{j,2} \| \leq 
C\, 2^{j}\,\|\Delta \tau\|_\infty ~.
\end{equation}
Taking advantage of this decay,
to prove (\ref{toprove1}), we decompose 
\begin{eqnarray*}
\|\sum_{j=-\infty}^{-\gamma} K_j^* K_j\|^{1/2} &\leq&
\| \sum_{j=-\infty}^{-\gamma} \tilde K_{j,1}^* \tilde K_{j,1} \|^{1/2}\\
& & + 
\sum_{j=-\infty}^{-\gamma} (\|\tilde K_{j,2}\|+ 2^{1/2}\,\| \tilde K_{j,2}\|^{1/2}  
\| \tilde K_{j,1}\|^{1/2})~
\end{eqnarray*}
and verify that 
\begin{equation}
\label{upperboundd40dfd8}
\|\tilde K_{j,1}\|\leq C\,  \|\nabla \tau \|_\infty ~~\mbox{and}~~
\| \sum_{j=-\infty}^{0} \tilde K_{j,1}^* \tilde K_{j,1} \|^{1/2} \leq 
C\,  \|\nabla \tau \|_\infty ~.
\end{equation}
The kernel of the self-adjoint operator $\tilde K_{j,1}^* \tilde K_{j,1}$  is:
\[
\tilde k_{j} (y,z) = \int \tilde k_{j,1}^* (x,y)\, \tilde k_{j,1} (x,z)\, dx
= a(y)\,a(z)\,\tilde h_j \star h_j (z-y)\,,
\]
with $\tilde h_j (u) = h_j^*(-u)$. It results that the kernel
of $\tilde K = \sum_{j \leq 0} \tilde K_{j,1}^* \tilde K_{j,1}$  is:
\[
\tilde k(y,z) = \sum_{j \leq 0} \tilde k_j (y,z)= a(y)\,a(z)\,\theta(z-y)~~\mbox{with}~~
\theta(x) = \sum_{j \leq 0} \tilde h_j \star h_j (x)~.
\]
Applying Young's inequality to $\|\tilde K f\|$ gives
\[
\|\tilde K\| \leq \sup_{u \in \R^d} |a(u)|^2\, \|\theta \|_1~.
\]
Since $\hat \theta(\omega) = \sum_{j \leq 0} |\hat h(2^{-j} \omega)|^2$ and
$\hat h(0) = \int h(x)\,dx = 0$ and $h$ is both regular with a polynomial decay,
we verify that $\|\theta\|_1 < \infty$.
Moreover, since 
$(1- \det(\Id - \nabla \tau(u))) \geq (1 - \|\nabla \tau\|_\infty)^d$ we have
$\sup_{u} |a(u)|\leq d\, \|\nabla \tau\|_\infty$ which proves
that $\|\tilde K\|^{1/2} \leq C\,\|\nabla \tau\|_\infty$.
Since $\|\tilde K_{j,1}\|^2 \leq \|\tilde K\|$ we get the
same inequality for $\|\tilde K_{j,1}\|^2$, which proves
the two upper bounds of (\ref{upperboundd40dfd8}).

The last sum $\sum_{j=0}^{\infty} K_j^* K_j$ carries the singular part
of the operator, which is isolated and evaluated separately by 
decomposing $K_j = K_{j,1} + K_{j,2}$, with a first kernel
\begin{equation}
\label{epsilonfdnwd013}
k_{j,1} (x,u) = h_j (x-u) - 
h_j ((\Id-\nabla \tau(u)) (x-u)) \,\det(\Id - \nabla \tau(u)) 
\end{equation}
satisfying  $K_{j,1} 1 = \int k_{j,1} (x,u) \,du = 0$ 
if $\int h(x)\,dx = 0$. The second kernel is
\begin{equation}
\label{epsilonfdnwd014}
k_{j,2} (x,u) = \det(\Id - \nabla \tau(u)) \Bigl(
h_j ((\Id-\nabla \tau(u)) (x-u)) - h_j (x - \tau(x) - u+\tau(u)) \Bigr)~.
\end{equation}
The sum $\sum_{j\geq 0} K_{j,1}^* K_{j,1}$ has a singular kernel along
its diagonal, and its norm 
is evaluated separately with the upper bound
\begin{equation}
\label{inequali}
\|\sum_{j=0}^{\infty} K_j^* K_j\|^{1/2} \leq 
\| \sum_{j=0}^{\infty} K_{j,1}^* K_{j,1} \|^{1/2} + 
\sum_{j=0}^{\infty} (\|K_{j,2}\|+ 2^{1/2} \| K_{j,2}\|^{1/2}  \| K_{j,1}\|^{1/2})~.
\end{equation}

We will prove that
\begin{equation}
\label{upperbounddf41}
\|K_{j,1} \| \leq C\, \|\nabla \tau\|_\infty
\end{equation}
and
\begin{equation}
\label{upperboundd8}
\|K_{j,2} \| \leq C\, \min(2^{-j} \|\rH \tau\|_\infty \,,\,\|\nabla \tau\|_\infty)~.
\end{equation}
It implies that
$\|K_{j} \| \leq C\, \|\nabla \tau\|_\infty$. Inserting this
inequality in (\ref{adsfnwonfs0}) yields the first lemma 
result (\ref{firsnfidfnw}) and it proves (\ref{toprove2}).
Equations (\ref{upperbounddf41}) and (\ref{upperboundd8}) also prove that
\begin{equation}
\label{upperboundd400}
\sum_{j=0}^{\infty} (\|K_{j,2}\|+ 2^{1/2} \| K_{j,2}\|^{1/2}  \| K_{j,1}\|^{1/2})
\leq C (\|\nabla \tau \|_\infty + \|\rH \tau \|_\infty )~.
\end{equation}
If $\int h(x)\,dx = 0$ then
thanks to the vanishing integrals of $k_{j,1}$ we will prove that
\begin{equation}
\label{upperboundd40} 
\| \sum_{j=0}^{\infty} K_{j,1}^* K_{j,1} \|^{1/2} \leq 
C\, (\|\nabla \tau \|_\infty + \|\rH \tau \|_\infty )~.
\end{equation}
Inserting (\ref{upperboundd400}) and (\ref{upperboundd40}) 
in (\ref{inequali}) proves (\ref{toprove3}).

Let us now first prove the upper bound (\ref{upperboundd8}) on $K_{j,2}$.
The kernel of $K_{j,2}$ is
\[
k_{j,2} (x,u) = \det(\Id - \nabla \tau(u)) \Bigl(
h_j ((\Id-\nabla \tau(u)) (x-u)) - h_j (x - \tau(x) - u+\tau(u)) \Bigr)~.
\]
A Taylor expansion of $h_j$ together with
a Taylor expansion of $\tau(x)$ gives
\begin{equation}
\label{epsilonfdnwd10}
\tau(x) - \tau(u) = \nabla \tau(u) (x-u) + \alpha(u,x-u)
\end{equation}
with
\begin{equation}
\label{epsilonfdnwd}
\alpha(u,z) = \int_{0}^1 t\,z\,\rH\tau(u+(1-t) z)\,z\,dt~,
\end{equation}
so
\begin{eqnarray}
\label{upperboundsfndd4}
 k_{j,2} (x,u) &=& - \det(\Id - \nabla \tau(u))\,\\
\nonumber
&&
\int_0^1 \nabla h_j 
\Bigl((\Id - t\,\nabla \tau(u))(x-u)+ (1-t)\,(\tau(u) - \tau(x))\Bigr)\, 
\alpha(u,x-u)\,dt~.
\end{eqnarray}

For $j \geq 0$, we prove that $\|K_{j,2} \|$ decays like $2^{-j}$.
Observe that $|\det(\Id - \nabla \tau(u))| \leq 2^d$.
Since $\nabla h_j (u) = 2^{j+dj} \nabla h(2^{j} u)$,
the change of variable $x' = 2^{j} (x-u)$ 
in (\ref{upperboundsfndd4}) gives
\begin{eqnarray*}
\int |k_{j,2} (x,u)|\, dx &\leq& 2^{d}\,\int \left| 
\,\int_0^1 \nabla h
\Bigl((\Id - t\,\nabla \tau(u)) x'\right. \\
& & \left. + (1-t)\,2^{j} (\tau(u) - \tau(2^{-j} x'+u))\Bigr)\, 2^{j}
\alpha(u,2^{-j} x')\,dt\right| \,dx'~.
\end{eqnarray*}
For any $0 \leq t \leq 1$
\[
\Bigl|(\Id - t\,\nabla \tau(u)) x'+ (1-t)\,2^{j} (\tau(2^{-j} x'+u) - \tau(u))\Bigr|  \geq |x'|\, (1 - \|\nabla \tau \|_\infty) \geq |x'|/2~.
\]
Equation (\ref{epsilonfdnwd}) also implies that
\begin{equation}
\label{boundsdfnsl1}
|2^{j}\,\alpha(u,2^{-j}x')| = 2^{-j} 
|\int_{0}^1 t\,x'\,\rH\tau(u+(1-t) 2^{-j} x')\,x'\,dt| \leq 
2^{-j}\, \|\rH \tau \|_\infty\,\frac {|x'|^2} 2~.
\end{equation}
Since $|\nabla h (u)| \leq C \, (1 + |u| )^{-d-2}$, with
the change of variable $x = x'/2$ we get
\begin{equation}
\label{bo8sdfhnsldfs1}
\int |k_{j,2} (x,u)|\, dx \leq C\,2^{-j} \,
\|\rH \tau \|_\infty \, .
\end{equation}

For $j \leq 0$, we use a maximum error bound on the remainder $\alpha$ 
of the Taylor approximation (\ref{epsilonfdnwd10}):
\[
|2^{j}\,\alpha(u,2^{-j}x')| \leq 2\, \|\nabla \tau \|_\infty\, |x'|~,
\]
which proves that $\int |k_{j,2} (x,u)| dx \leq C\,\|\nabla \tau \|_\infty$ and hence that
\begin{equation}
\label{bo8sdfhnsldfs2}
\int |k_{j,2} (x,u)|\, dx \leq C\, \min(2^{-j} \|\rH \tau\|_\infty \,,\,
\|\nabla \tau\|_\infty)~.
\end{equation}

Similarly, we compute $\int |k_{j,2} (x,u)|\, du$ with
the change of variable $u' = 2^{j} (x-u)$ which leads to the same bound
(\ref{bo8sdfhnsldfs2}).
Schur's lemma gives:
\begin{equation}
\label{boundsdfnsld0s1}
\|K_{j,2} \| \leq C\, \min(2^{-j} \|\rH \tau\|_\infty \,,\,\|\nabla \tau\|_\infty)
\end{equation}
which finishes the proof of (\ref{upperboundd8}).

Let us now compute the
upper bound (\ref{upperbounddf41}) on $K_{j,1}$.
Its kernel $k_{j,1}$ in (\ref{epsilonfdnwd013}) can be written 
$k_{j,1} (x,u) = 2^{dj}\, g(u,2^{j} (x-u))$ with
\begin{equation}
\label{boudnsf10}
g(u,v) = h(v) - h((\Id-\nabla \tau(u))v) \,\det(\Id - \nabla \tau(u)).
\end{equation}
A first-order Taylor decomposition of $h$ gives
\begin{eqnarray}
\label{boudnsf101}
g(u,v) &=& (1 - \det(\Id - \nabla \tau(u)))\, 
h((\Id-\nabla \tau(u)v) \\
\nonumber & &+
\int_0^1 \nabla h ((1-t)v + t(\Id - \nabla \tau(u)) v) \cdot \nabla \tau (u) v\, dt. 
\end{eqnarray}
Since
$\det(\Id - \nabla \tau(u)) \geq (1 - \|\nabla \tau\|_\infty)^d$ we get
$(1- \det(\Id - \nabla \tau(u))) \leq d\, \|\nabla \tau\|_\infty$.
Moreover $\|\nabla \tau\|_\infty \leq 1/2$ and $h(x)$ as well as 
its partial derivatives
have a decay which is $O((1 + |x|)^{-d-2})$, so
\begin{equation}
\label{boudnsf1}
|g(u,v)| \leq C\,\|\nabla \tau\|_\infty\, 
\Bigl(1 + |v| \Bigr)^{-d-2}~,
\end{equation}
so $k_{j,1} (x,u) = O\Bigl( 2^{dj}\, \|\nabla \tau\|_\infty\, 
(1 + 2^{j} |x-u| )^{-d-2}\Bigr)$.
Since
\[
\int |k_{j,1} (x,u)|\,du = O( \|\nabla \tau\|_\infty)~~\mbox{and}~~
\int |k_{j,1} (x,u)|\,dx = O( \|\nabla \tau\|_\infty)~,
\]
Schur's lemma (\ref{schs1eq2}) proves that $\|K_{j,1}\| = O( \|\nabla \tau\|_\infty)$
and hence (\ref{upperbounddf41}).

Let us now prove (\ref{upperboundd40}) when $\int h(x)\,dx = 0$.
The kernel of the self-adjoint operator $Q_j = K_{j,1}^* K_{j,1}$  is:
\begin{eqnarray}
\bar k_{j} (y,z) &=& \int k_{j,1}^* (x,y)\, k_{j,1} (x,z)\, dx
= \int 2^{2dj} g^*(y, 2^{j} (x-y))\,g(z, 2^{j} (x-z))\,dx
\nonumber\\
\label{kdnfonfsp}
&=& \int 2^{dj} g^*(y, x'+ 2^{j} (z-y))\,g(z, x'))\,dx'~.
\end{eqnarray}
The singular kernel $\bar k = \sum_j \bar k_{j}$ of $\sum_j Q_j$ 
almost satisfies the hypotheses of the T(1) theorem of 
David, Journ\'e and Semmes \cite{David}
but not quite because it does not
satisfy the decay condition
$|\bar k (y,z) - \bar k (y,z')| 
\leq C |z'-z|^\alpha\, |z-y|^{-d-\alpha}$ 
for some $\alpha > 0$. We bound this operator with
Cotlar's lemma \cite{Stein2} which proves that if $Q_j$ satisfies
\begin{equation}
\label{cotarls}
\forall j,l\,,\,\|Q^*_{j}\,Q_l \| \leq |\beta(j-l)|^2~~\mbox{and}~~
\|Q_{j}\,Q^*_l \| \leq |\beta(j-l)|^2~,
\end{equation}
then
\begin{equation}
\label{cotarls2}
\|\sum_j Q_j \| \leq \sum_{j} \beta(j)~. 
\end{equation}

Since $Q_j$ is self-adjoint, it is sufficient to bound
$\|Q_l\,Q_j\|$.
The kernel of $Q_l\,Q_j$ is computed from the kernel $\bar k_j$ of $Q_j$
\begin{equation}
\label{kdfnsdfjj}
\bar k_{l,j} (y,z) = \int \bar k_{j} (z,u)\,\bar k_{l} (y,u)\,du .
\end{equation}
An upper bound of $\|Q_l\,Q_j\|$ is obtained with
Schur's lemma (\ref{schs1eq2}) applied to $\bar k_{l,j}$.
Inserting (\ref{kdnfonfsp}) in (\ref{kdfnsdfjj}) gives
\begin{eqnarray}
\int |\bar k_{l,j} (y,z)|\,dy = 
\int \left| \right.\int && g(u,x)\,g(u,x')\,
2^{dl}\,
g^*(y,x+2^{l} (u-y))\,\nonumber\\
\label{boudnsf0}
& & \left.
2^{dj}\,g^*(z,x'+2^{j} (u-z))\,dx\,dx'\,du \right| \, dy~.
\end{eqnarray}
The parameters $j$ and $l$ have symmetrical roles and we can
thus suppose that $j\geq l$.

Since $\int h(x)\,dx = 0$ it results from (\ref{boudnsf10}) that
$\int g(u,v)\,dv = 0$ for all $u$. 
For $v = (v_n)_{n \leq d}$, one can thus write 
$g(u,v) = \frac {\partial \bar g (u,v)} {\partial v_1}$
and (\ref{boudnsf1}) implies that
\begin{equation}
\label{boudnsf01}
|\bar g(u,v)| \leq C\,\|\nabla \tau \|_\infty\,
\Bigl(1 + |v| \Bigr)^{-d-1}~.
\end{equation}
Let us make an integration by parts along the variable $u_1$ in
(\ref{boudnsf0}). 
Since all first and second-order derivatives of $h(x)$ have
a decay which is $O((1 + |x|)^{-d-2})$, we derive 
from (\ref{boudnsf10}) that for any $u = (u_n)_{n \leq d} \in \R^d$ and
$v = (v_n)_{n \leq d} \in \R^d$
\begin{equation}
\label{boudnsf2}
\left|
\frac{\partial g(u,v)} {\partial u_1}
\right| \leq C\,\|\rH \tau \|_\infty\,
\Bigl(1 + |v|\,(1 - \|\nabla \tau\|_\infty) \Bigr)^{-d-1}~,
\end{equation}
and from (\ref{boudnsf101}) 
\begin{equation}
\label{boudnsf3}
\left|
\frac{\partial g(u,v)} {\partial v_1}
\right| \leq C\,\|\nabla \tau\|_\infty\, 
\Bigl(1 + |v|\,(1 - \|\nabla \tau\|_\infty) \Bigr)^{-d-1}~.
\end{equation}
In the integration by part, 
integrating $2^{dj} g(z,x'+2^{j} (u-z))$ 
brings out a term proportional to $2^{-j}$ 
and differentiating
$g(u,x)\,g(u,x')\,2^{dl} g(y,x+2^{l} (u-y))$
brings out a term bounded by $2^{l}$. 
An upper bound of (\ref{boudnsf0}) 
is obtained by inserting (\ref{boudnsf1},\ref{boudnsf01},
\ref{boudnsf2},\ref{boudnsf3}), 
which prove that there exists $C$ such that
\begin{eqnarray*}
\int |\bar k_{l,j} (y,z)|\,dy &\leq& C^2 \,(2^{-j}\, \|\nabla \tau\|^3_\infty\, 
\|\rH \tau \|_\infty + 2^{l-j}\, \|\nabla \tau\|^4_\infty) \\
&\leq&
C^2\,2^{l-j}\, (\|\nabla \tau\|_\infty + \|\rH \tau \|_\infty )^4~.
\end{eqnarray*}
The same calculation proves the same bound on $\int |\bar k_{l,j} (y,z)|\,dz$
so Schur's lemma (\ref{schs1eq2}) implies that 
\[
\|Q_l\,Q_j \| \leq C^2\,2^{l-j}\, (\|\nabla \tau\|_\infty + \|\rH \tau \|_\infty )^4~.
\]
Applying Cotlar's lemma (\ref{cotarls}) with 
$\beta(j) = C\,2^{-|j|/2}\,(\|\nabla \tau\|_\infty + \|\rH \tau \|_\infty )^2$ proves that 
\begin{equation}
\label{upperboundd4}
\| \sum_{j=-\infty}^{+\infty} K_{j,1}^* K_{j,1} \| 
= \| \sum_j Q_j \| \leq 
C\,(\|\nabla \tau\|_\infty + \|\rH \tau \|_\infty )^2~,
\end{equation} 
which implies (\ref{upperboundd40}).

\section{Proof of Lemma \protect\ref{Cauchsdf}}
\label{Cauchsdfproof}

It results from (\ref{convadfndsoihfdsf}) that there exists
$\epsilon_J$ with $\lim_{J \rightarrow \infty} \epsilon_J = 0$ such that
\[
\sup_{p \in \cP_J -  \Omega_J^f} \Big\| {S_J [p] f}  -
\frac {\|S_J [p] f\|} 
{\|S_J [p] \delta\|}  {S_J [p] \delta}  \Big\|^2 \leq \frac{\epsilon_J} 2 \|S_J[p] f\|^2\,,
\]
and $\sum_{p \in  \Omega_J^f} \|{S_J [p] f}\|^2 \leq {\epsilon_J}  \|f\|^2/8$.
Since $\|S_J[\cP_J] f\|^2= \|f\|^2$, we get
\begin{equation}
\label{fid0f98sdf1}
\sum_{p \in \cP_J} \Big\| {S_J [p] f}  -
\frac {\|S_J [p] f\|} 
{\|S_J [p] \delta\|}  {S_J [p] \delta}  \Big\|^2 \leq \epsilon_J \| f\|^2.
\end{equation}

The set of all extensions of a $p \in \cP_J$ into
$\cP_{J\pl 1}$ is defined in (\ref{extenset}). It
can be rewritten $\cP_{J\pl 1}^{p} = \cP_{J\pl 1} \cap C_J(p)$,
and (\ref{pathconsdf9u8sd}) proves that
\[
\|S_J[p] f - S_J[p] h\|^2 \geq 
\sum_{p' \in \cP_{J\pl 1} \cap C_J(p)} 
\|S_{J\pl 1}[p'] f - S_{J\pl 1}[p'] h\|^2 .
\]
Iterating $k$ times on this result yields
\[
\|S_J[p] f - S_J [p] h\|^2 \geq 
\sum_{p' \in \cP_{J+k} \cap C_{J+k}(p)} 
\|S_{J+k}[p'] f - S_{J+k}[p'] h\|^2 .
\]
Applying it to
$f$ and $h = \mu_p \delta$
with $\mu_p = {\|S_J[p] f\|} /{\|S_J[p] \delta\|}$
gives
\[
\Big\|S_J[p] f - \mu_p \,S_J[p] \delta  \Big\|^2 \geq 
\sum_{p' \in \cP_{J+k} \cap C_{J+k}(p)} 
\Big\|S_{J+k} [p'] f - 
\mu_p\,S_{J+k}[p'] \delta  \Big\|^2 .
\]
Summing over $p \in \cP_J$ and
applying (\ref{fid0f98sdf1}) proves that
\[
\sum_{p \in \cP_J} \sum_{p' \in \cP_{J+k} \cap C_{J+k}(p)} 
 \Big\| {S_{J+k}[p'] f}  -
\frac {\|S_J[p] f\|} 
{\|S_J[p] \delta\|}  {S_{J+k}[p'] \delta}  \Big\|^2 \leq \epsilon_J \| f\|^2~,
\]
and hence
\[
\sum_{p \in \cP_J} \sum_{p' \in \cP_{J+k} \cap C_{J+k}(p)} 
\Big| \frac {\|S_{J+k}[p'] f\|}  {\|S_{J+k}[p'] \delta \|}  -
\frac {\|S_J[p] f\|} {\|S_J[p] \delta\|}  
\Big|^2 \,\|S_{J+k}[p'] \delta \|^2\leq \epsilon_J \| f\|^2~.
\]
If $q \in C_{J+k}(p')$
then 
$S_{J+k}(q) ={\|S_{J+k}[p'] f\|}/ {\|S_{J+k}[p'] \delta\|}$.
But $p'\in C_J(p)$ so
$q \in C_J (p)$ and hence
$S_{J}(q) ={\|S _J [p] f\|}/ {\|S _J [p] \delta\|}$.
Finally $\|S_{J+k}[p'] \delta \|^2 = \mu(C_{J+k}(p'))$
so the sum can be rewritten as a path integral
\[
\int_{\cP^{\infty}} | S_{J+k}f(q) -  S_J  f(q)|^2\, d\mu(q) \leq 
\epsilon_J \| f\|^2~,
\]
which proves that $\{S_J f \}_{J\in \N}$ is a Cauchy sequence in $\Ld(\cPinf,d\mu)$.

\section{Proof of Lemma \protect\ref{lemmarandom}}
\label{prooflemmarandom}

This appendix proves that
\begin{equation}
\label{Kautansinf00}
E(|K_\tau X |^2) \leq E(\|K_\tau \|^2)\, E(|X|^2)~,
\end{equation}
as well as a generalization
to sequence of operators, at the end of the appendix. 
The lemma result is proved by restricting $X$ to a finite hypercube
$I_T = \{ (x_1,...,x_d) \in \R^d~:~\forall i \leq d~,~|x_i| \leq T \}$,
whose indicator function $\One_{I_T}$ 
defines a finite energy process $X_T (x) = X(x)\, \One_{I_T} (x)$. 
We shall verify that 
$E(|K_\tau X(x)|^2)$ does not depend upon $x$ and that
\begin{equation}
\label{Kautansinf}
E(|K_\tau X(x)|^2) = 
\lim_{T \rightarrow \infty} \frac { E(\|K_\tau X_T \|^2)}
{(2\,T)^d} ~.
\end{equation}
Let first show how this result implies
(\ref{Kautansinf00}). The $\LD$ operator norm definition implies
\[
\|K_\tau X_T \|^2 =  \int
|K_\tau X_T (x)|^2\,dx \leq \|K_\tau \|^2\, 
\int |X_T (x)|^2\,dx ~.
\]
Since $X$ and $\tau $ are independent processes
\[
E( \|K_\tau X_T \|^2) \leq E(\|K_\tau \|^2)\, 
 E(|X|^2)\,(2 T)^d~.
\]
Applying (\ref{Kautansinf}) thus proves 
the lemma result (\ref{Kautansinf00}).

To prove (\ref{Kautansinf}), we first compute
\[
E(|K_\tau X(x)|^2) = E \Bigl( \iint k_\tau (x,u)\,k^*_\tau (x,u')\,X(u)\, X^*(u') du\,du' \Bigr)~.
\]
Since $X$ is stationary  $E(X(u)\,X^*(u')) = A_X (u-u')$, and
the lemma hypothesis supposes that 
$E(k_\tau (x,u)\,k^*_\tau (x,u') )= \bar k_\tau (x-u,x-u')$.
Since $X$ and $\tau $ are independent, the change of variable
$v = x-u$ and $v' = x-u'$ gives
\begin{eqnarray}
E(|K_\tau X(x)|^2) &=& 
\iint \bar k_\tau(x-u,x-u') \,A_X (u-u') \, du\,du' \nonumber \\
\label{Kautansinf2}
&=&  \iint \bar k_\tau(v,v') \,A_X (v-v') \, dv\,dv'~,
\end{eqnarray}
which proves that $E(|K_\tau X(x)|^2)$ does not depend upon $x$.
Similarly
\begin{equation}
\label{Kautansinf3}
E(|K_\tau X_T(x)|^2) 
=  \iint \bar k_\tau(v,v') \,A_X (v-v') \, 
\One_{I_T}(v-x)\,\One_{I_T}(v'-x)\,dv\,dv',
\end{equation}
and integrating along $x$ gives
\begin{equation}
\label{Kautansinf30}
(2T)^{-d}\,E(\|K_\tau X_T\|^2) 
=  \iint \bar k_\tau(v,v') \,A_X (v-v') \, (1 - \rho_T (v-v'))\,dv\,dv',
\end{equation}
with
\[
1 - \rho_T (v-v') = (2T)^{-d}
\int \One_{I_T}(v-x)\,\One_{I_T}(v'-x) \,dx = \prod_{i=1}^d 
\Bigl(1 - \frac{|v_i-v'_i|} {2T} \Bigr)\, \One_{I_T}(v-v') ~
\]
and hence
\begin{equation}
\label{Kautansinf40}
0 \leq \rho_T(v) \leq (2T)^{-1}\sum_{i=1}^d |v_i| \leq d\,(2T)^{-1} |v|~.
\end{equation}
Inserting (\ref{Kautansinf2}) in (\ref{Kautansinf30}) proves that
\begin{equation}
\label{Kautansinf4}
(2\,T)^{-d} E(\|K_\tau X_T\|^2) 
=  E(|K_\tau X(x)|^2) -  \iint \bar k_\tau(v,v') \,A_X (v-v') \, \rho_T(v-v')\,dv\,dv'~.
\end{equation}
Since $\iint |\bar k_\tau (v,v')| \,|v-v'| \,dv\,dv' < \infty$
and $A_X(v-v') \leq A_X(0) = E(|X|^2)$, it results from
(\ref{Kautansinf4}) and (\ref{Kautansinf40}) that 
\[
\lim_{T \rightarrow \infty} (2\,T)^{-d} E(\|K_\tau X_T\|^2)   =
E(|K_\tau X(x)|^2) ~,
\]
which proves (\ref{Kautansinf}).

Lemma \ref{lemmarandom} is extended to 
sequences of operators 
$\overline K_\tau = \{ K_{\tau,n} \}_{n\in I}$ with kernels 
$\{k_{\tau,n} \}_{n \in I}$, as follow. Let us denote
\begin{equation}
\label{operandfonsdf}
\|\overline K_\tau X \|^2 = \sum_{n\in I} |K_{\tau,n} X|^2 ~~\mbox{and}~~
\|\overline K_\tau f \|^2 = \sum_{n\in I} \|K_{\tau,n} f\|^2 ~.
\end{equation}
If each average bilinear kernel is stationary
\begin{equation}
\label{Kautansinf000hp}
E(k_{\tau,n} (x,u)\,k^*_{\tau,n} (x,u') ) = \bar k_{\tau,n} (x-u,x-u')~~
\end{equation}
and
\begin{equation}
\label{Kautansinf000hp2}
\iint |\sum_{n\in I} \bar k_{\tau,n} (v,v')| \,|v-v'|\,dv\,dv' < \infty~,
\end{equation}
then 
\begin{equation}
\label{Kautansinf000}
E(\|\overline K_\tau X \|^2) \leq E(\|\overline K_\tau \|^2)\, E(|X|^2)~.
\end{equation}
The proof of this extension follows the same derivations as the
proof of (\ref{Kautansinf00}) for a single operator. 
It just requires to replace the $\LD$ norm
$\|f\|^2$ by the norm $\sum_{n\in I} \|f_n\|^2$ over the space
of finite energy sequences $\{f_n \}_{n \in I}$ of $\LD$ functions
and the sup operator norms in $\LD$ by sup operator
norms on sequence of $\LD$ functions.

\section{Proof of Theorem \ref{Theoelas6}}
\label{theornadofns} 
This appendix proves that 
$E(\|[S_J[\cP_J]\,,\, L_\tau] X\|^2) \leq E(\|U[\cP_J] X \|_1)^2
\,B(\tau)$
with
\begin{equation}
\label{commnaf872}
B(\tau) = 
C\,E \Bigl\{ \Bigl(\|\nabla \tau \|_\infty (\log \frac{\|\Delta \tau\|_\infty}{\|\nabla \tau \|_\infty }
\vee 1)+ \|\rH \tau \|_\infty \Bigr)^2\Bigr\} ,
\end{equation}
and $E(\|U[\cP_J] X \|_1) = \sum_{m=0}^{+\infty}
\Bigl(\sum_{p \in \cLa^m_J} E( |U [p] X |^2)  \Bigr)^{1/2}$.

For this purpose, we shall first prove that if
for any stationary process $X$
\begin{equation}
\label{Updfnd09n}
E(\|[W_J , L_\tau] X\|^2) \leq B(\tau)\,E(|X|^2)
\end{equation}
where
\[
E( \|[W_J\,,\,L_\tau] X\|^2) = 
E(|[A_J\,,\,L_\tau] X|^2) + \sum_{\lambda \in \Lambda_J} 
E(|[W[\lambda] \,,\,L_\tau] X|^2)
\]
then
\begin{equation}
\label{comutarand0}
E(\|[S_J[\cP_J]\,,\, L_\tau] X\|^2)
 \leq B(\tau) \, E(\|U[\cP_J] X \|_1 )^2~.
\end{equation}

Since a modulus operator is nonexpansive and commutes with $L_\tau$,
with the same argument as in the proof of (\ref{err10dfi79}), 
we derive from (\ref{Updfnd09n}) that
\begin{equation}
\label{comutarand089hf8}
E(\|[ U_J , L_\tau] X\|^2) \leq B(\tau)\,E(|X|^2)~.
\end{equation}
The proof of Proposition \ref{confppfoidnsf} also shows that
$U_J$ is nonexpansive for the mean square norm
on processes. Since $S_J[\cP_J]$ is obtained by iterating on $U_J$ 
it results that
\[
E(\|[S_J[\cP_J ]\,,\, L_\tau] X\|^2)  \leq B(\tau)
\, E(\|U[\cP_J] X \|_1 )^2~.
\]
The proof of this inequality 
follows the same derivations as in Appendix \ref{Proofsdfn},
for $L = L_\tau$, by replacing
$f$ by $X$, $\|f\|^2$ by $E(|X|^2)$, 
$\|U[p] f\|^2$ by $E(|U[p] X|^2)$, 
and the $\LD$ sup operator norm
$\|[{U_J},L]\|$ by $B(\tau)$ which satisfies (\ref{comutarand089hf8})
for all $X$.

The proof of (\ref{commnaf872}) is ended by verifying that
\begin{equation}
\label{Updfnd09n2}
E(\|[W_J , L_\tau] X\|^2) \leq E(C^2(\tau))\,E(|X|^2)
\end{equation}
and hence $B(\tau) = E(C^2(\tau))$ with 
\[
C(\tau) =C\,
\Bigl(\|\nabla \tau \|_\infty (\log \frac{\|\Delta \tau\|_\infty}{\|\nabla \tau \|_\infty }
\vee 1) + \|\rH \tau \|_\infty\Bigr)~.
\]
The inequality (\ref{Updfnd09n2}) is derived from 
Lemma \ref{thedw0} which proves that 
the $\LD$ operator norm of the commutator $[W_J , L_\tau]$ satisfies
\begin{equation}
\label{commnaf7}
\|[W_J , L_\tau] \| 
\leq C(\tau)~,
\end{equation}
and by applying to
$\overline K_\tau = [W_J , L_\tau]  = \{[A_J,L_\tau]\,,\,[W[\la]\,,\,L_\tau] \}_{\la \in \Lambda_J}$ the extension 
(\ref{Kautansinf000}) of Lemma \ref{lemmarandom}.
This extension proves that if 
the kernels of the wavelet
commutator satisfy the conditions (\ref{Kautansinf000hp})
and (\ref{Kautansinf000hp2}) then
\[
E(\|[W_J , L_\tau] X\|^2) \leq E(\|[W  , L_\tau] \|^2)
\,E(|X|^2).
\]
Together with (\ref{commnaf7}) it proves (\ref{Updfnd09n2}).

To finish the proof we 
verify that the wavelet commutator kernels satisfy (\ref{Kautansinf000hp})
and (\ref{Kautansinf000hp2}).
If $Z_j f (x) = f \star h_j (x)$ with $h_j (x) = 2^{dj} h(2^{j} x)$
then the kernel of the integral commutator operator
$[Z_j , L_\tau] = Z_j L_\tau - L_\tau Z_j$ is
\begin{eqnarray}
\label{kerneldsfn}
 k_{\tau,j} (x,u)&=&  h_j (x-u-\tau(x))\\
 \nonumber
&& - h_j (x - u - \tau(u+\tau(\beta(u))))\,
|\det (\Id - \nabla \tau(u+\tau(\beta(u))))|^{-1}~
\end{eqnarray}
where $\beta$ is defined by $\beta(x) = x + \tau(\beta(x))$.
The kernel of $[A_J ,L_\tau]$ is $k_{\tau,J}$ with $h = \phi$, and the kernel of
$[W[\la],L_\tau]$ for $\lambda = 2^j\ga$ is $k_{\tau,j}$ with $h(x) = 
\psi(r^{-1} x)$.
Since $\tau$ and $\nabla \tau$ are jointly stationary, the joint probability
distribution of their
values at $x$ and $u + \tau(\beta(u))$ only depends upon $x-u$.
It results that 
$E(k_{\tau,j} (x,u)\,k_{\tau,j} (x,u') ) = \bar k_{\tau,j} (x-u,x-u')$ 
which proves the kernel stationarity (\ref{Kautansinf000hp}) for wavelet
commutators.

The second 
kernel hypothesis (\ref{Kautansinf000hp2}) is proved by showing that if 
$|h(x)| = O((1 + |x|)^{-d-2})$ then
\[
\iint |\sum_{j \geq -J} \bar k_{\tau,j} (v,v')| \,|v-v'|\,dv\,dv' < \infty~.
\]
Since $\bar k_{\tau,j} (v,v')= E(k_{\tau,j} (x,x-v)\,k_{\tau,j} (x,x-v'))$, it
is sufficient to prove that there exists $C$ such that for all $x$,
with probability $1$
\begin{equation}
\label{kerneldsfnf8s0}
I = \sum_{j \geq -J} \iint |k_{\tau,j} (x,x-v)|\,|k_{\tau,j} (x,x-v')|\,|v-v'|\,
dv\,dv' \leq C~.
\end{equation}
Since $h_j (x) = 2^{dj} h(2^{j} x)$ and $u + \tau(\beta(u)) = \beta(u)$, 
it results from (\ref{kerneldsfn}) that
$k_{\tau,j} (x,x-2^{-j} w)= 2^{dj}\,\tilde k_{\tau,j} (x,x-w)$ with
\begin{eqnarray}
\label{kerneldsfnf8s}
\tilde k_{\tau,j} (x,x-w) &=&
h (w-2^{j} \tau(x)) \\
\nonumber 
& & - h (w - 2^{j} \tau(\beta(x-2^{-j} w)))\,
|\det (\Id - \nabla \tau(\beta(x-2^{-j} w)))|^{-1}~.
\end{eqnarray}
The change of variable $w = 2^{j} v$ and $w' = 2^{-j} v'$ in 
(\ref{kerneldsfnf8s0}) shows that 
$I = \sum_{j\geq -J} 2^{-j} I_j$ with
\[
I_j =
\iint |\tilde k_{\tau,j} (x,x- w)|\,|\tilde k_{\tau,j} (x,x-w')|\,|w-w'|\,
dw\,dw' ~.
\]
Since $|h(w)| = O((1 + |w|)^{-d-2})$ and $\|\nabla \tau\|_\infty \leq 1/2$
with probability $1$,
by computing separately the integrals of each 
of the four terms of the
product $|\tilde k_{\tau,j} (x,x+ w)|\,|k_{\tau,j} (x,x+w')|\,|w-w'|$,
with change of variables, $y = w+2^{j} \tau(x)$ and
$z = w+ 2^{j} \tau(\beta(x+2^{-j}w))$, we verify that there exists $C'$ such
that $I_j \leq C'$ and hence that $I = \sum_{j\geq -J} 2^{-j} I_j \leq 2^{J+1} C'$
with probability 1. It proves (\ref{kerneldsfnf8s0}) and
hence the second kernel hypothesis  (\ref{Kautansinf000hp2}).\\

{\bf Acknowledgement} I would like to thank Joan Bruna, Mike Glinsky and Nir Soren for the many inspiring conversations in connection with image processing, physics and group theory.

\frenchspacing
\bibliographystyle{plain}

\begin{thebibliography}{99}
\bibitem{Joakim}
J. Anden, and S. Mallat, ``Multiscale scattering for audio classification'',
Proc. of ISMIR Conf., Florida, 2011.

\bibitem{Poggio}
J. Bouvrie, L. Rosasco, T. Poggio: ``On Invariance in Hierarchical Models". Proc. of NIPS 2009

\bibitem{Bruna}
J. Bruna, S. Mallat, ``Classification with scattering operators'', 
Proc. of IEEE CVPR, Colorado Springs, 2011.

\bibitem{Bruna2}
J. Bruna, S. Mallat, ``Invariant Scattering Convolution Networks'',
submitted to IEEE Trans. on PAMI, 2012.


\bibitem{David}
G. David, J.-L. Journ\'e, and  S. Semmes, ``Op\'erateurs de Calder\`on Zygumund, fonctions para-accr\'etives et interpolation,'' Rev. Mat. Iberoameri. 1:1-56, 1985.

\bibitem{frankel}
T. Frankel, ``The Geometry of Physics, An Introduction'', Cabridge Univ. Press, 2004.

\bibitem{Frazier}
M. Frazier, B. Jawerth, G. Weiss, ``Littlewood-Paley theory and the 
study of function spaces,'' CBMS ó AMS, (1991).


\bibitem{Geller1}
D. Geller and I. Z. Pesenson,
``Band-Limited Localized Parseval Frames and Besov Spaces on Compact Homogeneous Manifolds'', J. of Geometric Analysis, July, 2010.

\bibitem{Glinsky}
M. Glinsky, ``A new perspective on renormalization: invariant actions, a dynamical DNA'', June 2011, http://arxiv.org/abs/1106.4369

\bibitem{Hormander}
Lars H\"ormander, 
``Fourier Integral Operators, I,''Acta Mathematica, 127:79ó183, 1971.


\bibitem{LeCun1}
Y. LeCun, K. Kavukvuoglu and C. Farabet: ``Convolutional Networks and Applications in Vision", Proc. Int. Symposium on Circuits and Systems (ISCAS'10)'', IEEE, 2010.


\bibitem{slotine}
W. Lohmiller and J.J.E. Slotine ``On Contraction Analysis for Nonlinear Systems'', Automatica, 34(6), 1998.


\bibitem{MallatEUSIPCO}
S. Mallat.  ``Recursive Interferometric Representation'',
Proc. of EUSICO conference, Danemark, August 2010.

\bibitem{Meyerbook}
Y. Meyer, ``Wavelets and Operators'', Cambridge University Press, 1992.

\bibitem{Olver}
P. Olver, ``Equivalence, Invariants and Symmetry'' Cambridge University Press, 1995.

\bibitem{Donoho2}
I. Rahman, I. Drori ,  V. Stodden ,  D. Donoho,
"Multiscale Representations for Manifold-Valued Data," SIAM Multiscale Modelling and Simulation, vol. 4, no. 4, 2005, pp. 1201ó1232.


\bibitem{poggio2}
M. Riesenhuber, T. Poggio, ``Hierarchical models of object recognition 
in cortex,'' Nature Neuroscience, 2: 1019ó1025.

\bibitem{Sifre}
L. Sifre, S. Mallat, ``Combined scattering for rotation invariant texture analysis,''
subm. ESANN 2012.


\bibitem{Stein1}
E. Stein, ``Topics in harmonic analysis related to Littlewood-Paley theory'', 
Annals of  Mathematical Studies, Princeton University Press, 1970.

\bibitem{Stein2}
E. Stein, ``Harmonic Analysis'', 
Princeton Mathematical Series, Princeton University Press, 1993.

\bibitem{Trouve}
A. Trouv\'e, L. Younes; ``Local Geometry of Deformable Templates''; SIAM J. Math Anal. Volume 37, Issue 1, pp. 17-59; 2005.

\bibitem{Willard}
Stephen Willard, ``General Topology'' Addison-Wesley Publishing Company, Reading Massachusetts, 1970.




\end{thebibliography}

\end{document}